\magnification=\magstep1
\baselineskip =13pt
\overfullrule =0pt

\centerline {\bf HEISENBERG DOUBLES AND DERIVED CATEGORIES}

\vskip .7cm

\centerline {\bf M. Kapranov}

\vskip 1cm
Let $\cal A$ be an Abelian category of finite homological dimension in which all ${\rm Ext}_{\cal A}^i(A,B)$ are finite sets.
One can, following C.M. Ringel [R1-3], associate to $\cal A$ an algebra $R({\cal A})$,
a version of the Hall algebra construction. Its structure constants are suitably normalized
numbers of short exact sequences.
Ringel has shown that in the case ${\cal A} = {\rm Rep}_{F_q}(\Gamma)$, the category of ${\bf F}_q$-representations
of a Dynkin quiver $\Gamma$, the algebra $R({\cal A})$ is identified with  the ``nilpotent" subalgebra $U_q({\bf n}^+)$
in the $q$-quantization $U_q({\bf g})$ of the semisimple Lie algebra corresponding to $\Gamma$.
This discovery has lead to several substantial advances in quantum group theory [Lu 1]. 
He has also shown how to put an algebra structure on the space $B({\cal A})
 = {\bf C}[{\cal K}_0 {\cal A}] \otimes R({\cal A})$
(here ${\cal K}_0 {\cal A}$ is the Grothendieck group of $\cal A$) so that for
${\cal A} = {\rm Rep}_{F_q}(\Gamma)$
one has $B({\cal A}) = U_q({\bf b}^+)$, the ``Borel" part of $U_q({\bf g})$.

The idea of extending the Hall algebra formalism to triangulated categories such as the
derived category $D^b({\cal A})$, seems to have been voiced by several people independently.
It appears naturally if one tries to find a construction of the full quantum group
$U_q({\bf g})$ in terms of ${\cal A} = {\rm Rep}_{F_q}(\Gamma)$. Indeed, various
 Abelian subcategories
in $D^b({\cal A})$ obtained from $\cal A$ by repeated application of derived
 Bernstein-Gelfand-Gelfand
reflection functors [GM] (and the Hall algebras of these subcategories) look temptingly
 similar to Borel subalgebras ${\bf b}^w\i {\bf g}$ obtained from ${\bf b}^+$ by action of
 elements of the Weyl group. 
Unfortunately, a direct mimicking of the Hall algebra construction but with exact
triangles replacing exact sequences, fails to give an associative multiplication,
even though the octohedral axiom looks like the right tool to establish the associativity.
One way to get around this difficulty is, as it was done in [X1], to ``amalgamate"
the (associative) Hall algebras of various Abelian subcategories in $D^b({\cal A})$,
but it seems to be not clear whether the resulting algebra is indeed $U_q({\bf g})$ nor
that it is a Hopf algebra at all. 

\vskip .1cm

The aim of the present paper is to exhibit an algebra $L({\cal A})$ which, although defined
in terms of $\cal A$, is invariant under derived equivalences and can be thus called
the ``Hall algebra of the derived category". We call $L({\cal A})$ the lattice algebra
of $\cal A$. Its construction was suggested by the fact that $U_q({\bf g})$ can be obtained
from the Hopf algebra  $U_q({\bf b}^+)$ by the Drinfeld double construction,
while the Hopf algebra structure on $B({\cal A})$ can be described in purely categorical terms,
as follows from the recent work of Green [Gr] (see [X2], [Kap]). However,
it turned out that it is not the Drinfeld double which appears naturally
in the study of $D^b({\cal A})$, but rather the so-called Heisenberg double of [AF] [ST].
In fact, one can find counterparts in Hopf algebra theory of several different
versions of derived categories, as shown in the following table:
\vfill\eject

\halign{\quad #  \hfill &\vrule width .8pt\vbox to
 .6cm{}\quad#{}\hfill\cr
\hskip .5cm{\bf Categories related to an }& \hskip .5cm {\bf Algebras related to the}\cr
\hskip .5cm {\bf Abelian category $\cal A$}& \hskip .5cm  {\bf Hopf algebra $\Xi = B({\cal A})$}\cr
\noalign{\hrule height .8pt}
$D^{[-1,0]}({\cal A})$, the category of {\bf Z}-graded& $HD(\Xi)$, the Heisenberg double\cr
 complexes situated in degrees $0, -1$& \cr
\noalign{\hrule}
 $D^b({\cal A})$, the standard bounded derived & \hskip .5cm $L({\cal A})$, the lattice algebra
\cr \hskip 2cm category &
\cr
\noalign{\hrule}
$D^{(2)}({\cal A})$, the category of 2-periodic & $DD(\Xi)$, the Drinfeld double
\cr \hskip 1cm (${\bf Z}/2$-graded) complexes & \cr
\noalign{\hrule}
 }

\vskip .3cm

The relation of $D^{[-1,0]}({\cal A})$ to the Heisenberg double of $\Xi=B({\cal A})$
is the easiest to understand: the commutation relations in $HD(\Xi)$ involve
certain products of the structure constants for the multiplication and comultiplication
in $\Xi$, and they can be interpreted as numbers of some 4-term exact sequences in $\cal A$,
which are obviously related to exact triangles in $D^{[-1,0]}({\cal A})$.

The algebra $L({\cal A})$ is obtained by taking one copy of
${\bf C}[{\cal K}_0({\cal A})]$ and infinitely many copies of $R({\cal A})$, one for each site of
an infinite 1-dimensional lattice and then imposing Heisenberg double-type commutation
relations between copies of $R({\cal A})$ at adjacent sites and oscillator relations of
the form $AB=\lambda_{AB} BA$, $\lambda_{AB}\in {\bf R}$, between basis vectors
of non-adjacent copies.  This algebra is similar to the ``lattice Kac-Moody algebras"
of [AFS]. The reason for taking an infinite lattice is clear: the $n$th copy of $R({\cal A})$,
$n\in {\bf Z}$, corresponds to the Abelian subcategory ${\cal A}[n]\i D^b({\cal A})$. 

\vskip .1cm

The author would like to acknowledge financial support from NSF grants and A.P. Sloan
Research Fellowship as well as from the Max-Planck Institute f\"ur Mathematik in Bonn
which provided excellent conditions for working on this paper.

\vfill\eject

\centerline {\bf \S 1. Heisenberg doubles for Hall-Ringel algebras}

\vskip 1cm

\noindent {\bf (1.1) The Heisenberg double.}
Let $\Xi$ be a Hopf algebra over {\bf C}.
We denote by 
$$\Delta: \Xi\rightarrow \Xi\otimes\Xi,  \quad
\epsilon: \Xi\rightarrow {\bf C}, \quad S: \Xi\rightarrow
\Xi$$
the comultiplication, the counit and the antipode of $\Xi$ respectively. 
For $x\in\Xi$ let $r_x: \Xi\rightarrow\Xi$ be the operator of right  multiplication by $x$, i.e., $r_x(y) = yx$, and let
$D_x: \Xi^*\rightarrow \Xi^*$ be the dual to $r_x$, i.e., for $f\in\Xi^*$ the functional $D_x(f)\in\Xi^*$ takes $y\mapsto f(yx)$.
Then the correspondence $x\mapsto D_x$ gives an embedding of algebras $D: \Xi\rightarrow {\rm End}(\Xi^*)$.

\vskip .1cm

For $f\in\Xi^*$ let $l_f: \Xi^*\rightarrow\Xi^*$ be the operator of left multiplication by $f$ (with respect to the algebra
structure on $\Xi^*$ defined by the map dual to $\Delta$), i.e., $l_f(\phi) = f\phi$. Again, we get an embedding of algebras
$l: \Xi^*\rightarrow {\rm End}(\Xi^*)$. 

\vskip .1cm

The Heisenberg double $HD(\Xi)$ is defined as the subalgebra in ${\rm End}(\Xi^*)$ generated by the images $D(\Xi)$ and $l(\Xi^*)$.
It is known [ST] that the map
$$\Xi^*\otimes \Xi\rightarrow HD(\Xi), \quad f\otimes x\mapsto l_fD_x,$$
is an isomorphism of vector spaces. Thus to describe the structure of $HD(\Xi)$ completely, it is enough to
explain how to bring a product $D_xl_f$ to a linear combination of products of the form $l_gD_y$. We will do this in the
coordinate-dependent language, following [Kas]. 

\vskip .2cm

Let $\{e_i\}, i\in I$, be a basis of $\Xi$ and $\{e^i\}$ be the dual
(topological) basis of $\Xi^*$. Introduce the structure constants for
the  multiplication and comultiplication with respect to our basis:
$$ e_{i} e_{j} = \sum_k m_{i j }^k e_k, \quad \Delta (e_k) = 
\sum \mu_{k}^{ij} e_{i}\otimes 
 e_{j}. \leqno (1.1.1)$$
Then one easily finds that
$$D_{e_i}(e^k) = \sum_j m_{ji}^k e^j, \quad l_{e^i}(e^k) = \sum_j \mu_j^{ik} e^j. \leqno (1.1.2)$$

\vskip .1cm

\proclaim (1.1.3) Proposition. In $HD(\Xi)$ we have the identity
$$D_{e_i} l_{e^j} = \sum_{a,b,c} m_{ab}^j \mu_i^{bc} l_{e^a} D_{e_c}.$$

\noindent {\sl Proof:} By (1.1.2), our statement is equivalent to:
$$\sum_k \mu_k^{jr} m_{si}^k  = \sum_{a,b,c,d} m_{ab}^j \mu_i^{bc} m_{dc}^r \mu_s^{ad}, \quad \forall i,j,r,s.$$
This equality, however, expresses the coincidence of the coefficients at $e_j\otimes e_r$ in $\Delta(e_se_i)$ and
$\Delta(e_s)\Delta(e_i)$, and so it is true.

\vskip .2cm

In view of this proposition we will regard $HD(\Xi)$ as an abstract algebra
 generated by symbols $Z_i, Z^i, i\in I$ subject to the
relations:
$$Z_iZ_j = \sum_k m_{ij}^k Z_k, \quad Z^iZ^j = \sum_k \mu_k^{ij} Z^k, \leqno (1.1.4)$$
$$Z_iZ^j = \sum_{a,b,c} m_{ab}^j \mu_i^{bc} Z^a Z_c. \leqno (1.1.5)$$

Note that $HD(\Xi)$ is not a Hopf algebra.

\vskip .3cm

\noindent{\bf (1.2) Heisenberg double with respect to a Hopf pairing.}
It is convenient to introduce a version of the above formalism
which avoids dualizing possibly infinite-dimensional spaces. More precisely,
let $\Xi$ and $\Omega$ be Hopf algebras. A Hopf pairing of $\Xi$ and $\Omega$
is a bilinear map $\phi: \Xi\times\Omega\to {\bf C}$ satisfying the following
conditions:
$$\phi(\xi\xi', \omega) = \phi^{\otimes 2}(\xi\otimes\xi', \Delta(\omega)),\leqno (1.2.1)$$
$$ \phi(\xi, \omega\omega')=\phi^{\otimes 2} (\Delta(\xi), \omega\otimes
\omega')\leqno (1.2.2)$$
$$\phi(1,\omega) = \epsilon_{\Omega}(\omega), \quad \phi(\xi, 1) = \epsilon_\Xi (\xi).\leqno (1.2.3)$$
Ths conditions (1.2.1-2) simply mean that the multiplication and the comultiplication
are conjugate with respect to the pairing. We did not include here any conditions on the antipodes since
we will not need them. 

\vskip .2cm

If $\phi$ is a Hopf pairing of $\Xi$ and $\Omega$, we define the Heisenberg double  $HD(\Xi, \Omega, \phi)$
associated to $\phi$ to be the tensor product $\Omega\otimes_{\bf C} \Xi$ with the multiplication
given as follows: First, both $\Omega$ (realized as $\Omega\otimes 1$) and $\Xi$ (realized as $1\otimes \Xi$)
are required to be subalgebras. Second, for $\xi\in\Xi, \omega\in\Omega$ we impose the condition:
$$\xi\omega = ({\rm Id}\otimes \phi\otimes {\rm Id})(\Delta_{\Omega}(\omega)\otimes \Delta_\Xi (\xi)),
\leqno (1.3.4)$$
where ${\rm Id}\otimes \phi\otimes {\rm Id}: \Omega\otimes\Omega\otimes\Xi\otimes\Xi\to
\Omega\otimes\Xi$ is the map induced by the transposed pairing $\phi^{op}: \Omega\otimes\Xi\to {\bf C}$
on the second and third factors.

Thus, when $\Omega=\Xi^*$ and $\phi$ is the canonical pairing, we get the definition of (1.1). 

\vskip .3cm

\noindent {\bf (1.3) Hall and Ringel algebras.} We now describe a particular class of Hopf algebras, whose
Heisenberg doubles we will be interested in.

\vskip .1cm

\noindent Let $\cal A$ be an Abelian category. We will say that $\cal A$ has finite type, if
if, for any objects $A,B \in {\cal A}$ all the groups
$ {\rm Ext}^i_{\cal A}(A,B)$
have finite cardinality and are zero for almost all $i$. If $\cal A$ is of finite type, then, for any
objects $A,B,C\in {\cal A}$, the number of subobjects $A'\i A$ such that $A'\simeq A$ and $C/A'\simeq B$,
is finite. Denote this number $g_{AB}^C$.

\vskip .1cm

Let $H({\cal A})$ be the {\bf C}-vector space with basis $[A]$ parametrized by all the isomorphism
classes of objects $A\in{\cal A}$. The rule
$$[A]\circ [B] = \sum_C g_{AB}^C [C] \leqno (1.3.1)$$
makes $H({\cal A})$ into an associative algebra with unit $1=[0]$, see [R1-3]. This algebra
is called the Hall algebra of $\cal A$.

\vskip .1cm

Let ${\cal K}_0({\cal A})$ be the Grothendieck group of $\cal A$. For an object $A\in{\cal A}$ let $\bar A$
be its class in ${\cal K}_0({\cal A})$. The rule
$$\langle \bar A, \bar B\rangle = \sqrt{\prod_{i\geq 0} | {\rm Ext}^i_{\cal A}
 (A,B)|^{(-1)^i}}, \quad A,B\in {\cal A} \leqno (1.3.2)$$
extends uniquely (because of the behavior of Ext in exact sequences) to a bilinear form
${\cal K}_0({\cal A}) \otimes {\cal K}_0({\cal A})\to {\bf R}^*$ known as the Euler form. We will often write just
$\langle  A, B\rangle $ for $\langle \bar A, \bar B\rangle $, if $A,B$ are objects.
The symmetrization of the Euler form will be denoted by
$$(\alpha|\beta) = \langle \alpha,\beta\rangle \cdot \langle \beta,\alpha\rangle,
\quad \alpha, \beta\in {\cal K}_0({\cal A}). \leqno (1.3.3)$$
The twisted multiplication
$$[A]*[B] = \langle B, A\rangle \cdot [A]\circ [B] \leqno (1.3.4)$$
is still associative. We will denote $R({\cal A})$ and call the Ringel algebra of $\cal A$ the same
vector space as $H({\cal A})$ but with $*$ as multiplication.

\vskip .1cm

\noindent{\bf (1.3.5) Remark.} It was C.M. Ringel [R3] who first drew attention to the particular twist (1.3.4).
More generally, one can twist by any bilinear form on ${\cal K}_0({\cal A})$, and the associativity will
still be preserved. In fact, several such twists were used in earlier papers by G. Lusztig [Lu2-3],
without specially distinguishing the Euler form (1.3.2). 

\vskip .1cm

Let
${\bf C}[{\cal K}_0 {\cal A}]$
be the group algebra of
${\cal K}_0{\cal A}$,
with basis
$K_{\alpha}$,
$\alpha \in {\cal K}_0 {\cal A}$
and multiplication
$K_{\alpha} K_{\beta} = K_{\alpha +\beta}$.
Let us extend the algebra
$R({\cal A})$
by adding to it these symbols
$K_{\alpha}$
which we make commute with
$[A] \in R({\cal A})$
by the rule

$$[A]K_{\beta} = (\bar A|\beta) K_{\beta} [A]. \leqno (1.3.6)$$
Denote the resulting algebra
$B({\cal A})$.
So as a vector space
$B({\cal A}) \simeq {\bf C} [{\cal K}_0{\cal A}] \otimes_{\bf C} R({\cal A})$,
with
$K_{\alpha} \otimes [A] \mapsto K_{\alpha}A$
establishing the isomorphism. We will call $B({\cal A})$ the extended Ringel algebra of $\cal A$.

\vskip .2cm

Assume now that $\cal A$ satisfies two additional conditions. First, any object of $\cal A$ has only
finitely many subobjects. Second, ${\rm Ext}^i_{\cal A}(A,B)=0$ for all $A,B$ and all
$i>1$. We will state the second condition by saying that the homological dimension of $\cal A$ is less or
equal to 1 and write ${\rm hd}({\cal A})\leq 1$. The next statement follows from results of
Green [Gr], see [X1] [Kap] for a detailed deduction.

\proclaim (1.3.7) Theorem. $B({\cal A})$ is a Hopf algebra with respect to the
comultiplication given on generators by
$$\Delta(K_{\alpha}) =K_{\alpha} \otimes K_{\alpha},$$
$$\Delta([A]) = \sum_{A^{\prime} \subset A} \langle A/A^{\prime}, A^{\prime}\rangle  {| {\rm Aut} (A^{\prime})|\cdot |{\rm Aut} (A/A^{\prime})|\over |{\rm Aut} (A)|} [A^{\prime}] \otimes K_{A^{\prime}} [A/A^{\prime}],$$ 
the counit $\epsilon : B({\cal A}) \rightarrow {\bf C}$
given by
$$\epsilon (K_{\alpha} [A]) =   1, \quad {\rm if}\quad { A=0} \quad {\rm and}\quad 
0,  \quad {\rm if} \quad { A \ne 0} $$
and the antipode
$S : B(\cal A) \rightarrow B({\cal A})$
given by
$$S(K_{\alpha} [A]) = \sum^{\infty}_{n=1} (-1)^n \sum_{A_0 \subset \ldots \subset A_n=A} \prod^n_{i=1} \langle A_i/A_{i-1},A_{i-1}\rangle  {\prod^n_{j=0} |{\rm Aut} (A_j/A_{j-1})|
\over  |{\rm Aut} (A)|} \cdot $$
$$\cdot [A_0] * [A_1/A_0] * \ldots * [A_n/A_{n-1}] \cdot K^{-1}_{\alpha} K^{-1}_A $$
where 
$A_0 \subset \ldots \subset A_n =A$
runs over arbitrary chains of strict
$(A_i \ne A_{i+1})$
inclusions of length $n$.

\vskip .2cm

\noindent {\bf (1.4) The Hopf pairing on $B({\cal A})$.} The elements $K_\alpha [A]$ form a {\bf C}-basis
of $B({\cal A})$. Let us define a bilinear pairing $\phi: B({\cal A})\times B({\cal A})\to {\bf C}$ by
putting

$$\phi (K_{\alpha}[A], K_{\beta} [B]) = (\alpha|\beta) ([A], [B])  = 
{ {(\alpha|\beta) \delta_{[A],[B]}} \over {| {\rm Aut} (A)|}}. \leqno (1.4.1)$$

\proclaim (1.4.2) Proposition.
The pairing $\phi$ is a Hopf pairing on $B({\cal A})$.

\noindent {\sl Proof:}
We need to prove the equality (1.3.1) (the other equality (1.3.2) will then follow by symmetry).
In other words, we need to prove that
$$\phi (K_{\alpha}[A]K_{\beta}[B], K_{\gamma}[C]) = \phi^{\otimes 2} (K_{\alpha}[A]
 \otimes K_{\beta} [B], \Delta(K_{\gamma} [C])) \leqno (1.4.3)$$
 To prove this, notice that the left hand side is
$$ (\bar A|\beta) \phi(K_{\alpha+\beta}[A][B], K_{\gamma}[C]) =
 (\bar A|\beta)(\alpha +\beta|\gamma)\phi ([A] * [B], [C]) = $$
$$= (\bar A|\beta)(\alpha+\beta|\gamma) \sum_{C^{\prime} \subset C} 
\langle C/C^{\prime}, C^{\prime}\rangle { {| {\rm Aut}(C^{\prime})| 
\cdot | {\rm Aut}(C/C^{\prime})|}\over  {| {\rm Aut} (C)|}} \cdot $$
$$\cdot \phi (A,C^{\prime}) \cdot \phi (B,C/C^{\prime}), \leqno (1.4.4) $$
while the right hand side of (1.4.3) is

$$ \phi^{\otimes 2} \biggl(K_{\alpha} [A] \otimes K_{\beta} [B], \quad \sum_{C^{\prime}
 \subset C} \langle C/C^{\prime}, C^{\prime}\rangle { {| {\rm Aut}(C^{\prime})| 
\cdot | {\rm Aut}(C/C^{\prime})|}\over {| {\rm Aut} (C)|} }\cdot $$
$$\cdot K_{\gamma}[C^{\prime}] \otimes K_{\bar C+\gamma} 
[C^ {\prime\prime}]\biggl) = \leqno (1.4.5)$$
$$= \sum_{C^{\prime} \subset C} \langle C/C^{\prime}, C^{\prime}\rangle
 { {| {\rm Aut}(C^{\prime})| \cdot |{\rm Aut}(C/C^{\prime})|} \over  
{|{\rm Aut}(C)|}} (\alpha|\gamma)(\beta |\gamma)(\beta|\bar C^{\prime}) \cdot $$
$$\cdot \phi ([A], [C^{\prime}]) \cdot \phi ([B], [C/C^{\prime}]). $$
Notice now that in order that
$\phi (A,C^{\prime}) \ne 0$,
we should have 
$A \simeq C^{\prime}$,
and under this assumption the corresponding summands in (1.4.4) and 
(1.4.5) coincide. Proposition is proved.

\vskip .3cm

\noindent {\bf (1.5) The Heisenberg double of $B({\cal A})$.} Let ${\rm Heis}({\cal A})$
be the Heisenberg double \hfill\break 
$HD(B({\cal A}), B({\cal A}), \phi)$, see (1.2). We will denote its generators    as follows:
$$Z^+_A = 1\otimes [A], \quad Z^-_A = [A]\otimes 1,\quad K_\alpha = 1\otimes K_\alpha,
\quad K_\alpha^- = K_\alpha\otimes 1. \leqno (1.5.1)$$
Thus, the $Z^+_A$ together with the $K_\alpha$, form a copy of $B({\cal A})$ inside 
${\rm Heis}({\cal A})$, and the same for the $Z^-_A$ with the $K^-_\alpha$. To find the cross-relations
between the plus and minus generators more explicitly, let us introduce the following
notations. For any objects $A,B,M,N\in{\cal A}$ let ${\cal F}_{AB}^{MN}$ be the set of
all exact sequences
$$0\to M\to B\buildrel \psi\over\to A\to N\to 0
\leqno (1.5.2)$$
and by $\gamma_{AB}^{MN}$ the quotient $|{\cal F}_{AB}^{MN}|/|{\rm Aut}(A)|\cdot |{\rm Aut}(B)|$.
The following statement, with its proof, is an adaptation of Proposition 6.2.12 from [Kap]
which, however, used a more cumbersome approach.

\proclaim (1.5.3) Proposition. We have the following equalities in ${\rm Heis}({\cal A})$:
$$Z^+_AZ^-_B = \sum_{M,N} \langle \bar B -\bar M, \bar M\rangle\cdot \langle \bar N, \bar B-\bar M\rangle \cdot
\gamma_{AB}^{MN} Z^-_M K_{\bar B-\bar M} Z^+_N=\leqno (1.5.4)$$
$$= \sum_{M,N} \gamma_{AB}^{MN} 
\langle \bar B-\bar M, \bar M-\bar N\rangle \cdot Z^-_M Z^+_N K_{\bar B-\bar M},$$
$$Z^-_AK_\alpha = (A|\alpha)^{-1}K_\alpha Z^-_A, \quad Z^+_A K^-_\alpha = K_\alpha^- Z^+_A,
\quad K_\alpha^+ K_\beta^- = (\alpha|\beta) K_\beta^-K_\alpha^+.\leqno (1.5.5)$$

\noindent {\sl Proof:} By (1.3.4) and the definition of the $Z^\pm_A$, we have
$$Z^+_AZ^-_B = ({\rm Id}\otimes\phi\otimes {\rm Id}) (\Delta([B])\otimes\Delta([A])) =$$
$$= \sum_{M,I} \sum_{I',N} \langle I, M\rangle \cdot \langle N, I'\rangle g_{MI}^B \cdot g_{I'N}^A
{ |{\rm Aut}(M)|\cdot |{\rm Aut}(I)|\cdot |{\rm Aut}(I')|\cdot |{\rm Aut}(N)|\over
 |{\rm Aut}(B)|\cdot |{\rm Aut}(A)|} \times$$
$$\times \phi(K_M[I], [I']) \cdot [M]\otimes K_{I'}N.$$
Note that $\phi(K_M[I], [I']) = \delta_{[I], [I']}/|{\rm Aut}(I)|$, so we can just put $I'=I$. Further,
 for any three objects $A,B,C$ let ${\cal E}_{AB}^C$ be the
set of exact sequences
$$0\rightarrow A\rightarrow C\rightarrow B \rightarrow 0.$$
Thus
$$\gamma_{AB}^{MN} = { |{\cal F}_{AB}^{MN}|\over |{\rm Aut}(A)|\cdot
|{\rm Aut}(B)|}, \quad g_{AB}^C = {|{\cal E}_{AB}^{C}|\over |{\rm Aut}(A)|\cdot
|{\rm Aut}(B)|}.$$
Notice now that
$${\cal F}_{AB}^{MN} \quad = \quad \coprod_{I\in {\rm Ob}({\cal A})/{\rm Iso}}
\bigl( {\cal E}_{MI}^B \times {\cal E}_{IN}^A \bigr) \bigl/ {\rm Aut}(I),
\leqno (1.5.6)$$
with ${\rm Aut}(I)$ acting freely. This just means that any
 long exact sequence (1.5.2) can be split into two short sequences with $I =
 {\rm Im}(\psi)$. Let ${\cal F}_{AB}^{MN}(I)$ be the $I$th part of the disjoint union (1.5.6).
Then, by taking all the above equalities into account, we find:
$$Z^+_AZ^-_B = \sum_{M,I,N} \langle I, M\rangle \cdot \langle N, I\rangle { |{\cal F}_{AB}^{MN}(I)|\over
|{\rm Aut}(A)|\cdot |{\rm Aut}(B)|} Z^+_MK_IZ^-_N,$$
and to get the claimed equality (1.5.4), it remains to notice that $\bar I = \bar B -\bar M$ once 
 ${\cal F}_{AB}^{MN}(I)\neq\emptyset$. The equalities (1.5.5) are obtained in a straightforward way.
Proposition is proved.

\vfill\eject

\centerline {\bf \S 2. Heisenberg doubles and tilting.}

\vskip 1cm

\noindent {\bf (2.1) Generalities on derived categories.} Let $\cal A$ be an Abelian category. 
By $C^b({\cal A})$
we denote the category of bounded complexes $A^\bullet = (A^i, d_i = d_{i, A}: A^i\rightarrow A^{i+1})$ over $\cal A$.
The shifted complex $A^\bullet[n], n\in {\bf Z}$, is defined by $(A^\bullet [n])^i = A^{n+i}$, $d_{i, A[n]} = 
(-1)^n d_{i,A}$. The homology objects of a complex $A^\bullet$ are denoted by $H^i(A^\bullet) = {\rm Ker}(d_i)/
{\rm Im}(d_{i-1})$. We denote by $D^b({\cal A})$ the bounded derived category of $\cal A$. It is obtained from
$C^b({\cal A})$ by formally inverting quasi-isomorphisms. If $A,B$ are two objects of $\cal A$ (regarded as complexes concentrated in degree 0), then
$${\rm Hom}_{D^b({\cal A})} (A, B[i]) = {\rm Ext}^i_{\cal A}(A,B).\leqno (2.1.1)$$

For any triangulated category $\cal D$ and any morphism $f: X\rightarrow Y$ in $\cal D$ we will denote by
${\rm Cone}(f)$ the isomorphism class of third terms $Z$ of possible exact triangles
$$ X\buildrel f\over\longrightarrow Y\rightarrow Z \rightarrow X[1].$$

Let us say that $\cal A$ has homological dimension $d$ and write ${\rm hd}({\cal A})\leq d$ if ${\rm Ext}^i_{\cal A}(A,B) = 0$
for any $A,B\in {\cal A}$ and any $i>d$. We say that $\cal A$ has finite homological dimension (${\rm hd}({\cal A}) < \infty$)
if ${\rm hd}({\cal A})\leq d$ for some $d$.

\proclaim (2.1.2) Proposition. If ${\rm hd}({\cal A})\leq 1$, then
 each  object of $D^b({\cal A})$ 
is isomorphic
to the complex $H^\bullet(A^\bullet)$ formed by the cohomology of $A^\bullet$ and equipped with zero differential.

This proposition is interesting for us because it gives a very explicit description of $D^b({\cal A})$ as a category. Indeed, given any two complexes $A^\bullet, B^\bullet$ with zero differential,
 we have $A^\bullet = \bigoplus_{i\in {\bf Z}} A^{-i}[i]$ and similarly for
$B^\bullet$, so by (2.1.1)
$${\rm Hom}_{D^b({\cal A})}(A^\bullet, B^\bullet) = \bigoplus_i {\rm Hom}_{\cal A}(A^i, B^i) \quad \oplus\quad
\bigoplus_i {\rm Ext}^1_{\cal A}(A^i, B^{i-1})\leqno (2.1.3)$$
In other words, a morphism $f: A^\bullet\rightarrow B^\bullet$ is the same as a sequence  of components
$f_i^{\rm Hom} \in {\rm Hom}_{\cal A}(A^i, B^i)$ and $f_i^{\rm Ext}\in {\rm Ext}^1_{\cal A}(A^i, B^{i-1})$. In the sequel
(including proof of (2.1.2) we will
use this notation for the components of a morphism.

\vskip .2cm

\noindent {\sl Proof of (2.1.2):} The proposition for $D^b({\cal A})$ seems to be well known. A simple argument (pointed out to me by A. Bondal) is to show by induction
that any bounded complex is quasiisomorphic
to the sum of its last cohomology object and
its canonical truncation just below this object. Here
we include,  for completeness sake,
a different proof which
does nt use induction and so is applicable not just to $D^b({\cal A})$ but also to other types of derived categories
(unbounded, periodic etc.).

\vskip .1cm

Let $(A^\bullet, d)$ be a complex, and let $K^\nu = {\rm ker}(d_\nu)$, $I^\nu = {\rm Im}(d_{\nu -1})$. We have short exact sequences
$$0\rightarrow K^\nu\buildrel \epsilon_\nu\over\rightarrow A^\nu\buildrel \pi_\nu\over\rightarrow I^{\nu+1}\rightarrow 0,
\leqno (2.1.4)$$
(with $\pi_\nu$ induced by $d_\nu$) which fit together into an exact sequence of complexes
$$0\rightarrow (K^\bullet, 0) \buildrel\epsilon\over\rightarrow (A^\bullet, d)\buildrel \pi\over\rightarrow (I^\bullet [1], 0)\rightarrow 0.
\leqno (2.1.5)$$
This sequence gives rise to an exact triangle in $D^b({\cal A})$, in particular, we get the boundary map
$\delta: (I^\bullet [1], 0) \rightarrow (K^\bullet [1], 0)$ such that $(A^\bullet, d) \simeq {\rm Cone} (\delta)[-1]$.
We want to compare (2.1.5) with the short exact sequence
$$0\rightarrow (I^\bullet, 0) \buildrel \phi\over\rightarrow (K^\bullet, 0)\buildrel\psi\over\rightarrow (H^\bullet, 0)\rightarrow 0
\leqno (2.1.6)$$
defining the cohomology $H^\bullet = H^\bullet(A^\bullet)$. This sequence implies that $H^\bullet \simeq {\rm Cone}(\phi)$.
Note that $\phi$ has only ${\rm Hom}$-components $\phi_\nu = \phi_\nu^{\rm Hom}: I^\nu \hookrightarrow K^\nu$.
Note also that $\phi_\nu = \delta[-1]_{\nu}^{\rm Hom}$. Indeed, the latter map is just the boundary homomorphsm
in the long cohomology sequence of (2.1.5), and this homomorphism is straightforwardly found to be $\phi_\nu$. 
As to $\delta[-1]^{\rm Ext}_\nu$, it is the class of the extension (2.1.4) and may well be non-zero. However,
the situation is saved by the following fact.

\proclaim (2.1.7) Lemma. There exists an automorphism $W$ of $(K^\bullet, 0)$ in the derived category such that $W\delta[-1] = \phi$.

The lemma implies our proposition by the axiom TR2 of triangulated categories (a commutative square extends to a morphism
of triangles). 

\vskip .2cm

\noindent {\sl Proof of the lemma:} Since ${\rm hd}({\cal A}) = 1$ and $\phi_\nu$ is injective, the restriction map
$$\phi_\nu^*: {\rm Ext}^1_{\cal A}(K^\nu, K^{\nu-1}) \rightarrow {\rm Ext}^1_{\cal A}(I^\nu, K^{\nu-1})$$
is surjective. Let $w_\nu \in {\rm Ext}^1_{\cal A}(K^\nu, K^{\nu-1})$ be any element mapping into $\delta[-1]_\nu^{\rm Ext}$.
Now define $W: (K^\bullet, 0)\rightarrow (K^\bullet, 0)$ to have the components $W_\nu^{\rm Hom} = {\rm Id}$
and $W_\nu^{\rm Ext} = - w_\nu$. Then $W$ is an isomorphism since it is given (with respect to the decomposition
$K^\bullet = \bigoplus K^{-i}[i]$) by a triangular matrix with identities on the diagonal. One immediately sees that $W\delta[-1]=\phi$,
since the Ext-terms in the composition will cancel. Lemma and Proposition 2.1.2 are proved.

\proclaim (2.1.8) Corollary. If ${\rm hd}({\cal A})\leq 1$, then each indecomposable object of
$D^b({\cal A})$ has the form $A[i]$ where $i\in{\bf Z}$ and $A$ is an indecomposable object of $\cal A$. 

\vskip .2cm

\noindent {\bf (2.2) The category $D^{[-1,0]}({\cal A})$ and tiltings.} Let $\cal A$ be as before. Denote by
$D^{[-1,0]}({\cal A})$ the full subcategory in $D^b({\cal A})$ formed by complexes situated in degrees
$-1, 0$ only. Given two Abelian categories $\cal A$ and $\cal B$, we will call an equivalence
$F: D^b({\cal A})\to D^b({\cal B})$ of triangulated categories a {\it tilting}, if
$F({\cal A})\i D^{[-1,0]}({\cal B})$. This condition is satisfied for equivalences given
by the so-called tilting modules [Ha]. 

From now on we assume that all the Abelian categories we consider have homological dimension less or equal to 1
and satisfy all the finiteness conditions of \S 1. Thus, for any such category $\cal A$ we have the Hopf
algebra $B({\cal A})$ and its Heisenberg double ${\rm Heis}({\cal A})$. 
For two objects $A,B\in{\cal A}$ denote
$$[A,B] = |{\rm Hom}(A,B)|^{+1/2} \cdot |{\rm Ext}^1 (A,B)|^{+1/2}. \leqno (2.2.1)$$
Because of two plus signs in the exponents, this quantity does not descend to the
Grothendieck group. 

Let us associate to any object (complex with zero differential)
$A^\bullet = A^{-1}[1]\oplus A^0\in D^{[-1,0]}({\cal A})$ the following element of
${\rm Heis}({\cal A})$:
$$Z(A^\bullet) = { Z^-_{A^{-1}}K^{-1}_{A^{-1}}Z^+_{A^0}\over \langle A^{-1}, A^{-1}\rangle
\cdot [A^0, A^{-1}]} = { Z^-_{A^{-1}}K^{-1}_{A^{-1}}Z^+_{A^0}\over \langle A^{-1}, A^{-1}\rangle
\cdot \langle A^0, A^{-1}\rangle\cdot |{\rm Ext}^1(A^0, A^{-1})|}. \leqno (2.2.2)$$

Now we can formulate the main result of this section.

\proclaim (2.3) Theorem. If $F: D^b({\cal A})\to D^b({\cal B})$ is a tilting, then
the correspondence $[A]\to Z(F(A))$ gives an injective homomorphism of algebras $F_*: R({\cal A})
\to {\rm Heis}({\cal B})$.

Before starting the proof, we do some preliminary work in the next subsection.

\vskip .3cm

\noindent{\bf (2.4) Counting exact triangles.} 
If $G$ is a finite group acting on a finite set $X$, we will call the ratio $|X|/|G|$ the orbifold number of elements of $X$
 modulo $G$.
This number is the same as $\sum_{ \{x\} \in X/G} 1/|{\rm Stab}(x)|$, the sum being over $G$-orbits on $X$, and $x$ being one representative
chosen for each orbit $\{x\}$.

Let $\cal A$ be as above.  For any three objects $A^\bullet, B^\bullet, C^\bullet$ of $D^b({\cal A})$
 we denote by 
$g_{A^\bullet, B^\bullet}^{C^\bullet}$ the orbifold number of exact triangles
$$A^\bullet\rightarrow C^\bullet\rightarrow B^\bullet \buildrel\partial\over\rightarrow A^\bullet[1]\leqno (2.4.1)$$
modulo ${\rm Aut}(A^\bullet)\times {\rm Aut}(B^\bullet)$.

If $A^\bullet, B^\bullet, C^\bullet$ are three bounded {\bf Z}-graded
objects of $\cal A$ (i.e., complexes with zero differentials), we denote by $\gamma_{A^\bullet, B^\bullet}^{C^\bullet}$
the orbifold number of long exact sequences 
$$...\rightarrow A^i \rightarrow C^i \rightarrow  B^i\buildrel
 \partial_i\over\rightarrow A^{i+1}\rightarrow
...\leqno (2.4.3)$$
modulo $\prod_i {\rm Aut}(A^i)\times {\rm Aut}(B^i)$.

If $A,B,C\in {\cal A}$ are three objects considered as {\bf Z}-graded objects (complexes with zero differential)
concentrated in degree 0, then clearly $\gamma_{AB}^C = g_{AB}^C$ coincides with the number introduced in
(1.3). For general complexes with zero differential, $g_{A^\bullet, B^\bullet}^{C^\bullet}$
differs from $\gamma_{A^\bullet, B^\bullet}^{C^\bullet}$. We will need one particular case when these
numbers can be easily compared.

\proclaim (2.4.3) Proposition. Let $A,B,M,N$ be any objects of $A$. Then
$$g_{A, B[1]}^{M[1]\oplus N} =\gamma_{A, B[1]}^{M[1]\oplus N}  \cdot |{\rm Ext}^1(N,M)|.$$

\noindent {\sl Proof:} The number $\gamma_{A, B[1]}^{M[1]\oplus N}$ counts exact sequences 
$$0\to M\buildrel u\over\to B\buildrel\varphi\over\to A\buildrel v\over\to N\to 0\leqno (2.4.4)$$
and thus is equal to
$${ \left| \bigl\{ \varphi: B\to A| \,\,  {\rm Ker}(\varphi)\simeq M, \, {\rm Coker}(\varphi)\simeq N\bigr\}
\right| \cdot |{\rm Aut}(M)|  \cdot |{\rm Aut}(N)|\over   |{\rm Aut}(A)| \cdot |{\rm Aut}(B)| }.\leqno (2.4.5)$$
The number $g_{A^\bullet, B^\bullet}^{C^\bullet}$ counts exact triangles
$$A\buildrel\alpha\over\to M[1]\oplus N \buildrel\beta\over\to B[1]\buildrel \varphi[1]\over\to 
A[1].\leqno (2.4.6)$$
Of course, any such triangle gives rise to a sequence of the form (2.4.4), but the correspondence is not
bijective. 

More precisely, let us fix $\varphi: B\to A$.  By (2.1.2), in order that a triangle (2.4.6) with the boundary map $\varphi[1]$
exists, it is necessary and sufficient that ${\rm Ker}(\varphi)\simeq M$, ${\rm Coker}(\varphi)\simeq N$.
If $\varphi$ satisfies this property, then the axiom TR2 of triangulated categories
implies that a triangle (2.4.6) can be constructed uniquely modulo an isomorphism of triangles
identical on $A$ and $B[1]$. Such an isomorphism is the same as just an automorphism of $M[1]\oplus N$.
So
$$g_{A, B[1]}^{M[1]\oplus N} = 
 { \left| \bigl\{ \varphi: B\to A| \,\,  {\rm Ker}(\varphi)\simeq M, \, {\rm Coker}(\varphi)\simeq N\bigr\}
\right|\cdot |{\rm Aut}(M[1]\oplus N)|\over |{\rm Aut}(A)|\cdot |{\rm Aut}(B)|\cdot |{\rm Stab}(\varphi)|},
\leqno (2.4.7)$$
where ${\rm Stab}(\varphi)\i {\rm Aut}(M[1]\oplus N)$ is the subgroup of automorphisms which,
together with the identities of $A$ and $B[1]$, give an automorphism of the triangle (2.4.6).
In fact, we claim that ${\rm Stab}(\varphi) = \{{\rm Id}\}$. To see this,
note that ${\rm Aut}(M[1]\oplus N)$ is the block matrix group
$$\pmatrix {{\rm Aut}(M)& {\rm Ext}^1(N,M)\cr 0& {\rm Aut}(N)}.\leqno (2.4.8)$$
Let $\psi=(\psi_M, \psi_{NM}, \psi_N)$ be an element of ${\rm Stab}(\varphi)$, so
$\psi_M\in {\rm Aut}(M)$, $\psi_{NM}\in {\rm Ext}^1(N,M)$, $\psi_N\in {\rm Aut}(N)$. 
Then $\psi_M={\rm Id}$ since $u: M\to B$ is an injection, and $\psi_N={\rm Id}$ since
$v: A\to N$ is a surjection and since $\psi$ together with ${\rm Id}_A, {\rm id}_B$ forms a morphism of triangles.
Further, the commutativity of the diagram
$$\matrix{&A &\buildrel\alpha\over\longrightarrow& M[1]\oplus N &\cr
{\rm Id}&\big\downarrow &&\big\downarrow & \psi=\pmatrix{1&\psi_{MN}\cr 0&1} \cr
&A&\buildrel\alpha\over\longrightarrow& M[1]\oplus N }$$
means that $\psi_{MN}\circ v=0$ in ${\rm Ext}^1(A,M)$. But because ${\rm hd}({\cal A})\leq 1$,
the surjection $v: A\to N$ induces an injection ${\rm Ext}^1(N,M)\to
{\rm Ext}^1(A,M)$, so $\psi_{MN}=0$. This proves our claim that ${\rm Stab}(\varphi)=\{{\rm Id}\}$.
The proposition follows now by comparing (2.4.7) with (2.4.5) and taking into account the factorization
(2.4.8). 

\vskip .3cm

\noindent {\bf (2.5) Proof of Theorem 2.3.} 
 First of all let ${\cal A}_i \i {\cal A}$ be the full subcategory
of $A$ such that $F(A)\in {\cal B}[i], i=0, 1$. Notice that for $A_i\in {\cal A}_i$
we have, denoting $B_i = F(A_i)[-i]\in {\cal B}$:
$${\rm Hom}_{\cal A}(A_{1}, A_0) = {\rm Hom}_{D^b({\cal B})}(B_{1}[1], B_0) = 0,\leqno (2.5.1a)$$
$$ {\rm Ext}^1_{\cal A}(A_0, A_{1}) = {\rm Hom}_{D^b({\cal B})}(A_0, A_{1}[2]) = 
{\rm Ext}^2_{\cal B}(A_0, A_{1})=0.\leqno (2.5.1b) $$
Let us denote by $F_*: R({\cal A})\to {\rm Heis}({\cal A})$ the unique {\bf C}-linear
map taking $[A]$ to $Z(F(A))$. Its injectivity is clear from the behavior on basis vectors,
so the main task is to prove that $F_*$ is an algebra homomorphism,
i.e., that
$$F_*([A']* [A'']) = F_*([A']) F_*([A'']), \quad \forall A', A''\in {\cal A}. \leqno (2.5.2) $$
The proof will be done in two steps. The first is given in the next proposition.

\proclaim (2.5.3) Proposition. The equality (2.5.2) holds in the case when $A'=A'_i\in {\cal A}_i$,
$A''=A''_j\in {\cal A}_j$ for some $i,j\in \{0,1\}$.

The second step is to deduce the general case from these particular cases. Let us first
explain how this is done and then prove Proposition 2.5.3. Namely, by (2.1.8) each indecomposable
object of $\cal A$ lies in one of the ${\cal A}_i$. So each $A\in {\cal A}$
can be uniquely written as $A=A_0\oplus A_{1}$ with $A_i\in {\cal A}_i$. By (2.5.1)
we have
$$[A] = |{\rm Hom}(A_0, A_1)|^{-1/2} [A_{1}]*[A_0].\leqno (2.5.4)$$
This means that we have a kind of normal form for elements of $R({\cal A})$, i.e., the
map $R({\cal A}_{1})\otimes_{\bf C} R({\cal A}_0)\to R({\cal A})$ given by the
multiplication, is an isomorphism of vector spaces.
 So the second step is accomplished by the next easy lemma.

\proclaim (2.5.5) Lemma. Let $R$ be a {\bf C}-algebra and $R_0, R_{1}\i R$ be subalgebras
such that the multiplication induced an isomorphism of vector spaces
$R_1\otimes_{\bf C}R_0\to R$. Let $S$ be another algebra
and $\phi: R\to S$ be a {\bf C}-linear map. Suppose that the equality $\phi(a'a'')
=\phi(a')\phi(a'')$ holds whenever $a'\in R_i, a''\in R_j$ for some $i,j$.
Then it holds for any $a',a''$, i.e., $\phi$ is an algebra homomorphism.

\noindent {\sl Proof of (2.5.5):} It is enough, by linearity, to consider the case when
$a'=a'_{1}a'_0, a''=a''_{1}a''_0$ with $a'_i, a''_i\in R_i$. Then, by our assumptions,
$$\phi(a')\phi(a'') = \phi(a'_{1})\phi(a'_0)\phi(a''_{1})\phi(a''_0).$$
Let us write $a'_0a''_{1} = \sum_j \alpha^{(j)}_{1}\alpha^{(j)}_0$ with $\alpha^{(j)}_i\in R_i$,
using our assumption on $R$.
Then, by our assumptions on $\phi$,
$$\phi(a'_0)\phi(a''_{1}) = \phi (a'_0)a''_{1}) = \sum_j \phi(\alpha^{(j)}_{1}) 
\phi(\alpha^{(j)}_0),$$
and so
$$\phi(a')\phi(a'') = \sum_j \phi(a'_{1}\alpha^{(j)}_{-1}) \phi(\alpha^{(j)}_0a''_0) = 
\sum \phi(a'_{1}\alpha^{(j)}_{1}\alpha^{(j)}_0a''_0) =\phi(a'_{1}a'_0a''_{1}a''_0)=
\phi(a'a''),$$
as claimed.

\vskip .2cm

We now prove Proposition 2.5.3 by considering all four possibilities for
$i,j$.

\vskip .2cm

\noindent\underbar {Case (0,0)}: $A'=A'_0, A''=A''_0\in {\cal A}_0$. Let $F(A')=B',
F(A'')=B''$. Then $F_*(A')=Z^+_{B'}$ and $F_*(A'')=Z^+_{B''}$. Since $F$ is an embedding
of an admissible Abelian category, it establishes a bijection between short exact sequences
$$0\to A'\to A\to A''\to 0$$
in $\cal A$ and short exact sequences
$$0\to B'\to B\to B''\to 0$$
in $\cal B$ (since both kinds of sequences are interpreted as exact triangles in the same
triangulated category $D^b({\cal B})$). Thus the equality (2.5.2) holds.

\vskip .2cm

\noindent\underbar {Case (1,1)}: $A'=A'_{1}, A''=A''_{1}\in {\cal A}_{1}$.
Let $F(A')=B'[1], F(A'')=B''[1]$. If
$$0\to A'\to A\to A''\to 0$$
is a short exact sequence in $\cal A$, then $A\in {\cal A}_1$, and denoting $B=F(A)[-1]$,
we have $g_{A'A''}^A=g_{B'B''}^B$. By our definition,
$F_*([A']) = Z^-_{B'}K_{B'}^{-1}  \langle B', B'\rangle ^{-1}$, and similarly for
$F_*([A''])$. Therefore
$$F_*([A']) F_*([A'']) = {Z^-_{B'}K_{B'}^{-1} Z^-_{B''} K_{B''}^{-1}\over \langle B',B'\rangle\cdot
\langle B'',B''\rangle} =  {Z^-_{B'}Z^-_{B''} K^{-1}_{\bar{B}' + \bar{B}''}\over 
\langle B',B'\rangle\cdot \langle B'',B''\rangle\cdot (B'|B'')} =$$
$$={ \sum_{B} g_{B'B''}^B \langle B'',B'\rangle Z^-_B K^{-1}_B\over \langle B ,B \rangle}
=F_*([A'] *[A'']).$$
Here we used the fact that $\bar B= \bar B' +  \bar B''$ whenever  $g_{B'B''}^B\neq 0$.

\vskip .2cm

\noindent\underbar {Case (1,0)}: $A'=A'_1\in {\cal A}_0, A''=A''_0\in {\cal A}_1$. Let
$F(A')=B'[1], F(A'')=B''$. We have, by definition of $F_*$
$$ F_*([A'])F_*([A'']) = Z^-_{B'} K^{-1}_{B'} Z^+_{B''} \cdot \langle B', B'\rangle ^{-1},$$
while, by (2.5.4), (2.5.1)  and by definition of $F_*$, 
$$F([A']*[A'']) = |{\rm Hom}(A'', A')|^{1/2} F([A'\oplus A'']) = 
{|{\rm Ext}^1(B'', B')|^{1/2} Z^-_{B'}K^{-1}_{B'}Z^+_{B''}\over
\langle B',B'\rangle \cdot |{\rm Hom}(B'', B')|^{1/2}\cdot |{\rm Ext}^1(B'', B')|^{1/2}},$$
which is the same as the previous quantity once we recall that 
${\rm Hom}(B'', B') = {\rm Ext}^{-1}(A'', A')=0$.

\vskip .2cm

\noindent \underbar{Case (0,1)}: $A'=A'_0\in {\cal A}_0, A''=A''_1\in{\cal A}_1$. Let
$F(A')=B', F(A'')=B''[1]$. We want to verify that 
$$ F_*([A']) F_*([A'']) = F_*([A']* [A'']) = \sum_{A\in {\cal A}} g_{A'A''}^A F_*([A]).$$
Note that the $A$ entering the last sum, may not lie in any of the ${\cal A}_i$. However,
if $A$ is included into an exact sequence
$$0\to A'\to A\to A''\to 0,$$
then $F(A)$, which necessarily has the form
 $M[1]\oplus N$ for some $M,N\in {\cal B}$, is included into an exact triangle
$$B'\to M[1]\oplus N \to B''[1] \to B'[1],$$
and this correspondence is a bijection, i.e., $g_{A',A''}^A = g_{B', B''[1]}^{M[1]\oplus N}$.
If such a triangle exists, then in ${\cal K}_0({\cal B})$ we have
$\bar M-\bar N =  \bar B''-\bar B'$, as it follows from the corresponding 4-term exact
sequence. Now, by  (2.2.2), we have 
$$F_*([A'])F_*([A'']) ={ Z^+_{B'} Z^-_{B''}K^{-1}_{B''}\over \langle B'', B''\rangle }=$$
$$=\sum_{M,N\in {\cal B}} \gamma_{B', B''[1]}^{M[1]\oplus N}
 { \langle \bar B''-\bar M, \bar M\rangle\cdot
\langle \bar N, \bar B''-\bar M\rangle\over \langle B'', B''\rangle}
Z^-_M K_{\bar B''-\bar M} Z^+_N K^{-1}_{B''}=$$
$$=\sum_{M,N} \gamma_{B', B''[1]}^{M[1]\oplus N}
{\langle B'', M\rangle \cdot \langle N, B''\rangle \over \langle B'', B''\rangle \cdot 
\langle M, M\rangle \cdot \langle N, M\rangle \cdot (B''|N)} Z^-_M K_M^{-1}Z^+_N.$$
Note that the quantity represented by the fraction in the last expression, can be written as
$${ \langle \bar B'', \bar M-\bar N  \rangle \over \langle B'', B''  \rangle 
\cdot \langle M, M \rangle  \cdot  \langle N, M  \rangle } = 
{\langle   \bar B'', \bar B''-\bar B' \rangle \over \langle B'', B''  \rangle 
\cdot \langle M, M \rangle  \cdot  \langle N, M  \rangle }.$$
Therefore, applying Proposition 2.4.3, we find:
$$F_*([A'])F_*([A'']) = \sum_{M,N} \gamma_{B', B''[1]}^{M[1]\oplus N}
{\langle   \bar B'', \bar B''-\bar B' \rangle \over \langle B'', B''  \rangle 
\cdot \langle M, M \rangle  \cdot  \langle N, M  \rangle } Z^-_M K_M^{-1}Z^+_N =$$
$$= \langle B''[1], B'\rangle \sum_{M,N} g_{B', B''[1]}^{M[1]\oplus N}
{Z^-_M K_M^{-1}Z^+_N \over \langle M, M\rangle\cdot \langle N, M\rangle \cdot
|{\rm Ext}^1({N,M})| } = $$
$$=\langle A'', A'\rangle \sum_{A\in {\cal A}} g_{A', A''}^A F_*([A]) =
F_*([A']*[A'']).$$
This finishes the proof of Theorem 2.3.

\vfill\eject

\centerline {\bf \S 3. The lattice algebra and the full derived category.}

\vskip 1cm

\noindent {\bf (3.1) Definition of the lattice algebra.} Let $\cal A$ be an Abelian category 
with ${\rm hd}({\cal A})\leq 1$ satisfying all the finiteness conditions of \S 1. Let
$B({\cal A})\supset R({\cal A}) $ be its extended Ringel algebra and ${\rm Heis}({\cal A})$
be the Heisenberg  double of $B({\cal A})$. As we saw in \S 2, ${\rm Heis}({\cal A})$ is
naturally related to the subcategory $D^{[-1,0]}({\cal A})$ in the derived category $D^b({\cal A})$.
Hereby the two copies of $\cal A$ inside $D^{[-1,0]}({\cal A})$, given by complexes concentrated in degree 
$0$ (resp. $(-1)$), give rise to two copies of $R({\cal A})$ in the double. We now introduce an algebra
$L({\cal A})$, called the {\it lattice algebra} of $\cal A$ by taking not just two but infinitely many
copies of $R({\cal A})$ (one for each site of an infinite lattice) and by imposing Heisenberg double-like
commutation relations between algebras at adjacent sites. More precisely, $L({\cal A})$
is, by definition, generated by symbols $Z^{(m)}_A$ with $A\in {\rm Ob}({\cal A})/{\rm Iso}$, $m\in {\bf Z}$
and $K_\alpha, \alpha\in {\cal K}_0({\cal A})$ which are subject to the following relations:
$$ K_\alpha K_\beta = K_{\alpha+\beta},\quad Z^{(m)}_A K_\alpha = (A|\alpha)^{(-1)^m}
K_\alpha Z^{(m)}_A,\leqno (3.1.1)$$
$$ Z^{(m)}_A Z^{(m)}_B = \langle B,A\rangle \sum_C g_{AB}^C Z^{(m)}_C, \leqno (3.1.2)$$
$$Z^{(m+1)}_A Z^{(m)}_B = \sum_{M,N} \gamma_{AB}^{MN} \langle \bar B-\bar M, \bar M -\bar N\rangle
 \cdot Z^{(m)}_M  Z^{m+1}_N K_{\bar B-\bar M}^{(-1)^m},
\leqno (3.1.3) $$
$$ Z^{(m)}_A  Z^{(n)}_B = (A|B)^{(-1)^{m-n}(n-m+1)} Z^{(n)}_B  Z^{(m)}_A , \quad |m-n|\geq 2.\leqno (3.1.4)$$

It is clear from these relations that the rule
$$\Sigma(Z^{(m)}_A) = Z^{(m+1)}_A, \quad \Sigma(K_\alpha)=K_\alpha^{-1}, \leqno (3.1.5)$$
defines an automorphism $\Sigma: L({\cal A})\to L({\cal A})$ which we call the shift
(or suspension) automorphism. 

\vskip .3cm

\noindent {\bf (3.2) Compatibility of the relations.} We now want to show that the relations
(3.1.1-4) are compatible in the sense that any element can be brought to a unique normal
form in which the upper indices of the $Z^{(m)}_A$ are increasing. 
More precisely, for any sequence $(a_i)_{i\in {\bf Z}}$ of elements of a possibly non-commutative
algebra $S$, almost all equal to 1, we define their ordered product to be
$$\prod_i^{\longrightarrow} a_i := a_pa_{p+1} ... a_q \in S,\leqno (3.2.1)$$
where $p,q$ are such thar $a_i=1$ unless $p\leq i\leq q$.

\proclaim (3.2.2) Proposition. The map of vector spaces
$$\nu: {\bf C}[{\cal K}_0({\cal A})]\otimes\bigotimes_{m\in {\bf Z}} R({\cal A}),
\quad K_\alpha\otimes\bigotimes_m [A_m] \mapsto \biggl( \prod_m^{\rightarrow} Z^{(m)}_{A_m}\biggr)
K_\alpha,$$
is an isomorphism.

Here and elsewhere in the paper all infinite tensor products of algebras are understood in the
restricted sense: almost all factors in any decomposable tensor are required to be 1. 
We will refer to an element explicitly realized as the value on $\nu$ on some tensor,
as being brought to the normal form.

\vskip .2cm

\noindent {\sl Proof:} The map $\nu$ is clearly surjective. To see the injectivity,
we need to establish the following two lemmas.

\proclaim (3.2.3) Lemma. For any $A,B,C\in {\cal A}$ and $m\in {\bf Z}$ the two possible ways of
bringing $Z^{(m+1)}_A Z^{(m)}_B Z^{(m-1)}_C$ to the normal form by using (3.1.3), lead to
the same answer.

\proclaim (3.2.4) Lemma. If $|m-n|\geq 2, |m+1-n|\geq 2$, then for any $A,B,C\in {\cal A}$
the multiplicative commutators with $Z^{(n)}_C$ of the left and the right hand sides of
(3.1.3) (the commutators being prescribed by (3.1.4)), are the same.

\vskip .2cm

\noindent {\sl Proof of Lemma 3.2.3:} The first way of bringing our monomial to the normal form is:
$$Z^{(m+1}_A Z^{(m)}_B Z^{(m-1)}_C = \sum_{M,N}
 \gamma_{AB}^{MN} \langle \bar B-\bar M, \bar M-\bar N\rangle Z^{(m)}_mZ^{(m+1)}_N K_{\bar B-\bar M}^{(-1)^{m+1}} Z^{(m-1)}_C = \leqno (3.2.5)$$
$$= \sum_{M,N}
 \gamma_{AB}^{MN} \langle \bar B-\bar M, \bar M-\bar N\rangle \cdot (\bar C| \bar B-\bar M)^{-1}
(\bar N, \bar C)^{-1} Z^{(m)}_M Z^{(m-1)}_C Z^{(m+1)}_N K_{\bar B-\bar M}^{(-1)^{m+1}}=$$
$$=\sum_{M,N,P,Q} \gamma_{AB}^{MN}\gamma_{MC}^{PQ} \langle \bar C-\bar P, \bar P-\bar Q\rangle
\cdot \langle \bar B-\bar M, \bar M-\bar N\rangle\cdot  (\bar C| \bar B-\bar M)^{-1}
(\bar N, \bar C)^{-1} (\bar C-\bar P|\bar N)\times$$
$$\times Z^{(m-1)}_P Z^{(m)}_Q Z^{(m+1)}_N K_{\bar C -\bar P}
^{(-1)^m} K_{\bar B-\bar M}^{(-1)^{m+1}},$$
while the other way is as follows:
$$Z^{(m+1)} Z^{(m)}_B Z^{(m-1)}_C = \sum_{P,U} \gamma_{PU}^{BC}
 \langle \bar C-\bar P, \bar P-\bar U\rangle \cdot Z^{(m+1)}_AZ^{(m-1)}_P
 Z^{(m)}_U K_{\bar C-\bar P}^{(-1)^m} =\leqno (3.2.6)$$
$$=\sum_{P.U} \sum_{P,U} \gamma_{PU}^{BC}
 \langle \bar C-\bar P, \bar P-\bar U\rangle (A|P)^{-1} Z^{(m-1)}_P Z^{(m+1)}_A Z^{(m)}_U 
K_{\bar C-\bar P}^{(-1)^m}=$$
$$=\sum_{P,U,Q,N} \gamma_{PU}^{BC}
\gamma_{AU}^{QN} \langle \bar C-\bar P, \bar P-\bar U\rangle (A|P)^{-1} \langle \bar U-\bar Q, \bar Q
-\bar N\rangle Z^{(m-1}_PZ^{(m)}_Q Z^{(m+1)}_N K_{\bar C-\bar P}^{(-1)^m}
 K_{\bar U-\bar Q}^{(-1)^{m+1}}.$$ 
Let us start to compare these two expressions. 
Our first remark is that the $\gamma$-quantities coincide.

\proclaim (3.2.7) Lemma. For any $A,B,C,P,Q,N\in {\cal A}$ we have the equality
$$\sum_M \gamma_{AB}^{MN}\gamma_{MC}^{PQ} = \sum_U \gamma_{BC}^{PU}\gamma_{AU}^{QN}.$$

\noindent {\sl Proof:} Both sides of the proposed equality have the following conceptual
meaning. They are equal to the orbifold number (modulo ${\rm Aut}(A)\times{\rm Aut}(B)
\times {\rm Aut}(C)$) of systems consisting, first, of a complex of length 2:
$$C\buildrel \psi\over\to B \buildrel \phi\over\to A, \quad \phi\psi=0,$$
and, second, of an identification of its cohomology, i.e., of isomorphisms
$$P\to {\rm Ker}(\psi), \quad Q\to {\rm Ker}(\phi)/{\rm Im}(\psi), \quad N\to
{\rm Coker}(\phi).$$
The object $M$ in the left hand side is ${\rm Im}(\psi)$, while $U$ in the right
hand side is ${\rm Im}(\phi)$. Q.E.D.

\vskip .2cm

Notice further that whenever a summand in any of the two sums in (3.2.7) is
non-zero, we have the equalities
$$\bar C-\bar P=\bar B-\bar U = \bar M-\bar Q, \quad \bar U-\bar Q = \bar A-\bar N =
\bar B -\bar M, \leqno (3.2.8)$$
which are obtained by applying  the fact that the Euler characteristic (in the Grothendieck group)
of a 4-term exact sequence is 0.

It follows from (3.2.8) that the ``$K$" factors in the end results of (3.2.5) and (3.2.6)
are the same. So it remains to compare the numerical factors given by the values of
the Euler form and its symmetrization.  We first compare the angle brackets, by noticing that
in virtue of (3.2.8),
$${ \langle \bar C-\bar P, \bar P-\bar U\rangle\cdot \langle \bar U-\bar Q, \bar Q-\bar N\rangle
\over \langle \bar C-\bar P, \bar P-\bar Q \rangle\cdot \langle \bar B-\bar M, \bar M-\bar N\rangle}
= { \langle \bar C-\bar P, \bar Q -\bar U\rangle\cdot \langle \bar U-\bar Q, \bar Q-\bar N\rangle
\over \langle \bar U-\bar Q, \bar M-\bar N\rangle}=$$
$$=\langle \bar C-\bar P, \bar Q-\bar U\rangle\cdot \langle \bar U-\bar Q, \bar Q-\bar M\rangle
=\langle \bar C-\bar P, \bar Q -\bar U\rangle \cdot \langle \bar U-\bar Q, \bar P-\bar C\rangle=$$
$$=  (\bar C-\bar P |\bar Q-\bar U)$$
and therefore the ratio of the total bracket contributions in the end results of (3.2.6)
and (3.2.5) is
$$ { (\bar C-\bar P|\bar Q-\bar U)\cdot (A|P)^{-1}\over (\bar C|\bar B-\bar M)^{-1}
 (\bar N|\bar C)^{-1} (\bar C-\bar P|\bar N)} = {(\bar C-\bar P|\bar Q-\bar U)\cdot 
(\bar C|\bar B-\bar M)\cdot
(P|N)\over (A|P)} =$$
$$(\bar C-\bar P|\bar Q-\bar U)\cdot (\bar C|\bar B-\bar M) \cdot (\bar P|\bar N-\bar A) =
(\bar C-\bar P|\bar Q-\bar U) \cdot (\bar C|\bar B-\bar M) \cdot (\bar P|\bar M-\bar B) =$$
$$= (\bar C-\bar P|\bar Q-\bar U) \cdot (\bar C-\bar P|\bar B-\bar M) = 1.$$
Lemma 3.2.3 is proved. 

\vskip .2cm

\noindent{\sl Proof of Lemma 3.2.4:} The product of $Z^{(m+1)}_A$ and $Z^{(m)}_B$ in any order
gives, when moved through $Z^{(n)}_C$, the factor
$$\bigl( (n-m)\bar A - (n-m+1)\bar B \bigl| \bar C\bigr) ^{(-1)^{m-n+1}}.$$
Thus it is enough to show that whenever $\gamma_{AB}^{MN}\neq 0$, one has the following
equality in ${\cal K}_0({\cal A})$:
$$(n-m)\bar A - (n-m+1)\bar B - (n-m)\bar N - (n-m+1)\bar M - (\bar B-\bar M),$$
where the last summand on the right comes from commuting $K_{\bar B-\bar M}^{(-1)^{m+1}}$
with $Z^{(n)}_C$. But, indeed, $\bar A-\bar B = \bar N-\bar M$, once $\gamma_{AB}^{MN}\neq 0$,
and the desired equality follows. 

This concludes the proof of Lemma 3.2.4 and Proposition 3.2.2. 

\vskip .3cm

\noindent {\bf (3.3) A basis in $L({\cal A})$.}
 It is  natural to label monomials in the $Z^{(i)}_A$
by (isomorphism classes of) graded objects of $\cal A$, 
i.e., by isomorphism classes of objects of $D^b({\cal A})$
which we represent as complexes with zero differential.

More precisely, to any graded object
$A^\bullet = \bigoplus A^{-i}[i]\in D^b({\cal A})$ we associate the monomial
$$Z(A^\bullet) = \prod_i^{\longrightarrow} {Z^{(i)}_{A^i} K_{A^i}^{(-1)^{i+1} i}
 \langle A^i, A^i\rangle ^i\over [A^i, A^{i-1}]}, \leqno (3.3.1)$$
where $[A^i, A^{i-1}]$ was defined in (2.2.1). The results of the previous subsection give: 

\proclaim (3.3.2) Proposition. The elements $Z(A^\bullet) K_\alpha $, for $A^\bullet
\in {\rm Ob}(D^b({\cal A}))/{\rm Iso}$ and $\alpha\in {\cal K}_0({\cal A})$, form a
{\bf C}-basis in $L({\cal A})$.

This, together with the homological interpretation of the quantities in Lemma (3.2.7),
suggests a deeper relation between $L({\cal A})$ and the derived category. More precisely, 
let $F: D^b({\cal A})\to D^b({\cal B})$ be any equivalence of triangulated categories
(we assume that both $\cal A$ and $\cal B$ satisfy our conditions of finiteness and homological dimension).
Then $F$ induces an isomorphism of Grothendieck groups $F_{\cal K}: {\cal K}_0({\cal A})\to
{\cal K}_0({\cal B})$ in a standard way. Our aim in this section is to prove the following result.
\vskip .1cm

\proclaim (3.4) Theorem. If $F$ is an equivalence of derived categories as above, then the
 correspondence
$$Z^{(p)}_A\to \Sigma^p(Z(F(A)), A\in {\cal A}, p\in {\bf Z}, \quad K_\alpha \to K_{F_{\cal K}(\alpha)},
\alpha\in {\cal K}_0({\cal A}),$$
defines an isomorphism of algebras $F_*: L({\cal A})\to L({\cal B})$.

\noindent {\sl Proof:} Our analysis is similar to that of (2.3). Namely, for $i\in {\bf Z}$
let ${\cal A}_i\i {\cal A}$ be the full subcategory of $A$ such that $F(A)\in {\cal B}[i]$. Then,
if $A_i\in {\cal A}_i, A_j\in {\cal A}_j$, we have
$${\rm Hom}_{\cal A}(A_i, A_j) = 0 \quad {\rm for}\quad j-i\notin \{0,1\}, \leqno (3.4.1a)$$
$${\rm Ext}^1_{\cal A}(A_i,A_j)=0 \quad {\rm for}\quad j-i\notin \{-1,0\}. \leqno (3.4.1b)$$
Also, each $A\in {\cal A}$ can be written uniquely as $A=\bigoplus A_i$ with $A_i\in {\cal A}_i$
and in the algebra $R({\cal A})$ we have the equality
$$[A] \quad =\quad \prod_i^{\longrightarrow} [A_i] \cdot |{\rm Hom}(A_i, A_{i-1})|^{-1/2}.
\leqno (3.4.2)$$
This means that the ordered product map defines an isomorphism of {\bf C}-vector spaces
$$\bigotimes_i R({\cal A}_i)\to R({\cal A}),$$
where the tensor product on the left is the restricted one. Denote by $R({\cal A}_i[-j])$
the subalgebra in $L({\cal A})$ spanned by $Z^{(j)}_A$ with $A\in {\cal A}_i$. 
Proposition 3.3.2 implies then the following.

\proclaim (3.4.3) Proposition. The map of vector spaces
$${\bf C}[{\cal K}_0({\cal A})] \otimes\bigotimes_{i,j\in {\bf Z}} R({\cal A}_i[-j])
\to L({\cal A}), \quad K_\alpha\otimes\bigotimes Z_{A_i}^{(j)} \mapsto 
\biggl(\prod_j^{\rightarrow}
\prod_i^{\rightarrow} Z_{A_i}^{(j)}\biggr)K_\alpha  ,$$
is an isomorphism.

 The proof of Theorem 3.4 consists basically of checking the relations and it is convenient to
first prove the following particular case, generalizing Theorem 2.3.

\proclaim (3.5) Proposition. For $F$ as above the rule $[A]\to Z(F(A))$ defines an injective
homomorphism of algebras $F_*: R({\cal A})\to L({\cal B})$. 

\noindent {\sl Proof:} 
 By an argument
similarly to Lemma 2.5.5, it is enough to prove the following partial statement.

\proclaim (3.5.1) Proposition. For any $i,j\in {\bf Z}$ and any $A'\in {\cal A}_{-i}$,
$A''\in {\cal A}_{-j}$, one has the equality $F_*([A']*[A'']) = F_*([A']) F_*([A''])$.

The minus sign is chosen for convenience, because $A\in {\cal A}_i$ means that $F(A)$, as a complex,
is situated in degree $(-i)$.
\vskip .1cm

\noindent {\sl Proof of (3.5.1):} We will work out several cases, similarly to the proof
of Theorem 2.3.

\vskip .2cm

\noindent\underbar{Case 1:} $i=j$. 
Let $A', A''\in {\cal A}_{-i}$ and let $F(A')=B'[-i], F(A'')=B''[-i]$.
If $A$ is such that $g_{A'A''}^A\neq 0$, then $A\in {\cal A}_{-i}$ and, denoting $B=F(A)$, we have
$g_{A'A''}^A=g_{B'B''}^B$. Thus

$$F_*([A'])F_*([A'']) = Z_{B'}^{(i)} K_{B'}^{(-1)^{i+1} i} Z_{B''}^{(i)}
K_{B''}^{(-1)^{i+1} i} \langle B', B'\rangle ^i \cdot \langle B'', B''\rangle ^i =$$
$$=Z^{(i)}_{B'}Z^{(i)}_{B''} K_{\bar B' +\bar B''}^{(-1)^{i+1} i} 
\langle B', B'\rangle ^i \cdot \langle B'', B''\rangle ^i \cdot (B'|B'')^i =$$
$$=\sum_{B\in {\cal B}} \langle B'', B'\rangle g_{B'B''}^B Z^{(i)}_{B} K_B^{(-1)^{i+1} i} \langle B, B\rangle ^i =
F_*([A']*[A'']),$$
where we used the equality $\bar B = \bar B' + \bar B''$ holding each time when $g_{B'B''}^B \neq 0$. 

\vskip .1cm

\noindent\underbar{Case 2:} $j=i+1$. Let $F(A')=B'[-i], F(A'')=B''[-i-1]$.
We have
$$F_*([A'])F_*([A'']) = Z_{B'}^{(i)} K_{B'}^{(-1)^{i+1} i} Z_{B''}^{(i+1)} K_{B''}^{ (-1)^{i+2} (i+1)}
\cdot \langle B', B'\rangle ^i \cdot \langle B'', B''\rangle ^{i+1},$$
while
$$F_*([A']*[A'']) = |{\rm Hom}(A'', A')|^{1/2} F([A'\oplus A'']) =$$ 
$$= { |{\rm Hom}(A'', A')|^{1/2} \cdot \langle B', B'\rangle ^i \cdot \langle B'', B''\rangle ^{i+1} \over
|{\rm Hom}(B'', B')|^{1/2}\cdot |{\rm Ext}^1 (B'', B')|^{1/2} }
 Z^{(i)}_{B'} K_{B'}^{(-1)^{i+1} i} Z^{(i+1)}_{B''}
K_{B''}^{(-1)^{i+2} (i+1)},$$
which is exactly the same once we recall that
$${\rm Ext}^1(B'', B') = {\rm Hom}(A'', A'), \quad {\rm Hom}(B'', B') = {\rm Ext}^{-1}(A'', A')=0.$$

\vskip .1cm

\noindent\underbar{Case 3:} $j=i-1$ and $i$ is even. 
 Let $F(A')=B'[-i], F(A'')=B''[-i+1]$. Then
$$F_*([A'])F_*([A'']) = Z^{(i)}_{B'} K_{B'}^{-i}Z^{(i-1)}_{B''} K_{B''}^{i-1} \langle B', B'\rangle ^i
\langle B'', B''\rangle ^{i-1}=$$
$$Z^{(i)}_{B'}Z^{(i-1)}_{B''} K_{-i\bar B' + (i-1) \bar B''}  \langle B', B'\rangle ^i
\langle B'', B''\rangle ^{i-1} (B'|B'')^{-i}=$$
$$= Z^{(i)}_{B'}Z^{(i-1)}_{B''} K_{-i\bar B' + (i-1) \bar B''} 
\langle \bar B'-\bar B'', \bar B'-\bar B''\rangle ^i \langle B'', B''\rangle ^{-1}=$$
$$\sum_{M,N} \gamma_{B', B''}^{MN} { \langle \bar B''-\bar M, \bar M\rangle \langle N, B''\rangle
\langle \bar B'-\bar B'', \bar B'-\bar B''\rangle ^i\over
\langle B'', B''\rangle } Z^{(i-1}_M K_{\bar B''-\bar M} Z^{(i)}_N K_{-i\bar B' + (i-1) \bar B''}.$$
We can transform the fraction in this expression in the same way as in the proof of Case (0,1)
of Proposition 2.5.3, getting:
$$ \langle B''[-i+1], B'[-i]\rangle \sum_{M,N} g_{B'[-i], B''[-i+1]}^{N[-i]\oplus M[-i+1]}
{Z^{(i-1)}_M K_M^{-1}Z^{(i)}_N K_{i (\bar B''-\bar B')} \langle \bar B''-\bar B', 
\bar B''-\bar B'\rangle ^i \over
\langle M, M\rangle\cdot \langle N,M\rangle \cdot |{\rm Ext}^1(N,M)|}.$$
Notice that for any non-zero summand we have $\bar B''-\bar B' =\bar M-\bar N$. Therefore
the last expression equals
$$\langle B''[-i+1], B'[-i]\rangle \sum_{M,N} g_{B'[-i], B''[-i+1]}^{N[-i]\oplus M[-i+1]}
{Z^{(i-1)}_M K_M^{i-1} Z^{(i)}_N K_N^{-i} \langle \bar N-\bar M, \bar N-\bar M\rangle ^i (N|M)^i
\over \langle M, M\rangle \cdot \langle N,M\rangle \cdot |{\rm Ext}^1(N,M)|}=$$
$$=\langle B''[-i+1], B'[-i]\rangle \sum_{M,N} g_{B'[-i], B''[-i+1]}^{N[-i]\oplus M[-i+1]}
{ Z^{(i-1)}_M K_M^{i-1} Z^{(i)}_N K_N^{-i} \langle M, M\rangle^{i-1} \langle N,N\rangle^i
\over [N,M] }=$$
$$= \langle A'', A'\rangle \sum_{A\in {\cal A}} g_{A'A''}^A F_*([A]) = F_*([A']*[A'']).$$

\vskip .2cm

\noindent\underbar{Case 4:} $j=i-1$ and $i$ is odd. Keeping the conventions for $A', A'', B', B''$
the same as before, we have$$F_*([A']) F_*([A'']) = Z^{(i)}_{B'} K_{B'}^i Z_{B''}^{(i-1)}
 K_{B''}^{-i+1} \langle B', B'\rangle^i \langle B'', B''\rangle^{i-1}=$$
$$=Z^{(i)}_{B'} Z^{(i-1)}_{B''} K_{i\bar B' - (i-1)\bar B''} \langle B', B'\rangle ^i 
\langle B'', B''\rangle^{i-1} (B'|B'')^{-i}=$$
$$= Z_{B'}^{(i)} Z_{B''}^{(i-1)} K_{i\bar B'-(i-1)\bar B''} \langle \bar B'-\bar B'', \bar B'-\bar B''
\rangle ^i \langle B'', B''\rangle ^{-1}=$$
$$=\sum_{M,N} \gamma_{B', B''}^{MN} { \langle \bar B''-\bar M, \bar M\rangle \cdot 
\langle N, B''\rangle \cdot \langle \bar B'-\bar B'', \bar B'-\bar B''\rangle ^i \over
\langle B'', B''\rangle}
Z^{(i-1)}_M K_{\bar M-\bar B''} Z^{(i)}_N K_{i\bar B'-(i-1)\bar B''}=$$
$$= \langle B''[-i+1], B'[-i]\rangle \sum_{M,N} g_{B'[-i], B''[-i+1]}^{N[-i]\oplus M[-i+1]}
{Z^{(i-1)}_M K_M Z^{(i)}_N K_{i(\bar B'-\bar B'')} \langle \bar B'-\bar B'', \bar B'-\bar B''\rangle ^i 
\over \langle M,M\rangle\cdot \langle N,M\rangle \cdot |{\rm Ext}^1(N,M)|} =$$
$$= \langle B''[-i+1], B'[-i]\rangle \sum_{M,N} g_{B'[-i], B''[-i+1]}^{N[-i]\oplus M[-i+1]}
{ Z^{(i-1)}_M K_M^{-i+1} Z_N^{(i)} K_N^i \langle \bar N-\bar M, \bar N-\bar M\rangle ^i
(N|M)^i \over \langle M,M\rangle [N,M] }=$$
$$= \langle B''[-i+1], B'[-i]\rangle \sum_{M,N} g_{B'[-i], B''[-i+1]}^{N[-i]\oplus M[-i+1]}
{ Z^{(i-1)}_M K_M^{-i+1} Z_N^{(i)} K_N^i \langle M,M\rangle^{i-1} \langle N,N\rangle^i
\over [N,M] },$$
and the argument is finished as in Case 3.

\vskip .2cm

\noindent\underbar{Case 5:} $|i-j|\geq 2$. Let $F(A')=B'[-i], F(A'')=B''[-j]$.
Then
$${\rm Hom}(A', A'') = {\rm Ext}^1(A', A'') = {\rm Hom}(A'', A') = {\rm Ext}^1(A'', A')=0.
\leqno (3.5.2)$$
Thus $[A']*[A''] = [A'\oplus A'']$. On the other hand, (3.5.2) implies that $(A'|A'') =
(B'|B'')=0$ and therefore
$$F_*([A'])F_*([A'']) = Z^{(i)}_{B'} K_{B'}^{(-1)^{i+1} i} Z_{B''}^{(j)} K_{B''}^{(-1)^{j+1} j}
\langle B',B'\rangle ^i \cdot \langle B'', B''\rangle ^j =$$
$$= F_*([A'\oplus A'']) =
F_*([A']*[A'']).$$
Proposition 3.5 is proved.

\vskip .3cm

\noindent {\bf (3.6) End of the proof of Theorem 3.4.} Let us define a linear operator
$F_*: L({\cal A})\to L({\cal B})$ by postulating its values on generators to be as stated in
the theorem and extending it to products of generators by using the normal form
of Proposition 3.2.2. In other words, we put
$$F_*\left( \biggl( \prod_i^{\to} Z^{(m)}_{A_m} \biggr) K_\alpha\right) = 
\biggl( \prod_i^{\to} F_*( Z^{(m)}_{A_m}) \biggr) K_{F_{\cal K}(\alpha)}.$$
 It is clear that
$F_*$ is bijective, so we need only to prove that it is an algebra homomorphism, i.e., that
it preserves the relations (3.1.1-4) in $L({\cal A})$. For (3.1.1) it is clear.
The preservation of (3.1.2) is the content of Proposition 3.5. The condition
that $F_*$ preserves (3.1.3-4) can be stated after some change of notation
as follows:

\proclaim (3.6.1) Proposition. Let $A', A''\in {\cal A}$. Then
$$F_*(Z^{(p+1)}_{A'})F_*(Z^{(p)}_{A''})=$$
$$\sum_{M,N\in {\cal A}}\gamma_{A' A''}^{MN} 
\langle \bar A''-\bar M, \bar M\rangle\cdot \langle \bar N, \bar A''-\bar M\rangle\cdot
F_*(Z^{(p)}_M)K^{(-1)^p}_{\overline{F(A'')}-\overline{F(M)}} F_*(Z^{(p+1)}_N).$$

\proclaim (3.6.2) Proposition. For $A', A''\in {\cal A}$ and $q\leq p-2$ we have
$$F_*(Z^{(p)}_A) F_*(Z^{(q)}_B) = (A|B)^{(-1)^{p-q}(q-p+1)} F_*(Z^{(q)}_B)F_*(Z^{(p)}_A).$$

Before proceeding to prove these statements, let us note the following.

\proclaim (3.6.3) Lemma. To establish Propositions 3.6.1-2 in full generality, it is enough
to prove them under the assumption that $A'\in {\cal A}_{-i}$, $A''\in {\cal A}_{-j}$
for some $i,j$.

\noindent {\sl
 Proof of the lemma:} Indeed, the multiplication law in $L({\cal A})$
with respect to the basis $K_\alpha Z(A^\bullet)$ consists of bringing the product
of two basis vectors to the normal form using (3.1.1-4). So our propositions, together
with what has been already proved, just say that
$$F_*(Z^{(p)}_{A'}) F_*(Z^{(q)}_{A''}) = F_*(Z^{(p)}_{A'}Z^{(q)}_{A''}).
\leqno (3.6.4)$$
If this is known each time when $A'\in {\cal A}_{-i}$, $A''\in {\cal A}_{-j}$,
then Proposition 3.4.3 together with an obvious modification of Lemma 2.5.5
give that (3.6.4) is valid in general. 

\vskip .3cm

\noindent {\bf (3.7) Proof of Proposition 3.6.1 when $A'\in {\cal A}_{-i}$, $A''\in {\cal A}_{-j}$.}
Again, we have to consider several cases.

\vskip .1cm

\noindent \underbar{Case 1:} $j=i$. Let $F(A')=B'[-i], F(A'')=B''[-i]$.
 In this case $F$ establishes a bijection between exact sequences
$$0\to M\to A''\to A'\to N\to 0$$
in $\cal A$ and exact sequences
$$0\to C\to B''\to B'\to D\to 0$$
in $\cal B$, so that for $C,D$ corresponding to $M,N$ we have
 $\gamma_{A' A''}^{MN}=\gamma_{B'B''}^{CD}$. From this the statement follows rather directly, by
comparing the normal form of $F_*(Z^{(p+1)}_{A'})F_*(Z^{(p)}_{A''})$ with the
image under $F_*$ of the normal form of $Z^{(p+1)}_{A'}Z^{(p)}_{A''}$.

\vskip .1cm

\noindent \underbar{Case 2:} $j=i+1$. Let $F(A') = B'[-i], F(A'')=B''[-i-1]$. Then
$$F_*(Z^{(p+1)}_{A'} = Z^{(i+p+1)}_{B'} K_{B'}^{(-1)^{i+p} i} \langle B', B'\rangle ^i,
\quad
F_*(Z^{(p)}_{A''} = Z^{(i+p+1)}_{B''} K_{B''}^{(-1)^{i+p} (i+1)} \langle B'', B''\rangle ^{i+1}.$$
Let $F^{-1}: D^b({\cal B})\to D^b({\cal A})$ be an inverse equivalence, and $(F^{-1})_*;
L({\cal B})\to L({\cal A})$ be the corresponding linear operators. The 
The argument needed is to handle our case is identical to the reasoning (already made in
the proof of Proposition 3.5) that $(F^{-1})_*$ realizes $R({\cal B})$ inside $L({\cal A})$. 

\vskip .1cm

\noindent \underbar{Case 3:} $j=i-1$. Let $F(A')=B'[-i], F(A'')=b]][-i+1]$. Then
${\rm Hom}(A'', A')=$ ${\rm Ext}^{-1}(B'', B')=0$, so $\gamma_{A'A''}^{MN}=0$ always except the case
$M=A'', N=A'$, in which case the value is 1. This means that
$$Z^{(p+1)}_{A'}Z^{(p)}_{A''} = Z^{(p)}_{A''} Z^{(p+1)}_{A'},$$
and we have to verify that $F_*$ preserved this commutativity. We have:
$$F_*(Z^{(p+1)}_{A'})= Z^{(p+i+1)}_{B'} K_{B'}^{(-1)^{p+i} i} \langle B', B'\rangle ^i, \quad
F_*(Z^{(p)}_{A''})= Z^{(p+i-1)}_{B''} K_{B''}^{(-1)^{p+i} (i-1)} \langle B'', B''\rangle ^{i-1},$$
and these expressions indeed commute because 
$$ Z^{(p+i+1)}_{B'}Z^{(p+i-1)}_{B''} = (B'|B'')^{-1} Z^{(p+i-1)}_{B''}Z^{(p+i+1)}_{B'}$$
and because of the commutation relation of the $K$'s with the $Z$'s.

\vskip .1cm

\noindent \underbar{Case 4:} $|j-i|\geq 2$. If $F(A')=B'[-i], F(A'')=B''[-j]$, then all the Hom and
Ext between $A'$ and $A''$ in either order are 0, so $(A'|A'')=(B'|B'')=1$,
and the argument is the same as in the previous case, only simpler because we do
not have to care about $(B'|B'')$. 

\vskip .3cm

\noindent {\bf (3.8) Proof of Proposition 3.6.2 when $A'\in {\cal A}_{-i}, A''\in {\cal A}_{-j}$.}
We consider several cases as to the relative position of $p+i$ and $q+j$. We denote
$F(A')=B'[-i], F(A'')=B''[-j]$.

\vskip .1cm

\noindent \underbar{Case 1:} $|(p+i)-(q+j)|\geq 2$. In this case
$$F_*(Z^{(p)}_{A'}) = Z^{(p+i)}_{B'} K_{B'}^{(-1)^{p+i+1} i} \langle B', B'\rangle ^i, \quad
F_*(Z^{(q)}_{A''}) = Z^{(q+j)}_{B''} K_{B'}^{(-1)^{q+j+1} j} \langle B'', B''\rangle ^j,$$
and by our assumption we have
$$F_*(Z^{(p)}_{A'})F_*(Z^{(q)}_{A''}) = (B'|B'')^\lambda F_*(Z^{(q)}_{A''})F_*(Z^{(p)}_{A'}),$$
where
$$\lambda = (-1)^{p+i-q-j} (q+j-p-1-1+i-j) = (-1)^{p+i-q-j}(q-p-1),$$
and once we take into account that $(B'|B'')=(A'|A'')^{(-1)^{i-j}}$,
we get the claimed statement.

\vskip .1cm

\noindent \underbar{Case 2:} $p+i=q+j$. This implies that $|i-j|\geq 2$ and thus there is neither Hom
nor Ext between $A', A''$ in either direction, hence $(A'|A'')=0$. So $Z^{(p)}_{A'}$
commutes with $Z^{(q)}_{A''}$. On the other hand, the vanishing of Hom and Ext implies that
$$Z^{(p+i)}_{B'} Z^{(p+i)}_{B''} = Z^{(p+i)}_{B''}Z^{(p+i)}_{B'} = Z^{(p+i)}_{B'\oplus B''},$$
so we are done.

\vskip .1cm

\noindent \underbar{Case 3:} $(p+i)-(q+j)=1$. Thus $i-j=q-p+1$. Since $|q-p|\geq 2$,
we will have $|i-j|\geq 2$ always except the case $q-p=-2$, when $i-j=-1$. In any event,
${\rm Hom}_{\cal B}(B'', B') = {\rm Ext}^{i-j}_{\cal A}(A'', A')=0$.
On the other hand,
$$F_*(Z^{(p)}_{A'}) = Z^{(p+i)}_{B'} K_{B'}^{(-1)^{p+i+1} i} \langle B', B'\rangle ^i, \quad
F_*(Z^{(q)}_{A''}) = Z^{(p+i-1)}_{B''} K_{B''}^{(-1)^{p+i} j} \langle B'', B''\rangle ^j,$$
and bringing their product to the normal form involves quantities $\gamma_{B'B''}^{MN}$
which vanish unless $B'=N, B''=M$. So $Z^{(p+i)}_{B'}$ and $Z^{(p+i-1)}_{B''}$
commute, and therefore
$$F_*(Z^{(p)}_{A'})F_*(Z^{(q)}_{A''})  = (B'|B'')^{-i+j} F_*(Z^{(q)}_{A''})F_*(Z^{(p)}_{A'}),$$
which is exactly what  we need, once we recall that $i-j=q-p+1$ and thus
$$(B'|B'')^{j-i}=(A'|A'')^{(-1)^{i-j} (j-i)}=
(A'|A'')^{(-1)^{q-p+1}(p-q-1)} =
(A'|A'')^{(-1)^{q-p}(q-p+1)}.$$

\vskip .1cm

\noindent \underbar{Case 4:} $(p+i)-(q+j)=-1$ is treated in an similar way. Theorem 3.4
is completely proved.

\vfill\eject

\centerline {\bf \S 4. Examples and discussion.}

\vskip 1cm

\noindent {\bf (4.1) Example: the ``universal cover" of the quantum group.}
Let us illustrate the construction of the lattice algebra on the classical example
of representations of quivers [R1-3].
 Let $\cal G$ be a semisimple
simply laced complex Lie algebra and $\Gamma$ be its Dynkin graph. Thus vertices of $\Gamma$
are identified with the  simple roots of $\cal G$ and for two such vertices $i\neq j$ the entry $a_{ij}$
of the Cartan matrix of $\cal G$ is minus the number of edges joining $i$ and $j$. Suppose
that an orientation of $\Gamma$ is chosen and let ${\cal A}={\rm Rep}_{{\bf F}_q}(\Gamma)$ be the
category of representations of $\Gamma$ over the finite field ${\bf F}_q$. Recall that such a representation $V$
is a rule which associates to any vertex $i$ a finite-dimensional ${\bf F}_q$-vector space
$V_i$ and to any edge $i\buildrel e\over\rightarrow j$ a linear operator
$V_e: V_i\to V_j$. As shown by Ringel ({\it loc. cit.}),
 the algebra $B({\cal A})$ is in this case isomorphic
to  $U_q({\bf b}^+)$, a natural ``Borel" subalgebra in $U_q({\cal G})$, the quantized enveloping algebra of $\cal G$.
More precisely, $U_q({\cal G})$ is generated by the symbols $E_i^\pm, K_i^{\pm 1}$ for $i\in {\rm Vert}(\Gamma)$
subject to the relations:
$$ E_i^\pm K_j = q^{\pm a_{ij}} K_j E_i^\pm, \quad K_iK_j=K_jK_i,\leqno (4.1.1)$$
$$\sum_{\nu=0}^{1-a_{ij}}
 {1-a_{ij}\choose \nu}_q (E_i^\pm)^\nu E_j^{\pm} (E_i^\pm)^{1-a_{ij}-\nu}=0, \quad i\neq j\leqno (4.1.2)$$
$$ [E^+_i, E^-_j] = {\delta_{ij} (K_i-K_i^{-1})\over q-1},\leqno (4.1.3)$$
and $B({\cal A})$ is isomorphic to the subalgebra generated by the $E^+_i$ and $K_i^{\pm 1}$.
Explicitly, $E_i^+$ corresponds to the element $[V(i)]$ where $V(i)$ is the respesentation associating ${\bf F}_q$
to the $i$th vertex and 0 to all other vertices. Similarly, $K_i$ corresponds to the element $K_{V(i)}$
of $B({\cal A})$. 

From this and the form of the comultiplication in $U_q({\bf b}^+)$ it is easy to deduce the following fact.

\proclaim (4.1.4) Proposition. For ${\cal A}={\rm Rep}_{{\bf F}_q}(\Gamma)$ the algebra $L({\cal A})$
is generated by the symbols $Z^{(m)}_i$, $m\in {\bf Z}, i\in {\rm Vert}(\Gamma)$ and $K_i^{\pm 1}$, 
$i\in {\rm Vert}(\Gamma)$ subject only to the following relations:
$$Z_i^{(m)}K_j = q^{(-1)^m a_{ij}}K_j Z_i^{(m)}, \quad K_iK_j = K_jK_i,\leqno (4.1.5)$$
$$\sum_{\nu=0}^{1-a_{ij}}
 {1-a_{ij}\choose \nu}_q (Z_i^{(m)})^\nu Z_j^{(m)} (Z_i^{(m)})^{1-a_{ij}-\nu}=0,
 \quad i\neq j \in {\rm Vert}(\Gamma), m\in {\bf Z},\leqno (4.1.6)$$
$$[Z^{(m)}_i, Z^{(m-1)}_j] = { \delta_{ij} K_i^{(-1)^m}\over q-1},\leqno (4.1.7)$$
$$Z^{(m)}_i  Z^{(n)}_j =  q^{(-1)^{m-n} (n-m+1)a_{ij}} Z^{(n)}_jZ^{(m)}_i \quad \forall i,j, \quad {\rm if}\quad |m-n|\geq 2.\leqno (4.1.8)$$

This shows that $L({\cal A})$ can be viewed as the ``universal cover", or the {\bf Z}-periodic version
of $U_q({\cal G})$. Note that the right hand side of the relation (4.1.7) is just one summand of the right hand
side of the similar relation (4.1.3): informally, in $L({\cal A})$ the other summand  is still present
but pertains to a different pair of generators. Let us also note the similarity of $L({\cal A})$
with so-called ``lattice Kac-Moody algebras" of [AFS].

\vskip .3cm

\noindent {\bf (4.2) ``Naive" lattice algebras.} The compatibility of the relations in
$L({\cal A})$ is not quite obvious, in particular because of the oscillator-type relations
(3.1.4) between copies of $R({\cal A})$ associated to non-adjacent lattice sites. 
So it may be useful to compare $L({\cal A})$ with the following general construction
which produces algebras in which such compatibility holds for free.

\vskip .1cm

Let $\Xi_m, m\in {\bf Z}$, be Hopf algebras and $\phi_m: \Xi_m\times \Xi_{m+1}\to {\bf C}$
be Hopf pairings. Define the {\it naive lattice algebra} $N=N(\{\Xi_m, \phi_m\})$
to be generated by elements of all the algebras $\Xi_m$ so that inside each $\Xi_m$
the elements are multiplied according to the multiplication law there while for
 elements of different
algebras we impose the relations:
$$\xi_{m+1}\xi_m = ({\rm Id}\otimes \phi_m \otimes {\rm Id})(\Delta_{\Xi_m}(\xi_m)
\otimes \Delta_{\Xi_{m+1}}(\xi_{m+1})), \leqno (4.2.1)$$
$$\xi_m\xi_{m'}=\xi_{m'}\xi_m, \quad |m-m'|\geq 2.$$
Thus if we put $\Xi_m$ at the $m$th site of a lattice, then the adjacent algebras
form a Heisenberg double while non-adjacent algebras commute. 

\proclaim (4.2.3) Proposition. The ordered product map $\bigotimes_{m\in {\bf Z}}\Xi_m\to N$
is always an isomorphism.

\noindent {\sl Proof:} We need only to verify that the two ways of bringing
any element   $\xi_{m+1}\xi_m\xi_{m-1}$ to the normal form $\sum_i \xi_{m-1}^{(i)}\xi_m^{(i)}
\xi_{m+1}^{(i)}$ by using (4.2.1), give the same result. But this easily follows
from the coassociativity of $\Delta_{\Xi_m}$.

\vskip .2cm

Applying this construction ot the case when for each $m$ we take $\Xi_m=B({\cal A})$ and $\phi_m
=\phi$ to be the Hopf pairing of (1.4), we get an algebra $N({\cal A})$ similar to $L({\cal A})$.
It is generated by symbols $Y^{(m)}_A, m\in {\bf Z}$, $A\in {\cal A}$ as well as $K^{(m)}_\alpha$,
$m\in{\bf Z}, \alpha\in {\cal K}_0({\cal A})$ with relations which are easy to find from
Proposition 1.5.3. In particular, for each $m$ the $K^{(m)}_\alpha$ form a copy of
${\bf C}[ {\cal K}_0({\cal A})]$, but adjacent copies do not commute. Because of this,
$N({\cal A})$ is not invariant under derived equivalence.

\vskip .3cm

\noindent {\bf (4.3) The bracket-free algebra $F({\cal A})$.} Most of the trouble in dealing with
the algebra $L({\cal A})$ comes from manipulating products of brackets, i.e., of values of the
Euler form. So it is tempting to define another, simpler algebra, by just dropping all these brackets.
More precisely, let $F({\cal A})$ be the algebra generated by synbols $X^{(m)}_A$,
$A\in {\cal A}, m\in {\bf Z}$ subject to the following relations:
$$X^{(m)}_A X^{(m)}_B = \sum_C g_{AB}^C X^{(m)}_C, \leqno (4.3.1)$$
$$X^{(m+1)}_A X^{(m)}_B =\sum_{M,N} \gamma_{AB}^{MN} X^{(m)}_M X^{(m+1)}_N, \leqno (4.3.2)$$
$$X^{(m)}_A X^{(n)}_B = X^{(n)}_BX^{(m)}_A, \quad |m-n|\geq 2.\leqno (4.3.3)$$
Note that these relations are compatible, as it follows from Lemma 3.2.7. In other words, 
the elements
$$X(A^\bullet) = \prod_m^{\to} X^{(m)}_{A^m}, \quad A^\bullet = \bigoplus_m A^m[-m], A^m\in {\cal A}$$
form a basis in $F({\cal A})$. The multiplication law in this basis can be described very nicely:
$$X(A^\bullet)X(B^\bullet) = \sum_{C^\bullet} \gamma_{A^\bullet B^\bullet}^{C^\bullet}
X(C^\bullet),\leqno (4.3.4)$$
where $\gamma_{A^\bullet B^\bullet}^{C^\bullet}$ is the orbifold number of long exact sequences
(2.4.3). Also, the procedure of bringing a maximally non-normal product to the normal form
can be nicely described in homological terms:
$$X^{(n)}_{A^{-n}} X^{(n-1)}_{A^{-n+1}} ... X^{(0)}_{A^0} = \sum_{d: A^\bullet\to A^\bullet[1]\atop
d^2=0} X(H^\bullet_d(A^\bullet)) \prod_m { |{\rm Aut}(H^m_d(A^\bullet))|\over
|{\rm Aut}(A^m)|}, \leqno (4.3.5)$$
where the sum is over all differentials making $A^\bullet$ into a complex, and 
$H^m_d(A^\bullet)$ is the $m$th cohomology with respect to $d$.

\vskip .1cm

Analogs of (4.3.4-5) can be obtained for the algebra $L({\cal A})$ as well, but they will be
encumbered by a lot of extra factors. However, these factors seem necessary to ensure the
invariance of the algebra under derived equivalence. In fact, the reason why $L({\cal A})$
possesses such invariance, is a subtle matching of two discrepancies. First is the discrepancy between 
the number
$g_{A^\bullet B^\bullet}^{C^\bullet}$
 of exact triangles in $D^{[-1,0]}({\cal A})$ and the number $\gamma_{A^\bullet B^\bullet}^{C^\bullet}$
of corresponding long exact sequences, which (for a particular case) 
was determined in Proposition 2.4.3 to be the factor
$|{\rm Ext}^1(N,M)|$.  The second is the discrepancy between the Hall multiplication 
$\circ$ in $H({\cal A})$ and its modification $*$ obtained by multiplying with $\langle B,A\rangle$,
see (1.3.4). When one term of a short exact sequence in $\cal A$ becomes shifted by 1
under a derived equivalence (so that we get an exact triangle in $D^{[-1,0]}({\cal A})$),
the difference between $*$ and $\circ$ will correspond
to the difference between $g$ and $\gamma$.

\vfill\eject

\centerline{\bf References.}

\vskip 1cm

\noindent [AF] A. Yu. Alexeev, L.D. Faddeev, Quantum $T^*G$ as 
a toy model for conformal field theory, {\it Comm. math. Phys.}, 
{\bf  141} (1991), 413-443.

\vskip .3cm

\noindent [AFS]  A. Yu. Alexeev, L.D. Faddeev, M. A. Semenov-Tian-Shansky, Hidden quantum
groups inside Kac-Moody algebras, Lecture Notes in Math. {\bf 1510}, p. 148-158,
Springer-Verlag 1992.

\vskip .3cm

\noindent [CP] V. Chari, A. Pressley, A Guide to Quantum groups,
Cambridge Univ. Press, 1995. 

\vskip .3cm

\noindent [GM] S.I. Gelfand, Y.I. Manin, Methods of Homological Algebra, Springer-Verlag, 1996.

\vskip .3cm

\noindent [Gr] J.A. Green, Hall algebras, hereditary algebras 
and quantum groups,
{\it Invent. Math.} {\bf  120} (1995), 361-377. 

\vskip .3cm

\noindent [Ha] D. Happel, Triangulated Categories in Representation Theory of Finite-Dimensional
Algebras (London Math. Soc. Lect. Note Series {\bf 119}), Cambridge Univ. Press, 1988.

\vskip .3cm

\noindent [J] A. Joseph, Quantum Groups and Their Primitive Ideals
(Ergebnisse der math. {\bf 29}), Springer-Verlag 1995.

\vskip .3cm

\noindent [Kap] M. Kapranov, Eisenstein series and quantum affine algebras, preprint alg-geom/9604018.

\vskip .3cm

\noindent [Kas] R.M. Kashaev, Heisenberg double and the pentagon relation, 
preprint q-alg/9503005.

\vskip .3cm

\noindent [Lu1] G. Lusztig, Introduction to Quantum Groups
(Progress in Math. 
{\bf  110}), Birkhauser, Boston, 1993. 

\vskip .3cm

\noindent [Lu2] G. Lusztig, Canonical bases arising from quantized enveloping 
algebras, {\it J. AMS}, {\bf  3} (1990), 447-498.

\vskip .3cm

\noindent [Lu3] G. Lusztig, Quivers, perverse sheaves and quantized 
enveloping algebras, {\it J. AMS}, {\bf  4} (1991), 365-421.

\vskip .3cm

\noindent [R1] C.M. Ringel, Hall algebras and quantum groups, 
{\it Invent. Math.}
{\bf  101} (1990), 583-592.

\vskip .3cm

\noindent [R2] C.M. Ringel, The composition algebra of a cyclic quiver,
{\it Proc. London Math. Soc.} (3) {\bf  66} (1993), 507-537.

\vskip .3cm

\noindent [R3] C.M. Ringel, Hall algebras revisited, in:
 ``Israel Math. Conf. Proc." vol. 7 (1993), p. 171-176. 

\vskip .3cm

\noindent [ST] M.A. Semenov-Tian-Shansky, Poisson Lie groups,
quantum duality principle and the quantum double,
{\it Contemporary math.}, {\bf  178}, p. 219-248,
Amer. Math. Soc, 1994.

\vskip .3cm

\noindent [X1] J. Xiao, Hall algebra in a root category, Preprint 95-070, Univ. of Bielefeld,
1995.

\vskip .3cm

\noindent [X2] J. Xiao, Drinfeld double and Green-Ringel theory of Hall algebras, Preprint
95-071, Univ. of Bielefeld, 1995.

\vskip 2cm

\noindent {\sl Author's address: Department of mathematics, Northwestern University, Evanston IL 
60208 USA \hfill\break
email: kapranov@math.nwu.edu}

\bye

\\
Title: Heisenberg doubles and derived categories
Author: M. Kapranov
Comments: 30 pages, plain TEX
Report-no: MPI-96-167
\\
Let A be an abelian category of finite type and homological dimension 1.
Then by results of Green R(A), the extended Hall-Ringel algebra of A,
has a natural Hopf algebra structure. We consider its Heisenberg double
Heis(A) and study its relation with D(A), the derived category of A.
We show that Heis(A) can be viewed as a "Hall algebra" of D^{0,1}(A),
the subcategory of complexes situated in degrees 0 and 1, in the following
sense: if B is the heart of a t-structure on D(A) lying in D^{0,1}(A),
then R(B) is naturally a subalgebra in Heis(A). Further, we define
a new algebra L(A) called the lattice algebra of A, obtained by
taking infinitely many copies of R(A), one for each site of
an infinite 1-dimensional lattice and imposing Heisenberg double-type
relations between copies at adjacent sites and oscillator-type
relations between copies at non-adjacent sites. This algebra
serves as the "Hall algebra" of the full derived category D(A)
in the following sense: any derived equivalence D(A)-->D(B)
induces an isomorphism of lattice algebras L(A)-->L(B).
//
\magnification=\magstep1
\baselineskip =13pt
\overfullrule =0pt

\centerline {\bf HEISENBERG DOUBLES AND DERIVED CATEGORIES}

\vskip .7cm

\centerline {\bf M. Kapranov}

\vskip 1cm
Let $\cal A$ be an Abelian category of finite homological dimension in which all ${\rm Ext}_{\cal A}^i(A,B)$ are finite sets.
One can, following C.M. Ringel [R1-3], associate to $\cal A$ an algebra $R({\cal A})$,
a version of the Hall algebra construction. Its structure constants are suitably normalized
numbers of short exact sequences.
Ringel has shown that in the case ${\cal A} = {\rm Rep}_{F_q}(\Gamma)$, the category of ${\bf F}_q$-representations
of a Dynkin quiver $\Gamma$, the algebra $R({\cal A})$ is identified with  the ``nilpotent" subalgebra $U_q({\bf n}^+)$
in the $q$-quantization $U_q({\bf g})$ of the semisimple Lie algebra corresponding to $\Gamma$.
This discovery has lead to several substantial advances in quantum group theory [Lu 1]. 
He has also shown how to put an algebra structure on the space $B({\cal A})
 = {\bf C}[{\cal K}_0 {\cal A}] \otimes R({\cal A})$
(here ${\cal K}_0 {\cal A}$ is the Grothendieck group of $\cal A$) so that for
${\cal A} = {\rm Rep}_{F_q}(\Gamma)$
one has $B({\cal A}) = U_q({\bf b}^+)$, the ``Borel" part of $U_q({\bf g})$.

The idea of extending the Hall algebra formalism to triangulated categories such as the
derived category $D^b({\cal A})$, seems to have been voiced by several people independently.
It appears naturally if one tries to find a construction of the full quantum group
$U_q({\bf g})$ in terms of ${\cal A} = {\rm Rep}_{F_q}(\Gamma)$. Indeed, various
 Abelian subcategories
in $D^b({\cal A})$ obtained from $\cal A$ by repeated application of derived
 Bernstein-Gelfand-Gelfand
reflection functors [GM] (and the Hall algebras of these subcategories) look temptingly
 similar to Borel subalgebras ${\bf b}^w\i {\bf g}$ obtained from ${\bf b}^+$ by action of
 elements of the Weyl group. 
Unfortunately, a direct mimicking of the Hall algebra construction but with exact
triangles replacing exact sequences, fails to give an associative multiplication,
even though the octohedral axiom looks like the right tool to establish the associativity.
One way to get around this difficulty is, as it was done in [X1], to ``amalgamate"
the (associative) Hall algebras of various Abelian subcategories in $D^b({\cal A})$,
but it seems to be not clear whether the resulting algebra is indeed $U_q({\bf g})$ nor
that it is a Hopf algebra at all. 

\vskip .1cm

The aim of the present paper is to exhibit an algebra $L({\cal A})$ which, although defined
in terms of $\cal A$, is invariant under derived equivalences and can be thus called
the ``Hall algebra of the derived category". We call $L({\cal A})$ the lattice algebra
of $\cal A$. Its construction was suggested by the fact that $U_q({\bf g})$ can be obtained
from the Hopf algebra  $U_q({\bf b}^+)$ by the Drinfeld double construction,
while the Hopf algebra structure on $B({\cal A})$ can be described in purely categorical terms,
as follows from the recent work of Green [Gr] (see [X2], [Kap]). However,
it turned out that it is not the Drinfeld double which appears naturally
in the study of $D^b({\cal A})$, but rather the so-called Heisenberg double of [AF] [ST].
In fact, one can find counterparts in Hopf algebra theory of several different
versions of derived categories, as shown in the following table:
\vfill\eject

\halign{\quad #  \hfill &\vrule width .8pt\vbox to
 .6cm{}\quad#{}\hfill\cr
\hskip .5cm{\bf Categories related to an }& \hskip .5cm {\bf Algebras related to the}\cr
\hskip .5cm {\bf Abelian category $\cal A$}& \hskip .5cm  {\bf Hopf algebra $\Xi = B({\cal A})$}\cr
\noalign{\hrule height .8pt}
$D^{[-1,0]}({\cal A})$, the category of {\bf Z}-graded& $HD(\Xi)$, the Heisenberg double\cr
 complexes situated in degrees $0, -1$& \cr
\noalign{\hrule}
 $D^b({\cal A})$, the standard bounded derived & \hskip .5cm $L({\cal A})$, the lattice algebra
\cr \hskip 2cm category &
\cr
\noalign{\hrule}
$D^{(2)}({\cal A})$, the category of 2-periodic & $DD(\Xi)$, the Drinfeld double
\cr \hskip 1cm (${\bf Z}/2$-graded) complexes & \cr
\noalign{\hrule}
 }

\vskip .3cm

The relation of $D^{[-1,0]}({\cal A})$ to the Heisenberg double of $\Xi=B({\cal A})$
is the easiest to understand: the commutation relations in $HD(\Xi)$ involve
certain products of the structure constants for the multiplication and comultiplication
in $\Xi$, and they can be interpreted as numbers of some 4-term exact sequences in $\cal A$,
which are obviously related to exact triangles in $D^{[-1,0]}({\cal A})$.

The algebra $L({\cal A})$ is obtained by taking one copy of
${\bf C}[{\cal K}_0({\cal A})]$ and infinitely many copies of $R({\cal A})$, one for each site of
an infinite 1-dimensional lattice and then imposing Heisenberg double-type commutation
relations between copies of $R({\cal A})$ at adjacent sites and oscillator relations of
the form $AB=\lambda_{AB} BA$, $\lambda_{AB}\in {\bf R}$, between basis vectors
of non-adjacent copies.  This algebra is similar to the ``lattice Kac-Moody algebras"
of [AFS]. The reason for taking an infinite lattice is clear: the $n$th copy of $R({\cal A})$,
$n\in {\bf Z}$, corresponds to the Abelian subcategory ${\cal A}[n]\i D^b({\cal A})$. 

\vskip .1cm

The author would like to acknowledge financial support from NSF grants and A.P. Sloan
Research Fellowship as well as from the Max-Planck Institute f\"ur Mathematik in Bonn
which provided excellent conditions for working on this paper.

\vfill\eject

\centerline {\bf \S 1. Heisenberg doubles for Hall-Ringel algebras}

\vskip 1cm

\noindent {\bf (1.1) The Heisenberg double.}
Let $\Xi$ be a Hopf algebra over {\bf C}.
We denote by 
$$\Delta: \Xi\rightarrow \Xi\otimes\Xi,  \quad
\epsilon: \Xi\rightarrow {\bf C}, \quad S: \Xi\rightarrow
\Xi$$
the comultiplication, the counit and the antipode of $\Xi$ respectively. 
For $x\in\Xi$ let $r_x: \Xi\rightarrow\Xi$ be the operator of right  multiplication by $x$, i.e., $r_x(y) = yx$, and let
$D_x: \Xi^*\rightarrow \Xi^*$ be the dual to $r_x$, i.e., for $f\in\Xi^*$ the functional $D_x(f)\in\Xi^*$ takes $y\mapsto f(yx)$.
Then the correspondence $x\mapsto D_x$ gives an embedding of algebras $D: \Xi\rightarrow {\rm End}(\Xi^*)$.

\vskip .1cm

For $f\in\Xi^*$ let $l_f: \Xi^*\rightarrow\Xi^*$ be the operator of left multiplication by $f$ (with respect to the algebra
structure on $\Xi^*$ defined by the map dual to $\Delta$), i.e., $l_f(\phi) = f\phi$. Again, we get an embedding of algebras
$l: \Xi^*\rightarrow {\rm End}(\Xi^*)$. 

\vskip .1cm

The Heisenberg double $HD(\Xi)$ is defined as the subalgebra in ${\rm End}(\Xi^*)$ generated by the images $D(\Xi)$ and $l(\Xi^*)$.
It is known [ST] that the map
$$\Xi^*\otimes \Xi\rightarrow HD(\Xi), \quad f\otimes x\mapsto l_fD_x,$$
is an isomorphism of vector spaces. Thus to describe the structure of $HD(\Xi)$ completely, it is enough to
explain how to bring a product $D_xl_f$ to a linear combination of products of the form $l_gD_y$. We will do this in the
coordinate-dependent language, following [Kas]. 

\vskip .2cm

Let $\{e_i\}, i\in I$, be a basis of $\Xi$ and $\{e^i\}$ be the dual
(topological) basis of $\Xi^*$. Introduce the structure constants for
the  multiplication and comultiplication with respect to our basis:
$$ e_{i} e_{j} = \sum_k m_{i j }^k e_k, \quad \Delta (e_k) = 
\sum \mu_{k}^{ij} e_{i}\otimes 
 e_{j}. \leqno (1.1.1)$$
Then one easily finds that
$$D_{e_i}(e^k) = \sum_j m_{ji}^k e^j, \quad l_{e^i}(e^k) = \sum_j \mu_j^{ik} e^j. \leqno (1.1.2)$$

\vskip .1cm

\proclaim (1.1.3) Proposition. In $HD(\Xi)$ we have the identity
$$D_{e_i} l_{e^j} = \sum_{a,b,c} m_{ab}^j \mu_i^{bc} l_{e^a} D_{e_c}.$$

\noindent {\sl Proof:} By (1.1.2), our statement is equivalent to:
$$\sum_k \mu_k^{jr} m_{si}^k  = \sum_{a,b,c,d} m_{ab}^j \mu_i^{bc} m_{dc}^r \mu_s^{ad}, \quad \forall i,j,r,s.$$
This equality, however, expresses the coincidence of the coefficients at $e_j\otimes e_r$ in $\Delta(e_se_i)$ and
$\Delta(e_s)\Delta(e_i)$, and so it is true.

\vskip .2cm

In view of this proposition we will regard $HD(\Xi)$ as an abstract algebra
 generated by symbols $Z_i, Z^i, i\in I$ subject to the
relations:
$$Z_iZ_j = \sum_k m_{ij}^k Z_k, \quad Z^iZ^j = \sum_k \mu_k^{ij} Z^k, \leqno (1.1.4)$$
$$Z_iZ^j = \sum_{a,b,c} m_{ab}^j \mu_i^{bc} Z^a Z_c. \leqno (1.1.5)$$

Note that $HD(\Xi)$ is not a Hopf algebra.

\vskip .3cm

\noindent{\bf (1.2) Heisenberg double with respect to a Hopf pairing.}
It is convenient to introduce a version of the above formalism
which avoids dualizing possibly infinite-dimensional spaces. More precisely,
let $\Xi$ and $\Omega$ be Hopf algebras. A Hopf pairing of $\Xi$ and $\Omega$
is a bilinear map $\phi: \Xi\times\Omega\to {\bf C}$ satisfying the following
conditions:
$$\phi(\xi\xi', \omega) = \phi^{\otimes 2}(\xi\otimes\xi', \Delta(\omega)),\leqno (1.2.1)$$
$$ \phi(\xi, \omega\omega')=\phi^{\otimes 2} (\Delta(\xi), \omega\otimes
\omega')\leqno (1.2.2)$$
$$\phi(1,\omega) = \epsilon_{\Omega}(\omega), \quad \phi(\xi, 1) = \epsilon_\Xi (\xi).\leqno (1.2.3)$$
Ths conditions (1.2.1-2) simply mean that the multiplication and the comultiplication
are conjugate with respect to the pairing. We did not include here any conditions on the antipodes since
we will not need them. 

\vskip .2cm

If $\phi$ is a Hopf pairing of $\Xi$ and $\Omega$, we define the Heisenberg double  $HD(\Xi, \Omega, \phi)$
associated to $\phi$ to be the tensor product $\Omega\otimes_{\bf C} \Xi$ with the multiplication
given as follows: First, both $\Omega$ (realized as $\Omega\otimes 1$) and $\Xi$ (realized as $1\otimes \Xi$)
are required to be subalgebras. Second, for $\xi\in\Xi, \omega\in\Omega$ we impose the condition:
$$\xi\omega = ({\rm Id}\otimes \phi\otimes {\rm Id})(\Delta_{\Omega}(\omega)\otimes \Delta_\Xi (\xi)),
\leqno (1.3.4)$$
where ${\rm Id}\otimes \phi\otimes {\rm Id}: \Omega\otimes\Omega\otimes\Xi\otimes\Xi\to
\Omega\otimes\Xi$ is the map induced by the transposed pairing $\phi^{op}: \Omega\otimes\Xi\to {\bf C}$
on the second and third factors.

Thus, when $\Omega=\Xi^*$ and $\phi$ is the canonical pairing, we get the definition of (1.1). 

\vskip .3cm

\noindent {\bf (1.3) Hall and Ringel algebras.} We now describe a particular class of Hopf algebras, whose
Heisenberg doubles we will be interested in.

\vskip .1cm

\noindent Let $\cal A$ be an Abelian category. We will say that $\cal A$ has finite type, if
if, for any objects $A,B \in {\cal A}$ all the groups
$ {\rm Ext}^i_{\cal A}(A,B)$
have finite cardinality and are zero for almost all $i$. If $\cal A$ is of finite type, then, for any
objects $A,B,C\in {\cal A}$, the number of subobjects $A'\i A$ such that $A'\simeq A$ and $C/A'\simeq B$,
is finite. Denote this number $g_{AB}^C$.

\vskip .1cm

Let $H({\cal A})$ be the {\bf C}-vector space with basis $[A]$ parametrized by all the isomorphism
classes of objects $A\in{\cal A}$. The rule
$$[A]\circ [B] = \sum_C g_{AB}^C [C] \leqno (1.3.1)$$
makes $H({\cal A})$ into an associative algebra with unit $1=[0]$, see [R1-3]. This algebra
is called the Hall algebra of $\cal A$.

\vskip .1cm

Let ${\cal K}_0({\cal A})$ be the Grothendieck group of $\cal A$. For an object $A\in{\cal A}$ let $\bar A$
be its class in ${\cal K}_0({\cal A})$. The rule
$$\langle \bar A, \bar B\rangle = \sqrt{\prod_{i\geq 0} | {\rm Ext}^i_{\cal A}
 (A,B)|^{(-1)^i}}, \quad A,B\in {\cal A} \leqno (1.3.2)$$
extends uniquely (because of the behavior of Ext in exact sequences) to a bilinear form
${\cal K}_0({\cal A}) \otimes {\cal K}_0({\cal A})\to {\bf R}^*$ known as the Euler form. We will often write just
$\langle  A, B\rangle $ for $\langle \bar A, \bar B\rangle $, if $A,B$ are objects.
The symmetrization of the Euler form will be denoted by
$$(\alpha|\beta) = \langle \alpha,\beta\rangle \cdot \langle \beta,\alpha\rangle,
\quad \alpha, \beta\in {\cal K}_0({\cal A}). \leqno (1.3.3)$$
The twisted multiplication
$$[A]*[B] = \langle B, A\rangle \cdot [A]\circ [B] \leqno (1.3.4)$$
is still associative. We will denote $R({\cal A})$ and call the Ringel algebra of $\cal A$ the same
vector space as $H({\cal A})$ but with $*$ as multiplication.

\vskip .1cm

\noindent{\bf (1.3.5) Remark.} It was C.M. Ringel [R3] who first drew attention to the particular twist (1.3.4).
More generally, one can twist by any bilinear form on ${\cal K}_0({\cal A})$, and the associativity will
still be preserved. In fact, several such twists were used in earlier papers by G. Lusztig [Lu2-3],
without specially distinguishing the Euler form (1.3.2). 

\vskip .1cm

Let
${\bf C}[{\cal K}_0 {\cal A}]$
be the group algebra of
${\cal K}_0{\cal A}$,
with basis
$K_{\alpha}$,
$\alpha \in {\cal K}_0 {\cal A}$
and multiplication
$K_{\alpha} K_{\beta} = K_{\alpha +\beta}$.
Let us extend the algebra
$R({\cal A})$
by adding to it these symbols
$K_{\alpha}$
which we make commute with
$[A] \in R({\cal A})$
by the rule

$$[A]K_{\beta} = (\bar A|\beta) K_{\beta} [A]. \leqno (1.3.6)$$
Denote the resulting algebra
$B({\cal A})$.
So as a vector space
$B({\cal A}) \simeq {\bf C} [{\cal K}_0{\cal A}] \otimes_{\bf C} R({\cal A})$,
with
$K_{\alpha} \otimes [A] \mapsto K_{\alpha}A$
establishing the isomorphism. We will call $B({\cal A})$ the extended Ringel algebra of $\cal A$.

\vskip .2cm

Assume now that $\cal A$ satisfies two additional conditions. First, any object of $\cal A$ has only
finitely many subobjects. Second, ${\rm Ext}^i_{\cal A}(A,B)=0$ for all $A,B$ and all
$i>1$. We will state the second condition by saying that the homological dimension of $\cal A$ is less or
equal to 1 and write ${\rm hd}({\cal A})\leq 1$. The next statement follows from results of
Green [Gr], see [X1] [Kap] for a detailed deduction.

\proclaim (1.3.7) Theorem. $B({\cal A})$ is a Hopf algebra with respect to the
comultiplication given on generators by
$$\Delta(K_{\alpha}) =K_{\alpha} \otimes K_{\alpha},$$
$$\Delta([A]) = \sum_{A^{\prime} \subset A} \langle A/A^{\prime}, A^{\prime}\rangle  {| {\rm Aut} (A^{\prime})|\cdot |{\rm Aut} (A/A^{\prime})|\over |{\rm Aut} (A)|} [A^{\prime}] \otimes K_{A^{\prime}} [A/A^{\prime}],$$ 
the counit $\epsilon : B({\cal A}) \rightarrow {\bf C}$
given by
$$\epsilon (K_{\alpha} [A]) =   1, \quad {\rm if}\quad { A=0} \quad {\rm and}\quad 
0,  \quad {\rm if} \quad { A \ne 0} $$
and the antipode
$S : B(\cal A) \rightarrow B({\cal A})$
given by
$$S(K_{\alpha} [A]) = \sum^{\infty}_{n=1} (-1)^n \sum_{A_0 \subset \ldots \subset A_n=A} \prod^n_{i=1} \langle A_i/A_{i-1},A_{i-1}\rangle  {\prod^n_{j=0} |{\rm Aut} (A_j/A_{j-1})|
\over  |{\rm Aut} (A)|} \cdot $$
$$\cdot [A_0] * [A_1/A_0] * \ldots * [A_n/A_{n-1}] \cdot K^{-1}_{\alpha} K^{-1}_A $$
where 
$A_0 \subset \ldots \subset A_n =A$
runs over arbitrary chains of strict
$(A_i \ne A_{i+1})$
inclusions of length $n$.

\vskip .2cm

\noindent {\bf (1.4) The Hopf pairing on $B({\cal A})$.} The elements $K_\alpha [A]$ form a {\bf C}-basis
of $B({\cal A})$. Let us define a bilinear pairing $\phi: B({\cal A})\times B({\cal A})\to {\bf C}$ by
putting

$$\phi (K_{\alpha}[A], K_{\beta} [B]) = (\alpha|\beta) ([A], [B])  = 
{ {(\alpha|\beta) \delta_{[A],[B]}} \over {| {\rm Aut} (A)|}}. \leqno (1.4.1)$$

\proclaim (1.4.2) Proposition.
The pairing $\phi$ is a Hopf pairing on $B({\cal A})$.

\noindent {\sl Proof:}
We need to prove the equality (1.3.1) (the other equality (1.3.2) will then follow by symmetry).
In other words, we need to prove that
$$\phi (K_{\alpha}[A]K_{\beta}[B], K_{\gamma}[C]) = \phi^{\otimes 2} (K_{\alpha}[A]
 \otimes K_{\beta} [B], \Delta(K_{\gamma} [C])) \leqno (1.4.3)$$
 To prove this, notice that the left hand side is
$$ (\bar A|\beta) \phi(K_{\alpha+\beta}[A][B], K_{\gamma}[C]) =
 (\bar A|\beta)(\alpha +\beta|\gamma)\phi ([A] * [B], [C]) = $$
$$= (\bar A|\beta)(\alpha+\beta|\gamma) \sum_{C^{\prime} \subset C} 
\langle C/C^{\prime}, C^{\prime}\rangle { {| {\rm Aut}(C^{\prime})| 
\cdot | {\rm Aut}(C/C^{\prime})|}\over  {| {\rm Aut} (C)|}} \cdot $$
$$\cdot \phi (A,C^{\prime}) \cdot \phi (B,C/C^{\prime}), \leqno (1.4.4) $$
while the right hand side of (1.4.3) is

$$ \phi^{\otimes 2} \biggl(K_{\alpha} [A] \otimes K_{\beta} [B], \quad \sum_{C^{\prime}
 \subset C} \langle C/C^{\prime}, C^{\prime}\rangle { {| {\rm Aut}(C^{\prime})| 
\cdot | {\rm Aut}(C/C^{\prime})|}\over {| {\rm Aut} (C)|} }\cdot $$
$$\cdot K_{\gamma}[C^{\prime}] \otimes K_{\bar C+\gamma} 
[C^ {\prime\prime}]\biggl) = \leqno (1.4.5)$$
$$= \sum_{C^{\prime} \subset C} \langle C/C^{\prime}, C^{\prime}\rangle
 { {| {\rm Aut}(C^{\prime})| \cdot |{\rm Aut}(C/C^{\prime})|} \over  
{|{\rm Aut}(C)|}} (\alpha|\gamma)(\beta |\gamma)(\beta|\bar C^{\prime}) \cdot $$
$$\cdot \phi ([A], [C^{\prime}]) \cdot \phi ([B], [C/C^{\prime}]). $$
Notice now that in order that
$\phi (A,C^{\prime}) \ne 0$,
we should have 
$A \simeq C^{\prime}$,
and under this assumption the corresponding summands in (1.4.4) and 
(1.4.5) coincide. Proposition is proved.

\vskip .3cm

\noindent {\bf (1.5) The Heisenberg double of $B({\cal A})$.} Let ${\rm Heis}({\cal A})$
be the Heisenberg double \hfill\break 
$HD(B({\cal A}), B({\cal A}), \phi)$, see (1.2). We will denote its generators    as follows:
$$Z^+_A = 1\otimes [A], \quad Z^-_A = [A]\otimes 1,\quad K_\alpha = 1\otimes K_\alpha,
\quad K_\alpha^- = K_\alpha\otimes 1. \leqno (1.5.1)$$
Thus, the $Z^+_A$ together with the $K_\alpha$, form a copy of $B({\cal A})$ inside 
${\rm Heis}({\cal A})$, and the same for the $Z^-_A$ with the $K^-_\alpha$. To find the cross-relations
between the plus and minus generators more explicitly, let us introduce the following
notations. For any objects $A,B,M,N\in{\cal A}$ let ${\cal F}_{AB}^{MN}$ be the set of
all exact sequences
$$0\to M\to B\buildrel \psi\over\to A\to N\to 0
\leqno (1.5.2)$$
and by $\gamma_{AB}^{MN}$ the quotient $|{\cal F}_{AB}^{MN}|/|{\rm Aut}(A)|\cdot |{\rm Aut}(B)|$.
The following statement, with its proof, is an adaptation of Proposition 6.2.12 from [Kap]
which, however, used a more cumbersome approach.

\proclaim (1.5.3) Proposition. We have the following equalities in ${\rm Heis}({\cal A})$:
$$Z^+_AZ^-_B = \sum_{M,N} \langle \bar B -\bar M, \bar M\rangle\cdot \langle \bar N, \bar B-\bar M\rangle \cdot
\gamma_{AB}^{MN} Z^-_M K_{\bar B-\bar M} Z^+_N=\leqno (1.5.4)$$
$$= \sum_{M,N} \gamma_{AB}^{MN} 
\langle \bar B-\bar M, \bar M-\bar N\rangle \cdot Z^-_M Z^+_N K_{\bar B-\bar M},$$
$$Z^-_AK_\alpha = (A|\alpha)^{-1}K_\alpha Z^-_A, \quad Z^+_A K^-_\alpha = K_\alpha^- Z^+_A,
\quad K_\alpha^+ K_\beta^- = (\alpha|\beta) K_\beta^-K_\alpha^+.\leqno (1.5.5)$$

\noindent {\sl Proof:} By (1.3.4) and the definition of the $Z^\pm_A$, we have
$$Z^+_AZ^-_B = ({\rm Id}\otimes\phi\otimes {\rm Id}) (\Delta([B])\otimes\Delta([A])) =$$
$$= \sum_{M,I} \sum_{I',N} \langle I, M\rangle \cdot \langle N, I'\rangle g_{MI}^B \cdot g_{I'N}^A
{ |{\rm Aut}(M)|\cdot |{\rm Aut}(I)|\cdot |{\rm Aut}(I')|\cdot |{\rm Aut}(N)|\over
 |{\rm Aut}(B)|\cdot |{\rm Aut}(A)|} \times$$
$$\times \phi(K_M[I], [I']) \cdot [M]\otimes K_{I'}N.$$
Note that $\phi(K_M[I], [I']) = \delta_{[I], [I']}/|{\rm Aut}(I)|$, so we can just put $I'=I$. Further,
 for any three objects $A,B,C$ let ${\cal E}_{AB}^C$ be the
set of exact sequences
$$0\rightarrow A\rightarrow C\rightarrow B \rightarrow 0.$$
Thus
$$\gamma_{AB}^{MN} = { |{\cal F}_{AB}^{MN}|\over |{\rm Aut}(A)|\cdot
|{\rm Aut}(B)|}, \quad g_{AB}^C = {|{\cal E}_{AB}^{C}|\over |{\rm Aut}(A)|\cdot
|{\rm Aut}(B)|}.$$
Notice now that
$${\cal F}_{AB}^{MN} \quad = \quad \coprod_{I\in {\rm Ob}({\cal A})/{\rm Iso}}
\bigl( {\cal E}_{MI}^B \times {\cal E}_{IN}^A \bigr) \bigl/ {\rm Aut}(I),
\leqno (1.5.6)$$
with ${\rm Aut}(I)$ acting freely. This just means that any
 long exact sequence (1.5.2) can be split into two short sequences with $I =
 {\rm Im}(\psi)$. Let ${\cal F}_{AB}^{MN}(I)$ be the $I$th part of the disjoint union (1.5.6).
Then, by taking all the above equalities into account, we find:
$$Z^+_AZ^-_B = \sum_{M,I,N} \langle I, M\rangle \cdot \langle N, I\rangle { |{\cal F}_{AB}^{MN}(I)|\over
|{\rm Aut}(A)|\cdot |{\rm Aut}(B)|} Z^+_MK_IZ^-_N,$$
and to get the claimed equality (1.5.4), it remains to notice that $\bar I = \bar B -\bar M$ once 
 ${\cal F}_{AB}^{MN}(I)\neq\emptyset$. The equalities (1.5.5) are obtained in a straightforward way.
Proposition is proved.

\vfill\eject

\centerline {\bf \S 2. Heisenberg doubles and tilting.}

\vskip 1cm

\noindent {\bf (2.1) Generalities on derived categories.} Let $\cal A$ be an Abelian category. 
By $C^b({\cal A})$
we denote the category of bounded complexes $A^\bullet = (A^i, d_i = d_{i, A}: A^i\rightarrow A^{i+1})$ over $\cal A$.
The shifted complex $A^\bullet[n], n\in {\bf Z}$, is defined by $(A^\bullet [n])^i = A^{n+i}$, $d_{i, A[n]} = 
(-1)^n d_{i,A}$. The homology objects of a complex $A^\bullet$ are denoted by $H^i(A^\bullet) = {\rm Ker}(d_i)/
{\rm Im}(d_{i-1})$. We denote by $D^b({\cal A})$ the bounded derived category of $\cal A$. It is obtained from
$C^b({\cal A})$ by formally inverting quasi-isomorphisms. If $A,B$ are two objects of $\cal A$ (regarded as complexes concentrated in degree 0), then
$${\rm Hom}_{D^b({\cal A})} (A, B[i]) = {\rm Ext}^i_{\cal A}(A,B).\leqno (2.1.1)$$

For any triangulated category $\cal D$ and any morphism $f: X\rightarrow Y$ in $\cal D$ we will denote by
${\rm Cone}(f)$ the isomorphism class of third terms $Z$ of possible exact triangles
$$ X\buildrel f\over\longrightarrow Y\rightarrow Z \rightarrow X[1].$$

Let us say that $\cal A$ has homological dimension $d$ and write ${\rm hd}({\cal A})\leq d$ if ${\rm Ext}^i_{\cal A}(A,B) = 0$
for any $A,B\in {\cal A}$ and any $i>d$. We say that $\cal A$ has finite homological dimension (${\rm hd}({\cal A}) < \infty$)
if ${\rm hd}({\cal A})\leq d$ for some $d$.

\proclaim (2.1.2) Proposition. If ${\rm hd}({\cal A})\leq 1$, then
 each  object of $D^b({\cal A})$ 
is isomorphic
to the complex $H^\bullet(A^\bullet)$ formed by the cohomology of $A^\bullet$ and equipped with zero differential.

This proposition is interesting for us because it gives a very explicit description of $D^b({\cal A})$ as a category. Indeed, given any two complexes $A^\bullet, B^\bullet$ with zero differential,
 we have $A^\bullet = \bigoplus_{i\in {\bf Z}} A^{-i}[i]$ and similarly for
$B^\bullet$, so by (2.1.1)
$${\rm Hom}_{D^b({\cal A})}(A^\bullet, B^\bullet) = \bigoplus_i {\rm Hom}_{\cal A}(A^i, B^i) \quad \oplus\quad
\bigoplus_i {\rm Ext}^1_{\cal A}(A^i, B^{i-1})\leqno (2.1.3)$$
In other words, a morphism $f: A^\bullet\rightarrow B^\bullet$ is the same as a sequence  of components
$f_i^{\rm Hom} \in {\rm Hom}_{\cal A}(A^i, B^i)$ and $f_i^{\rm Ext}\in {\rm Ext}^1_{\cal A}(A^i, B^{i-1})$. In the sequel
(including proof of (2.1.2) we will
use this notation for the components of a morphism.

\vskip .2cm

\noindent {\sl Proof of (2.1.2):} The proposition for $D^b({\cal A})$ seems to be well known. A simple argument (pointed out to me by A. Bondal) is to show by induction
that any bounded complex is quasiisomorphic
to the sum of its last cohomology object and
its canonical truncation just below this object. Here
we include,  for completeness sake,
a different proof which
does nt use induction and so is applicable not just to $D^b({\cal A})$ but also to other types of derived categories
(unbounded, periodic etc.).

\vskip .1cm

Let $(A^\bullet, d)$ be a complex, and let $K^\nu = {\rm ker}(d_\nu)$, $I^\nu = {\rm Im}(d_{\nu -1})$. We have short exact sequences
$$0\rightarrow K^\nu\buildrel \epsilon_\nu\over\rightarrow A^\nu\buildrel \pi_\nu\over\rightarrow I^{\nu+1}\rightarrow 0,
\leqno (2.1.4)$$
(with $\pi_\nu$ induced by $d_\nu$) which fit together into an exact sequence of complexes
$$0\rightarrow (K^\bullet, 0) \buildrel\epsilon\over\rightarrow (A^\bullet, d)\buildrel \pi\over\rightarrow (I^\bullet [1], 0)\rightarrow 0.
\leqno (2.1.5)$$
This sequence gives rise to an exact triangle in $D^b({\cal A})$, in particular, we get the boundary map
$\delta: (I^\bullet [1], 0) \rightarrow (K^\bullet [1], 0)$ such that $(A^\bullet, d) \simeq {\rm Cone} (\delta)[-1]$.
We want to compare (2.1.5) with the short exact sequence
$$0\rightarrow (I^\bullet, 0) \buildrel \phi\over\rightarrow (K^\bullet, 0)\buildrel\psi\over\rightarrow (H^\bullet, 0)\rightarrow 0
\leqno (2.1.6)$$
defining the cohomology $H^\bullet = H^\bullet(A^\bullet)$. This sequence implies that $H^\bullet \simeq {\rm Cone}(\phi)$.
Note that $\phi$ has only ${\rm Hom}$-components $\phi_\nu = \phi_\nu^{\rm Hom}: I^\nu \hookrightarrow K^\nu$.
Note also that $\phi_\nu = \delta[-1]_{\nu}^{\rm Hom}$. Indeed, the latter map is just the boundary homomorphsm
in the long cohomology sequence of (2.1.5), and this homomorphism is straightforwardly found to be $\phi_\nu$. 
As to $\delta[-1]^{\rm Ext}_\nu$, it is the class of the extension (2.1.4) and may well be non-zero. However,
the situation is saved by the following fact.

\proclaim (2.1.7) Lemma. There exists an automorphism $W$ of $(K^\bullet, 0)$ in the derived category such that $W\delta[-1] = \phi$.

The lemma implies our proposition by the axiom TR2 of triangulated categories (a commutative square extends to a morphism
of triangles). 

\vskip .2cm

\noindent {\sl Proof of the lemma:} Since ${\rm hd}({\cal A}) = 1$ and $\phi_\nu$ is injective, the restriction map
$$\phi_\nu^*: {\rm Ext}^1_{\cal A}(K^\nu, K^{\nu-1}) \rightarrow {\rm Ext}^1_{\cal A}(I^\nu, K^{\nu-1})$$
is surjective. Let $w_\nu \in {\rm Ext}^1_{\cal A}(K^\nu, K^{\nu-1})$ be any element mapping into $\delta[-1]_\nu^{\rm Ext}$.
Now define $W: (K^\bullet, 0)\rightarrow (K^\bullet, 0)$ to have the components $W_\nu^{\rm Hom} = {\rm Id}$
and $W_\nu^{\rm Ext} = - w_\nu$. Then $W$ is an isomorphism since it is given (with respect to the decomposition
$K^\bullet = \bigoplus K^{-i}[i]$) by a triangular matrix with identities on the diagonal. One immediately sees that $W\delta[-1]=\phi$,
since the Ext-terms in the composition will cancel. Lemma and Proposition 2.1.2 are proved.

\proclaim (2.1.8) Corollary. If ${\rm hd}({\cal A})\leq 1$, then each indecomposable object of
$D^b({\cal A})$ has the form $A[i]$ where $i\in{\bf Z}$ and $A$ is an indecomposable object of $\cal A$. 

\vskip .2cm

\noindent {\bf (2.2) The category $D^{[-1,0]}({\cal A})$ and tiltings.} Let $\cal A$ be as before. Denote by
$D^{[-1,0]}({\cal A})$ the full subcategory in $D^b({\cal A})$ formed by complexes situated in degrees
$-1, 0$ only. Given two Abelian categories $\cal A$ and $\cal B$, we will call an equivalence
$F: D^b({\cal A})\to D^b({\cal B})$ of triangulated categories a {\it tilting}, if
$F({\cal A})\i D^{[-1,0]}({\cal B})$. This condition is satisfied for equivalences given
by the so-called tilting modules [Ha]. 

From now on we assume that all the Abelian categories we consider have homological dimension less or equal to 1
and satisfy all the finiteness conditions of \S 1. Thus, for any such category $\cal A$ we have the Hopf
algebra $B({\cal A})$ and its Heisenberg double ${\rm Heis}({\cal A})$. 
For two objects $A,B\in{\cal A}$ denote
$$[A,B] = |{\rm Hom}(A,B)|^{+1/2} \cdot |{\rm Ext}^1 (A,B)|^{+1/2}. \leqno (2.2.1)$$
Because of two plus signs in the exponents, this quantity does not descend to the
Grothendieck group. 

Let us associate to any object (complex with zero differential)
$A^\bullet = A^{-1}[1]\oplus A^0\in D^{[-1,0]}({\cal A})$ the following element of
${\rm Heis}({\cal A})$:
$$Z(A^\bullet) = { Z^-_{A^{-1}}K^{-1}_{A^{-1}}Z^+_{A^0}\over \langle A^{-1}, A^{-1}\rangle
\cdot [A^0, A^{-1}]} = { Z^-_{A^{-1}}K^{-1}_{A^{-1}}Z^+_{A^0}\over \langle A^{-1}, A^{-1}\rangle
\cdot \langle A^0, A^{-1}\rangle\cdot |{\rm Ext}^1(A^0, A^{-1})|}. \leqno (2.2.2)$$

Now we can formulate the main result of this section.

\proclaim (2.3) Theorem. If $F: D^b({\cal A})\to D^b({\cal B})$ is a tilting, then
the correspondence $[A]\to Z(F(A))$ gives an injective homomorphism of algebras $F_*: R({\cal A})
\to {\rm Heis}({\cal B})$.

Before starting the proof, we do some preliminary work in the next subsection.

\vskip .3cm

\noindent{\bf (2.4) Counting exact triangles.} 
If $G$ is a finite group acting on a finite set $X$, we will call the ratio $|X|/|G|$ the orbifold number of elements of $X$
 modulo $G$.
This number is the same as $\sum_{ \{x\} \in X/G} 1/|{\rm Stab}(x)|$, the sum being over $G$-orbits on $X$, and $x$ being one representative
chosen for each orbit $\{x\}$.

Let $\cal A$ be as above.  For any three objects $A^\bullet, B^\bullet, C^\bullet$ of $D^b({\cal A})$
 we denote by 
$g_{A^\bullet, B^\bullet}^{C^\bullet}$ the orbifold number of exact triangles
$$A^\bullet\rightarrow C^\bullet\rightarrow B^\bullet \buildrel\partial\over\rightarrow A^\bullet[1]\leqno (2.4.1)$$
modulo ${\rm Aut}(A^\bullet)\times {\rm Aut}(B^\bullet)$.

If $A^\bullet, B^\bullet, C^\bullet$ are three bounded {\bf Z}-graded
objects of $\cal A$ (i.e., complexes with zero differentials), we denote by $\gamma_{A^\bullet, B^\bullet}^{C^\bullet}$
the orbifold number of long exact sequences 
$$...\rightarrow A^i \rightarrow C^i \rightarrow  B^i\buildrel
 \partial_i\over\rightarrow A^{i+1}\rightarrow
...\leqno (2.4.3)$$
modulo $\prod_i {\rm Aut}(A^i)\times {\rm Aut}(B^i)$.

If $A,B,C\in {\cal A}$ are three objects considered as {\bf Z}-graded objects (complexes with zero differential)
concentrated in degree 0, then clearly $\gamma_{AB}^C = g_{AB}^C$ coincides with the number introduced in
(1.3). For general complexes with zero differential, $g_{A^\bullet, B^\bullet}^{C^\bullet}$
differs from $\gamma_{A^\bullet, B^\bullet}^{C^\bullet}$. We will need one particular case when these
numbers can be easily compared.

\proclaim (2.4.3) Proposition. Let $A,B,M,N$ be any objects of $A$. Then
$$g_{A, B[1]}^{M[1]\oplus N} =\gamma_{A, B[1]}^{M[1]\oplus N}  \cdot |{\rm Ext}^1(N,M)|.$$

\noindent {\sl Proof:} The number $\gamma_{A, B[1]}^{M[1]\oplus N}$ counts exact sequences 
$$0\to M\buildrel u\over\to B\buildrel\varphi\over\to A\buildrel v\over\to N\to 0\leqno (2.4.4)$$
and thus is equal to
$${ \left| \bigl\{ \varphi: B\to A| \,\,  {\rm Ker}(\varphi)\simeq M, \, {\rm Coker}(\varphi)\simeq N\bigr\}
\right| \cdot |{\rm Aut}(M)|  \cdot |{\rm Aut}(N)|\over   |{\rm Aut}(A)| \cdot |{\rm Aut}(B)| }.\leqno (2.4.5)$$
The number $g_{A^\bullet, B^\bullet}^{C^\bullet}$ counts exact triangles
$$A\buildrel\alpha\over\to M[1]\oplus N \buildrel\beta\over\to B[1]\buildrel \varphi[1]\over\to 
A[1].\leqno (2.4.6)$$
Of course, any such triangle gives rise to a sequence of the form (2.4.4), but the correspondence is not
bijective. 

More precisely, let us fix $\varphi: B\to A$.  By (2.1.2), in order that a triangle (2.4.6) with the boundary map $\varphi[1]$
exists, it is necessary and sufficient that ${\rm Ker}(\varphi)\simeq M$, ${\rm Coker}(\varphi)\simeq N$.
If $\varphi$ satisfies this property, then the axiom TR2 of triangulated categories
implies that a triangle (2.4.6) can be constructed uniquely modulo an isomorphism of triangles
identical on $A$ and $B[1]$. Such an isomorphism is the same as just an automorphism of $M[1]\oplus N$.
So
$$g_{A, B[1]}^{M[1]\oplus N} = 
 { \left| \bigl\{ \varphi: B\to A| \,\,  {\rm Ker}(\varphi)\simeq M, \, {\rm Coker}(\varphi)\simeq N\bigr\}
\right|\cdot |{\rm Aut}(M[1]\oplus N)|\over |{\rm Aut}(A)|\cdot |{\rm Aut}(B)|\cdot |{\rm Stab}(\varphi)|},
\leqno (2.4.7)$$
where ${\rm Stab}(\varphi)\i {\rm Aut}(M[1]\oplus N)$ is the subgroup of automorphisms which,
together with the identities of $A$ and $B[1]$, give an automorphism of the triangle (2.4.6).
In fact, we claim that ${\rm Stab}(\varphi) = \{{\rm Id}\}$. To see this,
note that ${\rm Aut}(M[1]\oplus N)$ is the block matrix group
$$\pmatrix {{\rm Aut}(M)& {\rm Ext}^1(N,M)\cr 0& {\rm Aut}(N)}.\leqno (2.4.8)$$
Let $\psi=(\psi_M, \psi_{NM}, \psi_N)$ be an element of ${\rm Stab}(\varphi)$, so
$\psi_M\in {\rm Aut}(M)$, $\psi_{NM}\in {\rm Ext}^1(N,M)$, $\psi_N\in {\rm Aut}(N)$. 
Then $\psi_M={\rm Id}$ since $u: M\to B$ is an injection, and $\psi_N={\rm Id}$ since
$v: A\to N$ is a surjection and since $\psi$ together with ${\rm Id}_A, {\rm id}_B$ forms a morphism of triangles.
Further, the commutativity of the diagram
$$\matrix{&A &\buildrel\alpha\over\longrightarrow& M[1]\oplus N &\cr
{\rm Id}&\big\downarrow &&\big\downarrow & \psi=\pmatrix{1&\psi_{MN}\cr 0&1} \cr
&A&\buildrel\alpha\over\longrightarrow& M[1]\oplus N }$$
means that $\psi_{MN}\circ v=0$ in ${\rm Ext}^1(A,M)$. But because ${\rm hd}({\cal A})\leq 1$,
the surjection $v: A\to N$ induces an injection ${\rm Ext}^1(N,M)\to
{\rm Ext}^1(A,M)$, so $\psi_{MN}=0$. This proves our claim that ${\rm Stab}(\varphi)=\{{\rm Id}\}$.
The proposition follows now by comparing (2.4.7) with (2.4.5) and taking into account the factorization
(2.4.8). 

\vskip .3cm

\noindent {\bf (2.5) Proof of Theorem 2.3.} 
 First of all let ${\cal A}_i \i {\cal A}$ be the full subcategory
of $A$ such that $F(A)\in {\cal B}[i], i=0, 1$. Notice that for $A_i\in {\cal A}_i$
we have, denoting $B_i = F(A_i)[-i]\in {\cal B}$:
$${\rm Hom}_{\cal A}(A_{1}, A_0) = {\rm Hom}_{D^b({\cal B})}(B_{1}[1], B_0) = 0,\leqno (2.5.1a)$$
$$ {\rm Ext}^1_{\cal A}(A_0, A_{1}) = {\rm Hom}_{D^b({\cal B})}(A_0, A_{1}[2]) = 
{\rm Ext}^2_{\cal B}(A_0, A_{1})=0.\leqno (2.5.1b) $$
Let us denote by $F_*: R({\cal A})\to {\rm Heis}({\cal A})$ the unique {\bf C}-linear
map taking $[A]$ to $Z(F(A))$. Its injectivity is clear from the behavior on basis vectors,
so the main task is to prove that $F_*$ is an algebra homomorphism,
i.e., that
$$F_*([A']* [A'']) = F_*([A']) F_*([A'']), \quad \forall A', A''\in {\cal A}. \leqno (2.5.2) $$
The proof will be done in two steps. The first is given in the next proposition.

\proclaim (2.5.3) Proposition. The equality (2.5.2) holds in the case when $A'=A'_i\in {\cal A}_i$,
$A''=A''_j\in {\cal A}_j$ for some $i,j\in \{0,1\}$.

The second step is to deduce the general case from these particular cases. Let us first
explain how this is done and then prove Proposition 2.5.3. Namely, by (2.1.8) each indecomposable
object of $\cal A$ lies in one of the ${\cal A}_i$. So each $A\in {\cal A}$
can be uniquely written as $A=A_0\oplus A_{1}$ with $A_i\in {\cal A}_i$. By (2.5.1)
we have
$$[A] = |{\rm Hom}(A_0, A_1)|^{-1/2} [A_{1}]*[A_0].\leqno (2.5.4)$$
This means that we have a kind of normal form for elements of $R({\cal A})$, i.e., the
map $R({\cal A}_{1})\otimes_{\bf C} R({\cal A}_0)\to R({\cal A})$ given by the
multiplication, is an isomorphism of vector spaces.
 So the second step is accomplished by the next easy lemma.

\proclaim (2.5.5) Lemma. Let $R$ be a {\bf C}-algebra and $R_0, R_{1}\i R$ be subalgebras
such that the multiplication induced an isomorphism of vector spaces
$R_1\otimes_{\bf C}R_0\to R$. Let $S$ be another algebra
and $\phi: R\to S$ be a {\bf C}-linear map. Suppose that the equality $\phi(a'a'')
=\phi(a')\phi(a'')$ holds whenever $a'\in R_i, a''\in R_j$ for some $i,j$.
Then it holds for any $a',a''$, i.e., $\phi$ is an algebra homomorphism.

\noindent {\sl Proof of (2.5.5):} It is enough, by linearity, to consider the case when
$a'=a'_{1}a'_0, a''=a''_{1}a''_0$ with $a'_i, a''_i\in R_i$. Then, by our assumptions,
$$\phi(a')\phi(a'') = \phi(a'_{1})\phi(a'_0)\phi(a''_{1})\phi(a''_0).$$
Let us write $a'_0a''_{1} = \sum_j \alpha^{(j)}_{1}\alpha^{(j)}_0$ with $\alpha^{(j)}_i\in R_i$,
using our assumption on $R$.
Then, by our assumptions on $\phi$,
$$\phi(a'_0)\phi(a''_{1}) = \phi (a'_0)a''_{1}) = \sum_j \phi(\alpha^{(j)}_{1}) 
\phi(\alpha^{(j)}_0),$$
and so
$$\phi(a')\phi(a'') = \sum_j \phi(a'_{1}\alpha^{(j)}_{-1}) \phi(\alpha^{(j)}_0a''_0) = 
\sum \phi(a'_{1}\alpha^{(j)}_{1}\alpha^{(j)}_0a''_0) =\phi(a'_{1}a'_0a''_{1}a''_0)=
\phi(a'a''),$$
as claimed.

\vskip .2cm

We now prove Proposition 2.5.3 by considering all four possibilities for
$i,j$.

\vskip .2cm

\noindent\underbar {Case (0,0)}: $A'=A'_0, A''=A''_0\in {\cal A}_0$. Let $F(A')=B',
F(A'')=B''$. Then $F_*(A')=Z^+_{B'}$ and $F_*(A'')=Z^+_{B''}$. Since $F$ is an embedding
of an admissible Abelian category, it establishes a bijection between short exact sequences
$$0\to A'\to A\to A''\to 0$$
in $\cal A$ and short exact sequences
$$0\to B'\to B\to B''\to 0$$
in $\cal B$ (since both kinds of sequences are interpreted as exact triangles in the same
triangulated category $D^b({\cal B})$). Thus the equality (2.5.2) holds.

\vskip .2cm

\noindent\underbar {Case (1,1)}: $A'=A'_{1}, A''=A''_{1}\in {\cal A}_{1}$.
Let $F(A')=B'[1], F(A'')=B''[1]$. If
$$0\to A'\to A\to A''\to 0$$
is a short exact sequence in $\cal A$, then $A\in {\cal A}_1$, and denoting $B=F(A)[-1]$,
we have $g_{A'A''}^A=g_{B'B''}^B$. By our definition,
$F_*([A']) = Z^-_{B'}K_{B'}^{-1}  \langle B', B'\rangle ^{-1}$, and similarly for
$F_*([A''])$. Therefore
$$F_*([A']) F_*([A'']) = {Z^-_{B'}K_{B'}^{-1} Z^-_{B''} K_{B''}^{-1}\over \langle B',B'\rangle\cdot
\langle B'',B''\rangle} =  {Z^-_{B'}Z^-_{B''} K^{-1}_{\bar{B}' + \bar{B}''}\over 
\langle B',B'\rangle\cdot \langle B'',B''\rangle\cdot (B'|B'')} =$$
$$={ \sum_{B} g_{B'B''}^B \langle B'',B'\rangle Z^-_B K^{-1}_B\over \langle B ,B \rangle}
=F_*([A'] *[A'']).$$
Here we used the fact that $\bar B= \bar B' +  \bar B''$ whenever  $g_{B'B''}^B\neq 0$.

\vskip .2cm

\noindent\underbar {Case (1,0)}: $A'=A'_1\in {\cal A}_0, A''=A''_0\in {\cal A}_1$. Let
$F(A')=B'[1], F(A'')=B''$. We have, by definition of $F_*$
$$ F_*([A'])F_*([A'']) = Z^-_{B'} K^{-1}_{B'} Z^+_{B''} \cdot \langle B', B'\rangle ^{-1},$$
while, by (2.5.4), (2.5.1)  and by definition of $F_*$, 
$$F([A']*[A'']) = |{\rm Hom}(A'', A')|^{1/2} F([A'\oplus A'']) = 
{|{\rm Ext}^1(B'', B')|^{1/2} Z^-_{B'}K^{-1}_{B'}Z^+_{B''}\over
\langle B',B'\rangle \cdot |{\rm Hom}(B'', B')|^{1/2}\cdot |{\rm Ext}^1(B'', B')|^{1/2}},$$
which is the same as the previous quantity once we recall that 
${\rm Hom}(B'', B') = {\rm Ext}^{-1}(A'', A')=0$.

\vskip .2cm

\noindent \underbar{Case (0,1)}: $A'=A'_0\in {\cal A}_0, A''=A''_1\in{\cal A}_1$. Let
$F(A')=B', F(A'')=B''[1]$. We want to verify that 
$$ F_*([A']) F_*([A'']) = F_*([A']* [A'']) = \sum_{A\in {\cal A}} g_{A'A''}^A F_*([A]).$$
Note that the $A$ entering the last sum, may not lie in any of the ${\cal A}_i$. However,
if $A$ is included into an exact sequence
$$0\to A'\to A\to A''\to 0,$$
then $F(A)$, which necessarily has the form
 $M[1]\oplus N$ for some $M,N\in {\cal B}$, is included into an exact triangle
$$B'\to M[1]\oplus N \to B''[1] \to B'[1],$$
and this correspondence is a bijection, i.e., $g_{A',A''}^A = g_{B', B''[1]}^{M[1]\oplus N}$.
If such a triangle exists, then in ${\cal K}_0({\cal B})$ we have
$\bar M-\bar N =  \bar B''-\bar B'$, as it follows from the corresponding 4-term exact
sequence. Now, by  (2.2.2), we have 
$$F_*([A'])F_*([A'']) ={ Z^+_{B'} Z^-_{B''}K^{-1}_{B''}\over \langle B'', B''\rangle }=$$
$$=\sum_{M,N\in {\cal B}} \gamma_{B', B''[1]}^{M[1]\oplus N}
 { \langle \bar B''-\bar M, \bar M\rangle\cdot
\langle \bar N, \bar B''-\bar M\rangle\over \langle B'', B''\rangle}
Z^-_M K_{\bar B''-\bar M} Z^+_N K^{-1}_{B''}=$$
$$=\sum_{M,N} \gamma_{B', B''[1]}^{M[1]\oplus N}
{\langle B'', M\rangle \cdot \langle N, B''\rangle \over \langle B'', B''\rangle \cdot 
\langle M, M\rangle \cdot \langle N, M\rangle \cdot (B''|N)} Z^-_M K_M^{-1}Z^+_N.$$
Note that the quantity represented by the fraction in the last expression, can be written as
$${ \langle \bar B'', \bar M-\bar N  \rangle \over \langle B'', B''  \rangle 
\cdot \langle M, M \rangle  \cdot  \langle N, M  \rangle } = 
{\langle   \bar B'', \bar B''-\bar B' \rangle \over \langle B'', B''  \rangle 
\cdot \langle M, M \rangle  \cdot  \langle N, M  \rangle }.$$
Therefore, applying Proposition 2.4.3, we find:
$$F_*([A'])F_*([A'']) = \sum_{M,N} \gamma_{B', B''[1]}^{M[1]\oplus N}
{\langle   \bar B'', \bar B''-\bar B' \rangle \over \langle B'', B''  \rangle 
\cdot \langle M, M \rangle  \cdot  \langle N, M  \rangle } Z^-_M K_M^{-1}Z^+_N =$$
$$= \langle B''[1], B'\rangle \sum_{M,N} g_{B', B''[1]}^{M[1]\oplus N}
{Z^-_M K_M^{-1}Z^+_N \over \langle M, M\rangle\cdot \langle N, M\rangle \cdot
|{\rm Ext}^1({N,M})| } = $$
$$=\langle A'', A'\rangle \sum_{A\in {\cal A}} g_{A', A''}^A F_*([A]) =
F_*([A']*[A'']).$$
This finishes the proof of Theorem 2.3.

\vfill\eject

\centerline {\bf \S 3. The lattice algebra and the full derived category.}

\vskip 1cm

\noindent {\bf (3.1) Definition of the lattice algebra.} Let $\cal A$ be an Abelian category 
with ${\rm hd}({\cal A})\leq 1$ satisfying all the finiteness conditions of \S 1. Let
$B({\cal A})\supset R({\cal A}) $ be its extended Ringel algebra and ${\rm Heis}({\cal A})$
be the Heisenberg  double of $B({\cal A})$. As we saw in \S 2, ${\rm Heis}({\cal A})$ is
naturally related to the subcategory $D^{[-1,0]}({\cal A})$ in the derived category $D^b({\cal A})$.
Hereby the two copies of $\cal A$ inside $D^{[-1,0]}({\cal A})$, given by complexes concentrated in degree 
$0$ (resp. $(-1)$), give rise to two copies of $R({\cal A})$ in the double. We now introduce an algebra
$L({\cal A})$, called the {\it lattice algebra} of $\cal A$ by taking not just two but infinitely many
copies of $R({\cal A})$ (one for each site of an infinite lattice) and by imposing Heisenberg double-like
commutation relations between algebras at adjacent sites. More precisely, $L({\cal A})$
is, by definition, generated by symbols $Z^{(m)}_A$ with $A\in {\rm Ob}({\cal A})/{\rm Iso}$, $m\in {\bf Z}$
and $K_\alpha, \alpha\in {\cal K}_0({\cal A})$ which are subject to the following relations:
$$ K_\alpha K_\beta = K_{\alpha+\beta},\quad Z^{(m)}_A K_\alpha = (A|\alpha)^{(-1)^m}
K_\alpha Z^{(m)}_A,\leqno (3.1.1)$$
$$ Z^{(m)}_A Z^{(m)}_B = \langle B,A\rangle \sum_C g_{AB}^C Z^{(m)}_C, \leqno (3.1.2)$$
$$Z^{(m+1)}_A Z^{(m)}_B = \sum_{M,N} \gamma_{AB}^{MN} \langle \bar B-\bar M, \bar M -\bar N\rangle
 \cdot Z^{(m)}_M  Z^{m+1}_N K_{\bar B-\bar M}^{(-1)^m},
\leqno (3.1.3) $$
$$ Z^{(m)}_A  Z^{(n)}_B = (A|B)^{(-1)^{m-n}(n-m+1)} Z^{(n)}_B  Z^{(m)}_A , \quad |m-n|\geq 2.\leqno (3.1.4)$$

It is clear from these relations that the rule
$$\Sigma(Z^{(m)}_A) = Z^{(m+1)}_A, \quad \Sigma(K_\alpha)=K_\alpha^{-1}, \leqno (3.1.5)$$
defines an automorphism $\Sigma: L({\cal A})\to L({\cal A})$ which we call the shift
(or suspension) automorphism. 

\vskip .3cm

\noindent {\bf (3.2) Compatibility of the relations.} We now want to show that the relations
(3.1.1-4) are compatible in the sense that any element can be brought to a unique normal
form in which the upper indices of the $Z^{(m)}_A$ are increasing. 
More precisely, for any sequence $(a_i)_{i\in {\bf Z}}$ of elements of a possibly non-commutative
algebra $S$, almost all equal to 1, we define their ordered product to be
$$\prod_i^{\longrightarrow} a_i := a_pa_{p+1} ... a_q \in S,\leqno (3.2.1)$$
where $p,q$ are such thar $a_i=1$ unless $p\leq i\leq q$.

\proclaim (3.2.2) Proposition. The map of vector spaces
$$\nu: {\bf C}[{\cal K}_0({\cal A})]\otimes\bigotimes_{m\in {\bf Z}} R({\cal A}),
\quad K_\alpha\otimes\bigotimes_m [A_m] \mapsto \biggl( \prod_m^{\rightarrow} Z^{(m)}_{A_m}\biggr)
K_\alpha,$$
is an isomorphism.

Here and elsewhere in the paper all infinite tensor products of algebras are understood in the
restricted sense: almost all factors in any decomposable tensor are required to be 1. 
We will refer to an element explicitly realized as the value on $\nu$ on some tensor,
as being brought to the normal form.

\vskip .2cm

\noindent {\sl Proof:} The map $\nu$ is clearly surjective. To see the injectivity,
we need to establish the following two lemmas.

\proclaim (3.2.3) Lemma. For any $A,B,C\in {\cal A}$ and $m\in {\bf Z}$ the two possible ways of
bringing $Z^{(m+1)}_A Z^{(m)}_B Z^{(m-1)}_C$ to the normal form by using (3.1.3), lead to
the same answer.

\proclaim (3.2.4) Lemma. If $|m-n|\geq 2, |m+1-n|\geq 2$, then for any $A,B,C\in {\cal A}$
the multiplicative commutators with $Z^{(n)}_C$ of the left and the right hand sides of
(3.1.3) (the commutators being prescribed by (3.1.4)), are the same.

\vskip .2cm

\noindent {\sl Proof of Lemma 3.2.3:} The first way of bringing our monomial to the normal form is:
$$Z^{(m+1}_A Z^{(m)}_B Z^{(m-1)}_C = \sum_{M,N}
 \gamma_{AB}^{MN} \langle \bar B-\bar M, \bar M-\bar N\rangle Z^{(m)}_mZ^{(m+1)}_N K_{\bar B-\bar M}^{(-1)^{m+1}} Z^{(m-1)}_C = \leqno (3.2.5)$$
$$= \sum_{M,N}
 \gamma_{AB}^{MN} \langle \bar B-\bar M, \bar M-\bar N\rangle \cdot (\bar C| \bar B-\bar M)^{-1}
(\bar N, \bar C)^{-1} Z^{(m)}_M Z^{(m-1)}_C Z^{(m+1)}_N K_{\bar B-\bar M}^{(-1)^{m+1}}=$$
$$=\sum_{M,N,P,Q} \gamma_{AB}^{MN}\gamma_{MC}^{PQ} \langle \bar C-\bar P, \bar P-\bar Q\rangle
\cdot \langle \bar B-\bar M, \bar M-\bar N\rangle\cdot  (\bar C| \bar B-\bar M)^{-1}
(\bar N, \bar C)^{-1} (\bar C-\bar P|\bar N)\times$$
$$\times Z^{(m-1)}_P Z^{(m)}_Q Z^{(m+1)}_N K_{\bar C -\bar P}
^{(-1)^m} K_{\bar B-\bar M}^{(-1)^{m+1}},$$
while the other way is as follows:
$$Z^{(m+1)} Z^{(m)}_B Z^{(m-1)}_C = \sum_{P,U} \gamma_{PU}^{BC}
 \langle \bar C-\bar P, \bar P-\bar U\rangle \cdot Z^{(m+1)}_AZ^{(m-1)}_P
 Z^{(m)}_U K_{\bar C-\bar P}^{(-1)^m} =\leqno (3.2.6)$$
$$=\sum_{P.U} \sum_{P,U} \gamma_{PU}^{BC}
 \langle \bar C-\bar P, \bar P-\bar U\rangle (A|P)^{-1} Z^{(m-1)}_P Z^{(m+1)}_A Z^{(m)}_U 
K_{\bar C-\bar P}^{(-1)^m}=$$
$$=\sum_{P,U,Q,N} \gamma_{PU}^{BC}
\gamma_{AU}^{QN} \langle \bar C-\bar P, \bar P-\bar U\rangle (A|P)^{-1} \langle \bar U-\bar Q, \bar Q
-\bar N\rangle Z^{(m-1}_PZ^{(m)}_Q Z^{(m+1)}_N K_{\bar C-\bar P}^{(-1)^m}
 K_{\bar U-\bar Q}^{(-1)^{m+1}}.$$ 
Let us start to compare these two expressions. 
Our first remark is that the $\gamma$-quantities coincide.

\proclaim (3.2.7) Lemma. For any $A,B,C,P,Q,N\in {\cal A}$ we have the equality
$$\sum_M \gamma_{AB}^{MN}\gamma_{MC}^{PQ} = \sum_U \gamma_{BC}^{PU}\gamma_{AU}^{QN}.$$

\noindent {\sl Proof:} Both sides of the proposed equality have the following conceptual
meaning. They are equal to the orbifold number (modulo ${\rm Aut}(A)\times{\rm Aut}(B)
\times {\rm Aut}(C)$) of systems consisting, first, of a complex of length 2:
$$C\buildrel \psi\over\to B \buildrel \phi\over\to A, \quad \phi\psi=0,$$
and, second, of an identification of its cohomology, i.e., of isomorphisms
$$P\to {\rm Ker}(\psi), \quad Q\to {\rm Ker}(\phi)/{\rm Im}(\psi), \quad N\to
{\rm Coker}(\phi).$$
The object $M$ in the left hand side is ${\rm Im}(\psi)$, while $U$ in the right
hand side is ${\rm Im}(\phi)$. Q.E.D.

\vskip .2cm

Notice further that whenever a summand in any of the two sums in (3.2.7) is
non-zero, we have the equalities
$$\bar C-\bar P=\bar B-\bar U = \bar M-\bar Q, \quad \bar U-\bar Q = \bar A-\bar N =
\bar B -\bar M, \leqno (3.2.8)$$
which are obtained by applying  the fact that the Euler characteristic (in the Grothendieck group)
of a 4-term exact sequence is 0.

It follows from (3.2.8) that the ``$K$" factors in the end results of (3.2.5) and (3.2.6)
are the same. So it remains to compare the numerical factors given by the values of
the Euler form and its symmetrization.  We first compare the angle brackets, by noticing that
in virtue of (3.2.8),
$${ \langle \bar C-\bar P, \bar P-\bar U\rangle\cdot \langle \bar U-\bar Q, \bar Q-\bar N\rangle
\over \langle \bar C-\bar P, \bar P-\bar Q \rangle\cdot \langle \bar B-\bar M, \bar M-\bar N\rangle}
= { \langle \bar C-\bar P, \bar Q -\bar U\rangle\cdot \langle \bar U-\bar Q, \bar Q-\bar N\rangle
\over \langle \bar U-\bar Q, \bar M-\bar N\rangle}=$$
$$=\langle \bar C-\bar P, \bar Q-\bar U\rangle\cdot \langle \bar U-\bar Q, \bar Q-\bar M\rangle
=\langle \bar C-\bar P, \bar Q -\bar U\rangle \cdot \langle \bar U-\bar Q, \bar P-\bar C\rangle=$$
$$=  (\bar C-\bar P |\bar Q-\bar U)$$
and therefore the ratio of the total bracket contributions in the end results of (3.2.6)
and (3.2.5) is
$$ { (\bar C-\bar P|\bar Q-\bar U)\cdot (A|P)^{-1}\over (\bar C|\bar B-\bar M)^{-1}
 (\bar N|\bar C)^{-1} (\bar C-\bar P|\bar N)} = {(\bar C-\bar P|\bar Q-\bar U)\cdot 
(\bar C|\bar B-\bar M)\cdot
(P|N)\over (A|P)} =$$
$$(\bar C-\bar P|\bar Q-\bar U)\cdot (\bar C|\bar B-\bar M) \cdot (\bar P|\bar N-\bar A) =
(\bar C-\bar P|\bar Q-\bar U) \cdot (\bar C|\bar B-\bar M) \cdot (\bar P|\bar M-\bar B) =$$
$$= (\bar C-\bar P|\bar Q-\bar U) \cdot (\bar C-\bar P|\bar B-\bar M) = 1.$$
Lemma 3.2.3 is proved. 

\vskip .2cm

\noindent{\sl Proof of Lemma 3.2.4:} The product of $Z^{(m+1)}_A$ and $Z^{(m)}_B$ in any order
gives, when moved through $Z^{(n)}_C$, the factor
$$\bigl( (n-m)\bar A - (n-m+1)\bar B \bigl| \bar C\bigr) ^{(-1)^{m-n+1}}.$$
Thus it is enough to show that whenever $\gamma_{AB}^{MN}\neq 0$, one has the following
equality in ${\cal K}_0({\cal A})$:
$$(n-m)\bar A - (n-m+1)\bar B - (n-m)\bar N - (n-m+1)\bar M - (\bar B-\bar M),$$
where the last summand on the right comes from commuting $K_{\bar B-\bar M}^{(-1)^{m+1}}$
with $Z^{(n)}_C$. But, indeed, $\bar A-\bar B = \bar N-\bar M$, once $\gamma_{AB}^{MN}\neq 0$,
and the desired equality follows. 

This concludes the proof of Lemma 3.2.4 and Proposition 3.2.2. 

\vskip .3cm

\noindent {\bf (3.3) A basis in $L({\cal A})$.}
 It is  natural to label monomials in the $Z^{(i)}_A$
by (isomorphism classes of) graded objects of $\cal A$, 
i.e., by isomorphism classes of objects of $D^b({\cal A})$
which we represent as complexes with zero differential.

More precisely, to any graded object
$A^\bullet = \bigoplus A^{-i}[i]\in D^b({\cal A})$ we associate the monomial
$$Z(A^\bullet) = \prod_i^{\longrightarrow} {Z^{(i)}_{A^i} K_{A^i}^{(-1)^{i+1} i}
 \langle A^i, A^i\rangle ^i\over [A^i, A^{i-1}]}, \leqno (3.3.1)$$
where $[A^i, A^{i-1}]$ was defined in (2.2.1). The results of the previous subsection give: 

\proclaim (3.3.2) Proposition. The elements $Z(A^\bullet) K_\alpha $, for $A^\bullet
\in {\rm Ob}(D^b({\cal A}))/{\rm Iso}$ and $\alpha\in {\cal K}_0({\cal A})$, form a
{\bf C}-basis in $L({\cal A})$.

This, together with the homological interpretation of the quantities in Lemma (3.2.7),
suggests a deeper relation between $L({\cal A})$ and the derived category. More precisely, 
let $F: D^b({\cal A})\to D^b({\cal B})$ be any equivalence of triangulated categories
(we assume that both $\cal A$ and $\cal B$ satisfy our conditions of finiteness and homological dimension).
Then $F$ induces an isomorphism of Grothendieck groups $F_{\cal K}: {\cal K}_0({\cal A})\to
{\cal K}_0({\cal B})$ in a standard way. Our aim in this section is to prove the following result.
\vskip .1cm

\proclaim (3.4) Theorem. If $F$ is an equivalence of derived categories as above, then the
 correspondence
$$Z^{(p)}_A\to \Sigma^p(Z(F(A)), A\in {\cal A}, p\in {\bf Z}, \quad K_\alpha \to K_{F_{\cal K}(\alpha)},
\alpha\in {\cal K}_0({\cal A}),$$
defines an isomorphism of algebras $F_*: L({\cal A})\to L({\cal B})$.

\noindent {\sl Proof:} Our analysis is similar to that of (2.3). Namely, for $i\in {\bf Z}$
let ${\cal A}_i\i {\cal A}$ be the full subcategory of $A$ such that $F(A)\in {\cal B}[i]$. Then,
if $A_i\in {\cal A}_i, A_j\in {\cal A}_j$, we have
$${\rm Hom}_{\cal A}(A_i, A_j) = 0 \quad {\rm for}\quad j-i\notin \{0,1\}, \leqno (3.4.1a)$$
$${\rm Ext}^1_{\cal A}(A_i,A_j)=0 \quad {\rm for}\quad j-i\notin \{-1,0\}. \leqno (3.4.1b)$$
Also, each $A\in {\cal A}$ can be written uniquely as $A=\bigoplus A_i$ with $A_i\in {\cal A}_i$
and in the algebra $R({\cal A})$ we have the equality
$$[A] \quad =\quad \prod_i^{\longrightarrow} [A_i] \cdot |{\rm Hom}(A_i, A_{i-1})|^{-1/2}.
\leqno (3.4.2)$$
This means that the ordered product map defines an isomorphism of {\bf C}-vector spaces
$$\bigotimes_i R({\cal A}_i)\to R({\cal A}),$$
where the tensor product on the left is the restricted one. Denote by $R({\cal A}_i[-j])$
the subalgebra in $L({\cal A})$ spanned by $Z^{(j)}_A$ with $A\in {\cal A}_i$. 
Proposition 3.3.2 implies then the following.

\proclaim (3.4.3) Proposition. The map of vector spaces
$${\bf C}[{\cal K}_0({\cal A})] \otimes\bigotimes_{i,j\in {\bf Z}} R({\cal A}_i[-j])
\to L({\cal A}), \quad K_\alpha\otimes\bigotimes Z_{A_i}^{(j)} \mapsto 
\biggl(\prod_j^{\rightarrow}
\prod_i^{\rightarrow} Z_{A_i}^{(j)}\biggr)K_\alpha  ,$$
is an isomorphism.

 The proof of Theorem 3.4 consists basically of checking the relations and it is convenient to
first prove the following particular case, generalizing Theorem 2.3.

\proclaim (3.5) Proposition. For $F$ as above the rule $[A]\to Z(F(A))$ defines an injective
homomorphism of algebras $F_*: R({\cal A})\to L({\cal B})$. 

\noindent {\sl Proof:} 
 By an argument
similarly to Lemma 2.5.5, it is enough to prove the following partial statement.

\proclaim (3.5.1) Proposition. For any $i,j\in {\bf Z}$ and any $A'\in {\cal A}_{-i}$,
$A''\in {\cal A}_{-j}$, one has the equality $F_*([A']*[A'']) = F_*([A']) F_*([A''])$.

The minus sign is chosen for convenience, because $A\in {\cal A}_i$ means that $F(A)$, as a complex,
is situated in degree $(-i)$.
\vskip .1cm

\noindent {\sl Proof of (3.5.1):} We will work out several cases, similarly to the proof
of Theorem 2.3.

\vskip .2cm

\noindent\underbar{Case 1:} $i=j$. 
Let $A', A''\in {\cal A}_{-i}$ and let $F(A')=B'[-i], F(A'')=B''[-i]$.
If $A$ is such that $g_{A'A''}^A\neq 0$, then $A\in {\cal A}_{-i}$ and, denoting $B=F(A)$, we have
$g_{A'A''}^A=g_{B'B''}^B$. Thus

$$F_*([A'])F_*([A'']) = Z_{B'}^{(i)} K_{B'}^{(-1)^{i+1} i} Z_{B''}^{(i)}
K_{B''}^{(-1)^{i+1} i} \langle B', B'\rangle ^i \cdot \langle B'', B''\rangle ^i =$$
$$=Z^{(i)}_{B'}Z^{(i)}_{B''} K_{\bar B' +\bar B''}^{(-1)^{i+1} i} 
\langle B', B'\rangle ^i \cdot \langle B'', B''\rangle ^i \cdot (B'|B'')^i =$$
$$=\sum_{B\in {\cal B}} \langle B'', B'\rangle g_{B'B''}^B Z^{(i)}_{B} K_B^{(-1)^{i+1} i} \langle B, B\rangle ^i =
F_*([A']*[A'']),$$
where we used the equality $\bar B = \bar B' + \bar B''$ holding each time when $g_{B'B''}^B \neq 0$. 

\vskip .1cm

\noindent\underbar{Case 2:} $j=i+1$. Let $F(A')=B'[-i], F(A'')=B''[-i-1]$.
We have
$$F_*([A'])F_*([A'']) = Z_{B'}^{(i)} K_{B'}^{(-1)^{i+1} i} Z_{B''}^{(i+1)} K_{B''}^{ (-1)^{i+2} (i+1)}
\cdot \langle B', B'\rangle ^i \cdot \langle B'', B''\rangle ^{i+1},$$
while
$$F_*([A']*[A'']) = |{\rm Hom}(A'', A')|^{1/2} F([A'\oplus A'']) =$$ 
$$= { |{\rm Hom}(A'', A')|^{1/2} \cdot \langle B', B'\rangle ^i \cdot \langle B'', B''\rangle ^{i+1} \over
|{\rm Hom}(B'', B')|^{1/2}\cdot |{\rm Ext}^1 (B'', B')|^{1/2} }
 Z^{(i)}_{B'} K_{B'}^{(-1)^{i+1} i} Z^{(i+1)}_{B''}
K_{B''}^{(-1)^{i+2} (i+1)},$$
which is exactly the same once we recall that
$${\rm Ext}^1(B'', B') = {\rm Hom}(A'', A'), \quad {\rm Hom}(B'', B') = {\rm Ext}^{-1}(A'', A')=0.$$

\vskip .1cm

\noindent\underbar{Case 3:} $j=i-1$ and $i$ is even. 
 Let $F(A')=B'[-i], F(A'')=B''[-i+1]$. Then
$$F_*([A'])F_*([A'']) = Z^{(i)}_{B'} K_{B'}^{-i}Z^{(i-1)}_{B''} K_{B''}^{i-1} \langle B', B'\rangle ^i
\langle B'', B''\rangle ^{i-1}=$$
$$Z^{(i)}_{B'}Z^{(i-1)}_{B''} K_{-i\bar B' + (i-1) \bar B''}  \langle B', B'\rangle ^i
\langle B'', B''\rangle ^{i-1} (B'|B'')^{-i}=$$
$$= Z^{(i)}_{B'}Z^{(i-1)}_{B''} K_{-i\bar B' + (i-1) \bar B''} 
\langle \bar B'-\bar B'', \bar B'-\bar B''\rangle ^i \langle B'', B''\rangle ^{-1}=$$
$$\sum_{M,N} \gamma_{B', B''}^{MN} { \langle \bar B''-\bar M, \bar M\rangle \langle N, B''\rangle
\langle \bar B'-\bar B'', \bar B'-\bar B''\rangle ^i\over
\langle B'', B''\rangle } Z^{(i-1}_M K_{\bar B''-\bar M} Z^{(i)}_N K_{-i\bar B' + (i-1) \bar B''}.$$
We can transform the fraction in this expression in the same way as in the proof of Case (0,1)
of Proposition 2.5.3, getting:
$$ \langle B''[-i+1], B'[-i]\rangle \sum_{M,N} g_{B'[-i], B''[-i+1]}^{N[-i]\oplus M[-i+1]}
{Z^{(i-1)}_M K_M^{-1}Z^{(i)}_N K_{i (\bar B''-\bar B')} \langle \bar B''-\bar B', 
\bar B''-\bar B'\rangle ^i \over
\langle M, M\rangle\cdot \langle N,M\rangle \cdot |{\rm Ext}^1(N,M)|}.$$
Notice that for any non-zero summand we have $\bar B''-\bar B' =\bar M-\bar N$. Therefore
the last expression equals
$$\langle B''[-i+1], B'[-i]\rangle \sum_{M,N} g_{B'[-i], B''[-i+1]}^{N[-i]\oplus M[-i+1]}
{Z^{(i-1)}_M K_M^{i-1} Z^{(i)}_N K_N^{-i} \langle \bar N-\bar M, \bar N-\bar M\rangle ^i (N|M)^i
\over \langle M, M\rangle \cdot \langle N,M\rangle \cdot |{\rm Ext}^1(N,M)|}=$$
$$=\langle B''[-i+1], B'[-i]\rangle \sum_{M,N} g_{B'[-i], B''[-i+1]}^{N[-i]\oplus M[-i+1]}
{ Z^{(i-1)}_M K_M^{i-1} Z^{(i)}_N K_N^{-i} \langle M, M\rangle^{i-1} \langle N,N\rangle^i
\over [N,M] }=$$
$$= \langle A'', A'\rangle \sum_{A\in {\cal A}} g_{A'A''}^A F_*([A]) = F_*([A']*[A'']).$$

\vskip .2cm

\noindent\underbar{Case 4:} $j=i-1$ and $i$ is odd. Keeping the conventions for $A', A'', B', B''$
the same as before, we have$$F_*([A']) F_*([A'']) = Z^{(i)}_{B'} K_{B'}^i Z_{B''}^{(i-1)}
 K_{B''}^{-i+1} \langle B', B'\rangle^i \langle B'', B''\rangle^{i-1}=$$
$$=Z^{(i)}_{B'} Z^{(i-1)}_{B''} K_{i\bar B' - (i-1)\bar B''} \langle B', B'\rangle ^i 
\langle B'', B''\rangle^{i-1} (B'|B'')^{-i}=$$
$$= Z_{B'}^{(i)} Z_{B''}^{(i-1)} K_{i\bar B'-(i-1)\bar B''} \langle \bar B'-\bar B'', \bar B'-\bar B''
\rangle ^i \langle B'', B''\rangle ^{-1}=$$
$$=\sum_{M,N} \gamma_{B', B''}^{MN} { \langle \bar B''-\bar M, \bar M\rangle \cdot 
\langle N, B''\rangle \cdot \langle \bar B'-\bar B'', \bar B'-\bar B''\rangle ^i \over
\langle B'', B''\rangle}
Z^{(i-1)}_M K_{\bar M-\bar B''} Z^{(i)}_N K_{i\bar B'-(i-1)\bar B''}=$$
$$= \langle B''[-i+1], B'[-i]\rangle \sum_{M,N} g_{B'[-i], B''[-i+1]}^{N[-i]\oplus M[-i+1]}
{Z^{(i-1)}_M K_M Z^{(i)}_N K_{i(\bar B'-\bar B'')} \langle \bar B'-\bar B'', \bar B'-\bar B''\rangle ^i 
\over \langle M,M\rangle\cdot \langle N,M\rangle \cdot |{\rm Ext}^1(N,M)|} =$$
$$= \langle B''[-i+1], B'[-i]\rangle \sum_{M,N} g_{B'[-i], B''[-i+1]}^{N[-i]\oplus M[-i+1]}
{ Z^{(i-1)}_M K_M^{-i+1} Z_N^{(i)} K_N^i \langle \bar N-\bar M, \bar N-\bar M\rangle ^i
(N|M)^i \over \langle M,M\rangle [N,M] }=$$
$$= \langle B''[-i+1], B'[-i]\rangle \sum_{M,N} g_{B'[-i], B''[-i+1]}^{N[-i]\oplus M[-i+1]}
{ Z^{(i-1)}_M K_M^{-i+1} Z_N^{(i)} K_N^i \langle M,M\rangle^{i-1} \langle N,N\rangle^i
\over [N,M] },$$
and the argument is finished as in Case 3.

\vskip .2cm

\noindent\underbar{Case 5:} $|i-j|\geq 2$. Let $F(A')=B'[-i], F(A'')=B''[-j]$.
Then
$${\rm Hom}(A', A'') = {\rm Ext}^1(A', A'') = {\rm Hom}(A'', A') = {\rm Ext}^1(A'', A')=0.
\leqno (3.5.2)$$
Thus $[A']*[A''] = [A'\oplus A'']$. On the other hand, (3.5.2) implies that $(A'|A'') =
(B'|B'')=0$ and therefore
$$F_*([A'])F_*([A'']) = Z^{(i)}_{B'} K_{B'}^{(-1)^{i+1} i} Z_{B''}^{(j)} K_{B''}^{(-1)^{j+1} j}
\langle B',B'\rangle ^i \cdot \langle B'', B''\rangle ^j =$$
$$= F_*([A'\oplus A'']) =
F_*([A']*[A'']).$$
Proposition 3.5 is proved.

\vskip .3cm

\noindent {\bf (3.6) End of the proof of Theorem 3.4.} Let us define a linear operator
$F_*: L({\cal A})\to L({\cal B})$ by postulating its values on generators to be as stated in
the theorem and extending it to products of generators by using the normal form
of Proposition 3.2.2. In other words, we put
$$F_*\left( \biggl( \prod_i^{\to} Z^{(m)}_{A_m} \biggr) K_\alpha\right) = 
\biggl( \prod_i^{\to} F_*( Z^{(m)}_{A_m}) \biggr) K_{F_{\cal K}(\alpha)}.$$
 It is clear that
$F_*$ is bijective, so we need only to prove that it is an algebra homomorphism, i.e., that
it preserves the relations (3.1.1-4) in $L({\cal A})$. For (3.1.1) it is clear.
The preservation of (3.1.2) is the content of Proposition 3.5. The condition
that $F_*$ preserves (3.1.3-4) can be stated after some change of notation
as follows:

\proclaim (3.6.1) Proposition. Let $A', A''\in {\cal A}$. Then
$$F_*(Z^{(p+1)}_{A'})F_*(Z^{(p)}_{A''})=$$
$$\sum_{M,N\in {\cal A}}\gamma_{A' A''}^{MN} 
\langle \bar A''-\bar M, \bar M\rangle\cdot \langle \bar N, \bar A''-\bar M\rangle\cdot
F_*(Z^{(p)}_M)K^{(-1)^p}_{\overline{F(A'')}-\overline{F(M)}} F_*(Z^{(p+1)}_N).$$

\proclaim (3.6.2) Proposition. For $A', A''\in {\cal A}$ and $q\leq p-2$ we have
$$F_*(Z^{(p)}_A) F_*(Z^{(q)}_B) = (A|B)^{(-1)^{p-q}(q-p+1)} F_*(Z^{(q)}_B)F_*(Z^{(p)}_A).$$

Before proceeding to prove these statements, let us note the following.

\proclaim (3.6.3) Lemma. To establish Propositions 3.6.1-2 in full generality, it is enough
to prove them under the assumption that $A'\in {\cal A}_{-i}$, $A''\in {\cal A}_{-j}$
for some $i,j$.

\noindent {\sl
 Proof of the lemma:} Indeed, the multiplication law in $L({\cal A})$
with respect to the basis $K_\alpha Z(A^\bullet)$ consists of bringing the product
of two basis vectors to the normal form using (3.1.1-4). So our propositions, together
with what has been already proved, just say that
$$F_*(Z^{(p)}_{A'}) F_*(Z^{(q)}_{A''}) = F_*(Z^{(p)}_{A'}Z^{(q)}_{A''}).
\leqno (3.6.4)$$
If this is known each time when $A'\in {\cal A}_{-i}$, $A''\in {\cal A}_{-j}$,
then Proposition 3.4.3 together with an obvious modification of Lemma 2.5.5
give that (3.6.4) is valid in general. 

\vskip .3cm

\noindent {\bf (3.7) Proof of Proposition 3.6.1 when $A'\in {\cal A}_{-i}$, $A''\in {\cal A}_{-j}$.}
Again, we have to consider several cases.

\vskip .1cm

\noindent \underbar{Case 1:} $j=i$. Let $F(A')=B'[-i], F(A'')=B''[-i]$.
 In this case $F$ establishes a bijection between exact sequences
$$0\to M\to A''\to A'\to N\to 0$$
in $\cal A$ and exact sequences
$$0\to C\to B''\to B'\to D\to 0$$
in $\cal B$, so that for $C,D$ corresponding to $M,N$ we have
 $\gamma_{A' A''}^{MN}=\gamma_{B'B''}^{CD}$. From this the statement follows rather directly, by
comparing the normal form of $F_*(Z^{(p+1)}_{A'})F_*(Z^{(p)}_{A''})$ with the
image under $F_*$ of the normal form of $Z^{(p+1)}_{A'}Z^{(p)}_{A''}$.

\vskip .1cm

\noindent \underbar{Case 2:} $j=i+1$. Let $F(A') = B'[-i], F(A'')=B''[-i-1]$. Then
$$F_*(Z^{(p+1)}_{A'} = Z^{(i+p+1)}_{B'} K_{B'}^{(-1)^{i+p} i} \langle B', B'\rangle ^i,
\quad
F_*(Z^{(p)}_{A''} = Z^{(i+p+1)}_{B''} K_{B''}^{(-1)^{i+p} (i+1)} \langle B'', B''\rangle ^{i+1}.$$
Let $F^{-1}: D^b({\cal B})\to D^b({\cal A})$ be an inverse equivalence, and $(F^{-1})_*;
L({\cal B})\to L({\cal A})$ be the corresponding linear operators. The 
The argument needed is to handle our case is identical to the reasoning (already made in
the proof of Proposition 3.5) that $(F^{-1})_*$ realizes $R({\cal B})$ inside $L({\cal A})$. 

\vskip .1cm

\noindent \underbar{Case 3:} $j=i-1$. Let $F(A')=B'[-i], F(A'')=b]][-i+1]$. Then
${\rm Hom}(A'', A')=$ ${\rm Ext}^{-1}(B'', B')=0$, so $\gamma_{A'A''}^{MN}=0$ always except the case
$M=A'', N=A'$, in which case the value is 1. This means that
$$Z^{(p+1)}_{A'}Z^{(p)}_{A''} = Z^{(p)}_{A''} Z^{(p+1)}_{A'},$$
and we have to verify that $F_*$ preserved this commutativity. We have:
$$F_*(Z^{(p+1)}_{A'})= Z^{(p+i+1)}_{B'} K_{B'}^{(-1)^{p+i} i} \langle B', B'\rangle ^i, \quad
F_*(Z^{(p)}_{A''})= Z^{(p+i-1)}_{B''} K_{B''}^{(-1)^{p+i} (i-1)} \langle B'', B''\rangle ^{i-1},$$
and these expressions indeed commute because 
$$ Z^{(p+i+1)}_{B'}Z^{(p+i-1)}_{B''} = (B'|B'')^{-1} Z^{(p+i-1)}_{B''}Z^{(p+i+1)}_{B'}$$
and because of the commutation relation of the $K$'s with the $Z$'s.

\vskip .1cm

\noindent \underbar{Case 4:} $|j-i|\geq 2$. If $F(A')=B'[-i], F(A'')=B''[-j]$, then all the Hom and
Ext between $A'$ and $A''$ in either order are 0, so $(A'|A'')=(B'|B'')=1$,
and the argument is the same as in the previous case, only simpler because we do
not have to care about $(B'|B'')$. 

\vskip .3cm

\noindent {\bf (3.8) Proof of Proposition 3.6.2 when $A'\in {\cal A}_{-i}, A''\in {\cal A}_{-j}$.}
We consider several cases as to the relative position of $p+i$ and $q+j$. We denote
$F(A')=B'[-i], F(A'')=B''[-j]$.

\vskip .1cm

\noindent \underbar{Case 1:} $|(p+i)-(q+j)|\geq 2$. In this case
$$F_*(Z^{(p)}_{A'}) = Z^{(p+i)}_{B'} K_{B'}^{(-1)^{p+i+1} i} \langle B', B'\rangle ^i, \quad
F_*(Z^{(q)}_{A''}) = Z^{(q+j)}_{B''} K_{B'}^{(-1)^{q+j+1} j} \langle B'', B''\rangle ^j,$$
and by our assumption we have
$$F_*(Z^{(p)}_{A'})F_*(Z^{(q)}_{A''}) = (B'|B'')^\lambda F_*(Z^{(q)}_{A''})F_*(Z^{(p)}_{A'}),$$
where
$$\lambda = (-1)^{p+i-q-j} (q+j-p-1-1+i-j) = (-1)^{p+i-q-j}(q-p-1),$$
and once we take into account that $(B'|B'')=(A'|A'')^{(-1)^{i-j}}$,
we get the claimed statement.

\vskip .1cm

\noindent \underbar{Case 2:} $p+i=q+j$. This implies that $|i-j|\geq 2$ and thus there is neither Hom
nor Ext between $A', A''$ in either direction, hence $(A'|A'')=0$. So $Z^{(p)}_{A'}$
commutes with $Z^{(q)}_{A''}$. On the other hand, the vanishing of Hom and Ext implies that
$$Z^{(p+i)}_{B'} Z^{(p+i)}_{B''} = Z^{(p+i)}_{B''}Z^{(p+i)}_{B'} = Z^{(p+i)}_{B'\oplus B''},$$
so we are done.

\vskip .1cm

\noindent \underbar{Case 3:} $(p+i)-(q+j)=1$. Thus $i-j=q-p+1$. Since $|q-p|\geq 2$,
we will have $|i-j|\geq 2$ always except the case $q-p=-2$, when $i-j=-1$. In any event,
${\rm Hom}_{\cal B}(B'', B') = {\rm Ext}^{i-j}_{\cal A}(A'', A')=0$.
On the other hand,
$$F_*(Z^{(p)}_{A'}) = Z^{(p+i)}_{B'} K_{B'}^{(-1)^{p+i+1} i} \langle B', B'\rangle ^i, \quad
F_*(Z^{(q)}_{A''}) = Z^{(p+i-1)}_{B''} K_{B''}^{(-1)^{p+i} j} \langle B'', B''\rangle ^j,$$
and bringing their product to the normal form involves quantities $\gamma_{B'B''}^{MN}$
which vanish unless $B'=N, B''=M$. So $Z^{(p+i)}_{B'}$ and $Z^{(p+i-1)}_{B''}$
commute, and therefore
$$F_*(Z^{(p)}_{A'})F_*(Z^{(q)}_{A''})  = (B'|B'')^{-i+j} F_*(Z^{(q)}_{A''})F_*(Z^{(p)}_{A'}),$$
which is exactly what  we need, once we recall that $i-j=q-p+1$ and thus
$$(B'|B'')^{j-i}=(A'|A'')^{(-1)^{i-j} (j-i)}=
(A'|A'')^{(-1)^{q-p+1}(p-q-1)} =
(A'|A'')^{(-1)^{q-p}(q-p+1)}.$$

\vskip .1cm

\noindent \underbar{Case 4:} $(p+i)-(q+j)=-1$ is treated in an similar way. Theorem 3.4
is completely proved.

\vfill\eject

\centerline {\bf \S 4. Examples and discussion.}

\vskip 1cm

\noindent {\bf (4.1) Example: the ``universal cover" of the quantum group.}
Let us illustrate the construction of the lattice algebra on the classical example
of representations of quivers [R1-3].
 Let $\cal G$ be a semisimple
simply laced complex Lie algebra and $\Gamma$ be its Dynkin graph. Thus vertices of $\Gamma$
are identified with the  simple roots of $\cal G$ and for two such vertices $i\neq j$ the entry $a_{ij}$
of the Cartan matrix of $\cal G$ is minus the number of edges joining $i$ and $j$. Suppose
that an orientation of $\Gamma$ is chosen and let ${\cal A}={\rm Rep}_{{\bf F}_q}(\Gamma)$ be the
category of representations of $\Gamma$ over the finite field ${\bf F}_q$. Recall that such a representation $V$
is a rule which associates to any vertex $i$ a finite-dimensional ${\bf F}_q$-vector space
$V_i$ and to any edge $i\buildrel e\over\rightarrow j$ a linear operator
$V_e: V_i\to V_j$. As shown by Ringel ({\it loc. cit.}),
 the algebra $B({\cal A})$ is in this case isomorphic
to  $U_q({\bf b}^+)$, a natural ``Borel" subalgebra in $U_q({\cal G})$, the quantized enveloping algebra of $\cal G$.
More precisely, $U_q({\cal G})$ is generated by the symbols $E_i^\pm, K_i^{\pm 1}$ for $i\in {\rm Vert}(\Gamma)$
subject to the relations:
$$ E_i^\pm K_j = q^{\pm a_{ij}} K_j E_i^\pm, \quad K_iK_j=K_jK_i,\leqno (4.1.1)$$
$$\sum_{\nu=0}^{1-a_{ij}}
 {1-a_{ij}\choose \nu}_q (E_i^\pm)^\nu E_j^{\pm} (E_i^\pm)^{1-a_{ij}-\nu}=0, \quad i\neq j\leqno (4.1.2)$$
$$ [E^+_i, E^-_j] = {\delta_{ij} (K_i-K_i^{-1})\over q-1},\leqno (4.1.3)$$
and $B({\cal A})$ is isomorphic to the subalgebra generated by the $E^+_i$ and $K_i^{\pm 1}$.
Explicitly, $E_i^+$ corresponds to the element $[V(i)]$ where $V(i)$ is the respesentation associating ${\bf F}_q$
to the $i$th vertex and 0 to all other vertices. Similarly, $K_i$ corresponds to the element $K_{V(i)}$
of $B({\cal A})$. 

From this and the form of the comultiplication in $U_q({\bf b}^+)$ it is easy to deduce the following fact.

\proclaim (4.1.4) Proposition. For ${\cal A}={\rm Rep}_{{\bf F}_q}(\Gamma)$ the algebra $L({\cal A})$
is generated by the symbols $Z^{(m)}_i$, $m\in {\bf Z}, i\in {\rm Vert}(\Gamma)$ and $K_i^{\pm 1}$, 
$i\in {\rm Vert}(\Gamma)$ subject only to the following relations:
$$Z_i^{(m)}K_j = q^{(-1)^m a_{ij}}K_j Z_i^{(m)}, \quad K_iK_j = K_jK_i,\leqno (4.1.5)$$
$$\sum_{\nu=0}^{1-a_{ij}}
 {1-a_{ij}\choose \nu}_q (Z_i^{(m)})^\nu Z_j^{(m)} (Z_i^{(m)})^{1-a_{ij}-\nu}=0,
 \quad i\neq j \in {\rm Vert}(\Gamma), m\in {\bf Z},\leqno (4.1.6)$$
$$[Z^{(m)}_i, Z^{(m-1)}_j] = { \delta_{ij} K_i^{(-1)^m}\over q-1},\leqno (4.1.7)$$
$$Z^{(m)}_i  Z^{(n)}_j =  q^{(-1)^{m-n} (n-m+1)a_{ij}} Z^{(n)}_jZ^{(m)}_i \quad \forall i,j, \quad {\rm if}\quad |m-n|\geq 2.\leqno (4.1.8)$$

This shows that $L({\cal A})$ can be viewed as the ``universal cover", or the {\bf Z}-periodic version
of $U_q({\cal G})$. Note that the right hand side of the relation (4.1.7) is just one summand of the right hand
side of the similar relation (4.1.3): informally, in $L({\cal A})$ the other summand  is still present
but pertains to a different pair of generators. Let us also note the similarity of $L({\cal A})$
with so-called ``lattice Kac-Moody algebras" of [AFS].

\vskip .3cm

\noindent {\bf (4.2) ``Naive" lattice algebras.} The compatibility of the relations in
$L({\cal A})$ is not quite obvious, in particular because of the oscillator-type relations
(3.1.4) between copies of $R({\cal A})$ associated to non-adjacent lattice sites. 
So it may be useful to compare $L({\cal A})$ with the following general construction
which produces algebras in which such compatibility holds for free.

\vskip .1cm

Let $\Xi_m, m\in {\bf Z}$, be Hopf algebras and $\phi_m: \Xi_m\times \Xi_{m+1}\to {\bf C}$
be Hopf pairings. Define the {\it naive lattice algebra} $N=N(\{\Xi_m, \phi_m\})$
to be generated by elements of all the algebras $\Xi_m$ so that inside each $\Xi_m$
the elements are multiplied according to the multiplication law there while for
 elements of different
algebras we impose the relations:
$$\xi_{m+1}\xi_m = ({\rm Id}\otimes \phi_m \otimes {\rm Id})(\Delta_{\Xi_m}(\xi_m)
\otimes \Delta_{\Xi_{m+1}}(\xi_{m+1})), \leqno (4.2.1)$$
$$\xi_m\xi_{m'}=\xi_{m'}\xi_m, \quad |m-m'|\geq 2.$$
Thus if we put $\Xi_m$ at the $m$th site of a lattice, then the adjacent algebras
form a Heisenberg double while non-adjacent algebras commute. 

\proclaim (4.2.3) Proposition. The ordered product map $\bigotimes_{m\in {\bf Z}}\Xi_m\to N$
is always an isomorphism.

\noindent {\sl Proof:} We need only to verify that the two ways of bringing
any element   $\xi_{m+1}\xi_m\xi_{m-1}$ to the normal form $\sum_i \xi_{m-1}^{(i)}\xi_m^{(i)}
\xi_{m+1}^{(i)}$ by using (4.2.1), give the same result. But this easily follows
from the coassociativity of $\Delta_{\Xi_m}$.

\vskip .2cm

Applying this construction ot the case when for each $m$ we take $\Xi_m=B({\cal A})$ and $\phi_m
=\phi$ to be the Hopf pairing of (1.4), we get an algebra $N({\cal A})$ similar to $L({\cal A})$.
It is generated by symbols $Y^{(m)}_A, m\in {\bf Z}$, $A\in {\cal A}$ as well as $K^{(m)}_\alpha$,
$m\in{\bf Z}, \alpha\in {\cal K}_0({\cal A})$ with relations which are easy to find from
Proposition 1.5.3. In particular, for each $m$ the $K^{(m)}_\alpha$ form a copy of
${\bf C}[ {\cal K}_0({\cal A})]$, but adjacent copies do not commute. Because of this,
$N({\cal A})$ is not invariant under derived equivalence.

\vskip .3cm

\noindent {\bf (4.3) The bracket-free algebra $F({\cal A})$.} Most of the trouble in dealing with
the algebra $L({\cal A})$ comes from manipulating products of brackets, i.e., of values of the
Euler form. So it is tempting to define another, simpler algebra, by just dropping all these brackets.
More precisely, let $F({\cal A})$ be the algebra generated by synbols $X^{(m)}_A$,
$A\in {\cal A}, m\in {\bf Z}$ subject to the following relations:
$$X^{(m)}_A X^{(m)}_B = \sum_C g_{AB}^C X^{(m)}_C, \leqno (4.3.1)$$
$$X^{(m+1)}_A X^{(m)}_B =\sum_{M,N} \gamma_{AB}^{MN} X^{(m)}_M X^{(m+1)}_N, \leqno (4.3.2)$$
$$X^{(m)}_A X^{(n)}_B = X^{(n)}_BX^{(m)}_A, \quad |m-n|\geq 2.\leqno (4.3.3)$$
Note that these relations are compatible, as it follows from Lemma 3.2.7. In other words, 
the elements
$$X(A^\bullet) = \prod_m^{\to} X^{(m)}_{A^m}, \quad A^\bullet = \bigoplus_m A^m[-m], A^m\in {\cal A}$$
form a basis in $F({\cal A})$. The multiplication law in this basis can be described very nicely:
$$X(A^\bullet)X(B^\bullet) = \sum_{C^\bullet} \gamma_{A^\bullet B^\bullet}^{C^\bullet}
X(C^\bullet),\leqno (4.3.4)$$
where $\gamma_{A^\bullet B^\bullet}^{C^\bullet}$ is the orbifold number of long exact sequences
(2.4.3). Also, the procedure of bringing a maximally non-normal product to the normal form
can be nicely described in homological terms:
$$X^{(n)}_{A^{-n}} X^{(n-1)}_{A^{-n+1}} ... X^{(0)}_{A^0} = \sum_{d: A^\bullet\to A^\bullet[1]\atop
d^2=0} X(H^\bullet_d(A^\bullet)) \prod_m { |{\rm Aut}(H^m_d(A^\bullet))|\over
|{\rm Aut}(A^m)|}, \leqno (4.3.5)$$
where the sum is over all differentials making $A^\bullet$ into a complex, and 
$H^m_d(A^\bullet)$ is the $m$th cohomology with respect to $d$.

\vskip .1cm

Analogs of (4.3.4-5) can be obtained for the algebra $L({\cal A})$ as well, but they will be
encumbered by a lot of extra factors. However, these factors seem necessary to ensure the
invariance of the algebra under derived equivalence. In fact, the reason why $L({\cal A})$
possesses such invariance, is a subtle matching of two discrepancies. First is the discrepancy between 
the number
$g_{A^\bullet B^\bullet}^{C^\bullet}$
 of exact triangles in $D^{[-1,0]}({\cal A})$ and the number $\gamma_{A^\bullet B^\bullet}^{C^\bullet}$
of corresponding long exact sequences, which (for a particular case) 
was determined in Proposition 2.4.3 to be the factor
$|{\rm Ext}^1(N,M)|$.  The second is the discrepancy between the Hall multiplication 
$\circ$ in $H({\cal A})$ and its modification $*$ obtained by multiplying with $\langle B,A\rangle$,
see (1.3.4). When one term of a short exact sequence in $\cal A$ becomes shifted by 1
under a derived equivalence (so that we get an exact triangle in $D^{[-1,0]}({\cal A})$),
the difference between $*$ and $\circ$ will correspond
to the difference between $g$ and $\gamma$.

\vfill\eject

\centerline{\bf References.}

\vskip 1cm

\noindent [AF] A. Yu. Alexeev, L.D. Faddeev, Quantum $T^*G$ as 
a toy model for conformal field theory, {\it Comm. math. Phys.}, 
{\bf  141} (1991), 413-443.

\vskip .3cm

\noindent [AFS]  A. Yu. Alexeev, L.D. Faddeev, M. A. Semenov-Tian-Shansky, Hidden quantum
groups inside Kac-Moody algebras, Lecture Notes in Math. {\bf 1510}, p. 148-158,
Springer-Verlag 1992.

\vskip .3cm

\noindent [CP] V. Chari, A. Pressley, A Guide to Quantum groups,
Cambridge Univ. Press, 1995. 

\vskip .3cm

\noindent [GM] S.I. Gelfand, Y.I. Manin, Methods of Homological Algebra, Springer-Verlag, 1996.

\vskip .3cm

\noindent [Gr] J.A. Green, Hall algebras, hereditary algebras 
and quantum groups,
{\it Invent. Math.} {\bf  120} (1995), 361-377. 

\vskip .3cm

\noindent [Ha] D. Happel, Triangulated Categories in Representation Theory of Finite-Dimensional
Algebras (London Math. Soc. Lect. Note Series {\bf 119}), Cambridge Univ. Press, 1988.

\vskip .3cm

\noindent [J] A. Joseph, Quantum Groups and Their Primitive Ideals
(Ergebnisse der math. {\bf 29}), Springer-Verlag 1995.

\vskip .3cm

\noindent [Kap] M. Kapranov, Eisenstein series and quantum affine algebras, preprint alg-geom/9604018.

\vskip .3cm

\noindent [Kas] R.M. Kashaev, Heisenberg double and the pentagon relation, 
preprint q-alg/9503005.

\vskip .3cm

\noindent [Lu1] G. Lusztig, Introduction to Quantum Groups
(Progress in Math. 
{\bf  110}), Birkhauser, Boston, 1993. 

\vskip .3cm

\noindent [Lu2] G. Lusztig, Canonical bases arising from quantized enveloping 
algebras, {\it J. AMS}, {\bf  3} (1990), 447-498.

\vskip .3cm

\noindent [Lu3] G. Lusztig, Quivers, perverse sheaves and quantized 
enveloping algebras, {\it J. AMS}, {\bf  4} (1991), 365-421.

\vskip .3cm

\noindent [R1] C.M. Ringel, Hall algebras and quantum groups, 
{\it Invent. Math.}
{\bf  101} (1990), 583-592.

\vskip .3cm

\noindent [R2] C.M. Ringel, The composition algebra of a cyclic quiver,
{\it Proc. London Math. Soc.} (3) {\bf  66} (1993), 507-537.

\vskip .3cm

\noindent [R3] C.M. Ringel, Hall algebras revisited, in:
 ``Israel Math. Conf. Proc." vol. 7 (1993), p. 171-176. 

\vskip .3cm

\noindent [ST] M.A. Semenov-Tian-Shansky, Poisson Lie groups,
quantum duality principle and the quantum double,
{\it Contemporary math.}, {\bf  178}, p. 219-248,
Amer. Math. Soc, 1994.

\vskip .3cm

\noindent [X1] J. Xiao, Hall algebra in a root category, Preprint 95-070, Univ. of Bielefeld,
1995.

\vskip .3cm

\noindent [X2] J. Xiao, Drinfeld double and Green-Ringel theory of Hall algebras, Preprint
95-071, Univ. of Bielefeld, 1995.

\vskip 2cm

\noindent {\sl Author's address: Department of mathematics, Northwestern University, Evanston IL 
60208 USA \hfill\break
email: kapranov@math.nwu.edu}

\bye